# N-ALGEBRAIC STRUCTURES AND S-N-ALGEBRAIC STRUCTURES

**W. B. Vasantha Kandasamy**

**Florentin Smarandache**

**2005**

# N-ALGEBRAIC STRUCTURES AND S-N-ALGEBRAIC STRUCTURES

**W. B. Vasantha Kandasamy**
e-mail: vasantha@iitm.ac.in
web: http://mat.iitm.ac.in/~wbv

**Florentin Smarandache**
e-mail: smarand@unm.edu

**2005**

# CONTENTS

# PREFACE

In this book, for the first time we introduce the notions of N-groups, N-semigroups, N-loops and N-groupoids. We also define a mixed N-algebraic structure. We expect the reader to be well versed in group theory and have at least basic knowledge about Smarandache groupoids, Smarandache loops, Smarandache semigroups and bialgebraic structures and Smarandache bialgebraic structures.

The book is organized into six chapters. The first chapter gives the basic notions of S-semigroups, S-groupoids and S-loops thereby making the book self-contained. Chapter two introduces N-groups and their Smarandache analogues. In chapter three, N-loops and Smarandache N-loops are introduced and analyzed. Chapter four defines N-groupoids and S-N-groupoids. Since the N-semigroup structures are sandwiched between groups and groupoids, the study can be carried out without any difficulty. Mixed N-algebraic structures and S-mixed algebraic structures are given in chapter five. Some problems are suggested in chapter six.

It is pertinent to mention that several exercises and problems (Some in the form of proof to the theorems are given in all the chapters.) A reader who attempts to solve them will certainly gain a sound knowledge about these concepts. We have given 50 problems for the reader to solve in chapter 6.

The main aim of this book is to introduce new concepts and explain them with examples there by encouraging young mathematics to pursue research in this direction. Several theorems based on the definition can be easily proved with simple modification. Innovative readers can take up that job.

Also these notions find their applications in automaton theory and colouring problems. The N-semigroups and N-automaton can be applied to construct finite machines, which can perform multitasks, so their capability would be much higher than the usual automaton of finite machines constructed. We have suggested a list of references for further reading.

W.B.VASANTHA KANDASAMY
FLORENTIN SMARANDACHE
January 2005

## Chapter One

# INTRODUCTORY CONCEPTS

In this chapter the authors introduce the basic properties to make the book a self contained one. This chapter has 3 sections. In section one some properties of groups and S-semigroups are given. In section two basic notions about loops and S-loops are given. Third section gives some properties essential about groupoids and S-groupoids. The reader is expected to know basic properties about groups. The reader is expected to refer [35-40] for more about the concepts given in this chapter.

## 1.1 Groups, Smarandache semigroups and its basic properties

It is a well-known fact that groups are the only algebraic structures with a single binary operation that is mathematically so perfect that an introduction of a richer structure within it is impossible. Now we proceed on to define a group.

**DEFINITION 1.1.1:** *A non empty set of elements G is said to form a group if in G there is defined a binary operation, called the product and denoted by '•' such that*

i. $a, b \in G$ *implies that* $a \bullet b \in G$ *(closed)*
ii. $a, b, c \in G$ *implies* $a \bullet (b \bullet c) = (a \bullet b) \bullet c$ *(associative law)*
iii. *There exists an element* $e \in G$ *such that* $a \bullet e = e \bullet a = a$ *for all* $a \in G$ *(the existence of identity element in G).*

iv. *For every $a \in G$ there exists an element $a^{-1} \in G$ such that $a \bullet a^{-1} = a^{-1} \bullet a = e$ (the existence of inverse in G).*

**DEFINITION 1.1.2:** *A subgroup N of a group G is said to be a normal subgroup of G if for every $g \in G$ and $n \in N$, $g\, n\, g^{-1} \in N$.*

Equivalently by $gNg^{-1}$ we mean the set of all $gng^{-1}$, $n \in N$ then N is a normal subgroup of G if and only if $gNg^{-1} \subset N$ for every $g \in G$.

**THEOREM 1.1.1:** N is a normal subgroup of G if and only if $gNg^{-1} = N$ for every $g \in G$.

**DEFINITION 1.1.3:** *Let G be a group. $Z(G) = \{x \in G \mid gx = xg$ for all $g \in G\}$. Then Z(G) is called the center of the group G.*

**DEFINITION 1.1.4:** *Let G be a group, A, B be subgroups of G. If $x, y \in G$ define $x \sim y$ if $y = axb$ for some $a \in A$ and $b \in B$. We call the set $AxB = \{axb / a \in A, b \in B\}$ a double coset of A, B in G.*

**DEFINITION 1.1.5:**: Let G be a group. A and B subgroups of G, we say A and B are conjugate with each other if for some $g \in G$, $A = gBg^{-1}$.

Clearly if A and B are conjugate subgroups of G then o(A) = o(B).

**THEOREM: (LAGRANGE).** If G is a finite group and H is a subgroup of G then o(H) is a divisor of o(G).

**COROLLARY 1.1.1:** *If G is a finite group and $a \in G$, then o(a) / o(G).*

**COROLLARY 1.1.2:** *If G is a finite group and $a \in G$, then $a^{o(G)} = e$.*

In this section we give the two Cauchy's theorems one for abelian groups and the other for non-abelian groups. The main result on finite groups is that if the order of the group is n (n < $\propto$) if p is a prime dividing n by Cauchy's theorem we will always be able to pick up an element $a \in G$ such that $a^p = e$. In fact we can say Sylow's theorem is a partial extension of Cauchy's theorem for he says this finite group G has a subgroup of order $p^\alpha$($\alpha \geq 1$, p, a prime).

**THEOREM: (CAUCHY'S THEOREM FOR ABELIAN GROUPS).** *Suppose G is a finite abelian group and p / o(G), where p is a prime number. Then there is an element $a \neq e \in G$ such that $a^p = e$.*

**THEOREM: (CAUCHY):** *If p is a prime number and p | o(G), then G has an element of order p.*

Though one may marvel at the number of groups of varying types carrying many different properties, except for Cayley's we would not have seen them to be imbedded in the class of groups this was done by Cayley's in his famous theorem. Smarandache semigroups also has a beautiful analog for Cayley's theorem which is given by A(S) we mean the set of all one to one maps of the set S into itself. Clearly A(S) is a group having n! elements if $o(S) = n < \propto$, if S is an infinite set, A(S) has infinitely many elements.

**THEOREM: (CAYLEY)** *Every group is isomorphic to a subgroup of A(S) for some appropriate S.*

The Norwegian mathematician Peter Ludvig Mejdell Sylow was the contributor of Sylow's theorems. Sylow's theorems serve double purpose. One hand they form partial answers to the converse of Lagrange's theorem and on the other hand they are the complete extension of Cauchy's Theorem. Thus Sylow's work interlinks the works of two great mathematicians Lagrange and Cauchy. The following theorem is one, which makes use of Cauchy's theorem. It gives a nice partial converse to Lagrange's theorem and is easily understood.

**THEOREM: (SYLOW'S THEOREM FOR ABELIAN GROUPS)** *If G is an abelian group of order o(G), and if p is a prime number, such that $p^{\alpha} \mid o(G)$, $p^{\alpha+1} \nmid o(G)$, then G has a subgroup of order $p^{\alpha}$.*

**COROLLARY 1.1.3:** *If G is an abelian group of finite order and $p^{\alpha} \mid o(G)$, $p^{\alpha+1} \nmid o(G)$, then there is a unique subgroup of G of order $p^{\alpha}$.*

**DEFINITION 1.1.6:** *Let G be a finite group. A subgroup G of order $p^{\alpha}$, where $p^{\alpha} / o(G)$ but $p^{\alpha} \nmid o(G)$, is called a p-Sylow subgroup of G. Thus we see that for any finite group G if p is any prime which divides o(G); then G has a p-Sylow subgroup.*

**THEOREM (FIRST PART OF SYLOW'S THEOREM):** *If p is a prime number and $p^{\alpha}/ o(G)$ and $p^{\alpha+1} \nmid o(G)$, G is a finite group, then G has a subgroup of order $p^{\alpha}$.*

**THEOREM: (SECOND PART OF SYLOW'S THEOREM):** *If G is a finite group, p a prime and $p^{n} \mid o(G)$ but $p^{n+1} \nmid o(G)$, then any two subgroup of G of order $p^{n}$ are conjugate.*

**THEOREM: (THIRD PART OF SYLOW'S THEOREM):** *The number of p-Sylow subgroups in G, for a given prime, is of the form 1 + kp.*

Here we first recall the definition of Smarandache semigroups as given by Raul (1998) and introduce in this section concepts like Smarandache commutative semigroup, Smarandache weakly commutative semigroup, Smarandache cyclic and weakly cyclic semigroups.

**DEFINITION 1.1.7:** *The Smarandache semigroup (S-semigroup) is defined to be a semigroup A such that a proper subset of A is a group (with respect to the same induced operation).*

**DEFINITION 1.1.8:** *Let S be a S-semigroup. If every proper subset of A in S, which is a group is commutative then we say*

*the S-semigroup S to be a Smarandache commutative semigroup.*

Remark: It is important to note that if S is a commutative semigroup and if S is a S-semigroup then S is a Smarandache commutative semigroup. Here we are interested in finding whether there exists proper subsets of S-semigroups which are subgroups of which some of them are commutative and some non-commutative. This leads us to define:

**DEFINITION 1.1.9:** *Let S be S-semigroup, if S contains at least a proper subset A that is a commutative subgroup under the operations of S then we say S is a Smarandache weakly commutative semigroup.*

**DEFINITION 1.1.10:** *Let S be S-semigroup if every proper subset A of S which is a subgroup is cyclic then we say S is a Smarandache cyclic semigroup.*

**DEFINITION 1.1.11:** *Let S be a S-semigroup if there exists at least a proper subset A of S, which is a cyclic subgroup under the operations of S then we say S is a Smarandache weakly cyclic semigroup.*

**DEFINITION 1.1.12:** *Let S be a S-semigroup. If the number of distinct elements in S is finite, we say S is a finite S-semigroup otherwise we say S is a infinite S-semigroup.*

**DEFINITION 1.1.13:** *Let S be a S-semigroup. If A be a proper subset of S which is subsemigroup of S and A contains the largest group of S then we say A to be the Smarandache hyper subsemigroup of S.*

***Example 1.1.1:*** Let S(8) be the Smarandache symmetric semigroup of all mappings of the set S = (1, 2, 3, ... , 8). Now S(8) has a Smarandache hyper subsemigroup for take

$$S_8 \bigcup \left\{ \begin{pmatrix} 1 & 2 & 3 & \ldots & 8 \\ 1 & 1 & 1 & \ldots & 1 \end{pmatrix}, \right.$$

$$\left. \begin{pmatrix} 1 & 2 & 3 & \ldots & 8 \\ 2 & 2 & 2 & \ldots & 2 \end{pmatrix}, \begin{pmatrix} 1 & 2 & 3 & \ldots & 8 \\ 3 & 3 & 3 & \ldots & 3 \end{pmatrix}, \ldots, \begin{pmatrix} 1 & 2 & 3 & \ldots & 8 \\ 8 & 8 & 8 & \ldots & 8 \end{pmatrix} \right\}.$$

Clearly A is a subsemigroup of S(8) and has the largest group $S_8$ in it.

**THEOREM 1.1.2:** *Let S be a S-semigroup. Every Smarandache hyper subsemigroup is a Smarandache subsemigroup but every Smarandache subsemigroup is not a Smarandache hyper subsemigroup.*

**DEFINITION 1.1.14:** *Let S be a S-semigroup. We say S is a Smarandache simple semigroup if S has no proper subsemigroup A, which contains the largest subgroup of S or equivalently S has no Smarandache hyper subsemigroup.*

**THEOREM 1.1.3:** *$Z_p = \{0, 1, 2, \ldots , p – 1\}$ where p is a prime is a S-semigroup under multiplication modulo p. But $Z_p$ is a Smarandache simple semigroup.*

**DEFINITION 1.1.15:** *Let A be a S-semigroup. H is a proper subset of A ($H \subset A$) be a group under the operations of A. For any $a \in A$ define the Smarandache right coset $Ha = \{ha / h \in H\}$, Ha is called the Smarandache right coset of H in A.*

**DEFINITION 1.1.16:** *Let S be a S-semigroup. $H \subset S$ be a subgroup of S. We say aH is the Smarandache coset of H in S for $a \in S$ if $Ha = aH$, that is $\{ha / h \in H\} = \{ah / h \in H\}$.*

To prove the classical Cayley's theorem of group for S-semigroups we need the concept of S-semigroup homomorphism using which we will prove our result.

**DEFINITION 1.1.17:** *Let S and S' be any two S-semigroups. A map $\phi$ from S to S' is said to be a S-semigroup homomorphism if*

*$\phi$ restricted to a subgroup $A \subset S \to A' \subset S'$ is a group homomorphism. The S-semigroup homomorphism is an isomorphism if $\phi$: $A \to A'$ is one to one and onto. Similarly, one can define S-semigroup automorphism on S.*

**THEOREM: (CAYLEY'S THEOREM FOR S-SEMIGROUP)** *Every S-semigroup is isomorphic to a S-semigroup S(N); of mappings of a set N to itself, for some appropriate set N.*

For more about S-semigroups please refer [40].

## 1.2 Loops Smarandache Loops and their basic properties

We at this juncture like to express that books solely on loops are meagre or absent as, R.H.Bruck deals with loops on his book "*A Survey of Binary Systems*", that too published as early as 1958, [3]. Other two books are on "*Quasigroups and Loops*" one by H.O. Pflugfelder, 1990 which is introductory and the other book co-edited by Orin Chein, H.O. Pflugfelder and J.D. Smith in 1990 [10]. So we felt it important to recall almost all the properties and definitions related with loops [3]. As this book is on Smarandache loops so, while studying the Smarandache analogous of the properties the reader may not be running short of notions and concepts about loops.

**DEFINITION 1.2.1:** *A non-empty set L is said to form a loop, if on L is defined a binary operation called the product denoted by '$\bullet$' such that*

- *i. For all a, b $\in$ L we have $a \bullet b \in L$ (closure property).*
- *ii. There exists an element e $\in$ L such that $a \bullet e = e \bullet a = a$ for all a $\in$ L (e is called the identity element of L).*
- *iii. For every ordered pair (a, b) $\in L \times L$ there exists a unique pair (x, y) in L such that ax = b and ya = b.*

**DEFINITION 1.2.2:** *Let L be a loop. A non-empty subset H of L is called a subloop of L if H itself is a loop under the operation of L.*

**DEFINITION 1.2.3:** *Let L be a loop. A subloop H of L is said to be a normal subloop of L, if*

i. $xH = Hx$.
ii. $(Hx)y = H(xy)$.
iii. $y(xH) = (yx)H$

*for all* $x, y \in L$.

**DEFINITION 1.2.4:** *A loop L is said to be a simple loop if it does not contain any non-trivial normal subloop.*

**DEFINITION 1.2.5:** *The commutator subloop of a loop L denoted by L' is the subloop generated by all of its commutators, that is,* $\langle\{x \in L \,/\, x = (y, z) \text{ for some } y, z \in L\}\rangle$ *where for* $A \subseteq L$, $\langle A\rangle$ *denotes the subloop generated by A.*

**DEFINITION 1.2.6:** *If x, y and z are elements of a loop L an associator (x, y, z) is defined by,* $(xy)z = (x(yz))\,(x, y, z)$.

**DEFINITION 1.2.7:** *The associator subloop of a loop L (denoted by A(L)) is the subloop generated by all of its associators, that is* $\langle\{x \in L \,/\, x = (a, b, c) \text{ for some } a, b, c \in L\}\rangle$.

**DEFINITION 1.2.8:** *The centre Z(L) of a loop L is the intersection of the nucleus and the Moufang centre, that is* $Z(L) = C(L) \cap N(L)$.

**DEFINITION [26]**: *A normal subloop of a loop L is any subloop of L which is the kernel of some homomorphism from L to a loop.*

Further Pflugfelder [10] has proved the central subgroup Z(L) of a loop L is normal in L.

**DEFINITION [26]**: *Let L be a loop. The centrally derived subloop (or normal commutator- associator subloop) of L is*

*defined to be the smallest normal subloop $L' \subset L$ such that $L / L'$ is an abelian group.*

*Similarly nuclearly derived subloop (or normal associator subloop) of L is defined to be the smallest normal subloop $L_1 \subset L$ such that $L / L_1$ is a group. Bruck proves L' and $L_1$ are well defined.*

**DEFINITION [26]**: *The Frattini subloop $\phi(L)$ of a loop L is defined to be the set of all non-generators of L, that is the set of all $x \in L$ such that for any subset S of L, $L = \langle x, S\rangle$ implies $L = \langle S\rangle$. Bruck has proved as stated by Tim Hsu $\phi(L) \subset L$ and $L / \phi(L)$ is isomorphic to a subgroup of the direct product of groups of prime order.*

**DEFINITION [20]**: *Let L be a loop. The commutant of L is the set $(L) = \{a \in L / ax = xa \ \forall x \in L\}$. The centre of L is the set of all $a \in C(L)$ such that $a \bullet xy = ax \bullet y = x \bullet ay = xa \bullet y$ and $xy \bullet a = x \bullet ya$ for all $x, y \in L$. The centre is a normal subloop. The commutant is also known as Moufang Centre in literature.*

**DEFINITION [21]:** *A left loop $(B, \bullet)$ is a set B together with a binary operation '$\bullet$' such that (i) for each $a \in B$, the mapping $x \to a \bullet x$ is a bijection and (ii) there exists a two sided identity $1 \in B$ satisfying $1 \bullet x = x \bullet 1 = x$ for every $x \in B$. A right loop is defined similarly. A loop is both a right loop and a left loop.*

**DEFINITION [11]** : *A loop L is said to have the weak Lagrange property if, for each subloop K of L, $|K|$ divides $|L|$. It has the strong Lagrange property if every subloop K of L has the weak Lagrange property.*

**DEFINITION 1.2.9:** *A loop L is said to be power-associative in the sense that every element of L generates an abelian group.*

**DEFINITION 1.2.10:** *A loop L is diassociative loop if every pair of elements of L generates a subgroup.*

**DEFINITION 1.2.11:** *A loop L is said to be a Moufang loop if it satisfies any one of the following identities:*

i. *(xy) (zx) = (x(yz))x*
ii. *((xy)z)y = x(y(zy))*
iii. *x(y(xz) = ((xy)x)z*

*for all x, y, z* $\in$ *L.*

**DEFINITION 1.2.12:** *Let L be a loop, L is called a Bruck loop if x(yx)z = x(y(xz)) and* $(xy)^{-1} = x^{-1}y^{-1}$ *for all x, y, z* $\in$ *L.*

**DEFINITION 1.2.13:** *A loop (L, •) is called a Bol loop if ((xy)z)y = x((yz)y) for all x, y, z* $\in$ *L.*

**DEFINITION 1.2.14:** *A loop L is said to be right alternative if (xy)y = x(yy) for all x, y* $\in$ *L and L is left alternative if (xx)y = x(xy) for all x, y* $\in$ *L. L is said to be an alternative loop if it is both a right and left alternative loop.*

**DEFINITION 1.2.15:** *A loop (L, •) is called a weak inverse property loop (WIP-loop) if (xy)z = e imply x(yz) = e for all x, y, z* $\in$ *L.*

**DEFINITION 1.2.16:** *A loop L is said to be semi alternative if (x, y, z) = (y, z, x) for all x, y, z* $\in$ *L, where (x, y, z) denotes the associator of elements x, y, z* $\in$ *L.*

**THEOREM (MOUFANG'S THEOREM):** *Every Moufang loop G is diassociative more generally, if a, b, c are elements in a Moufang loop G such that (ab)c = a(bc) then a, b, c generate an associative loop.*

The proof is left for the reader; for assistance refer Bruck R.H. [3].

**DEFINITION 1.2.17:** *Let L be a loop, L is said to be a two unique product loop (t.u.p) if given any two non-empty finite subsets A and B of L with |A| + |B| > 2 there exist at least two*

*distinct elements x and y of L that have unique representation in the from x = ab and y = cd with a, c $\in$ A and b, d $\in$ B.*

*A loop L is called a unique product (u.p) loop if, given A and B two non-empty finite subsets of L, then there always exists at least one x $\in$ L which has a unique representation in the from x = ab, with a $\in$ A and b $\in$ B.*

**DEFINITION 1.2.18:** *Let (L, $\bullet$) be a loop. The principal isotope (L, $*$) of (L, $\bullet$) with respect to any predetermined a, b $\in$ L is defined by x $*$ y = XY, for all x, y $\in$ L, where Xa = x and bY = y for some X, Y $\in$ L.*

**DEFINITION 1.2.19:** *Let L be a loop, L is said to be a G-loop if it is isomorphic to all of its principal isotopes.*

The main objective of this section is the introduction of a new class of loops with a natural simple operation. As to introduce loops several functions or maps are defined satisfying some desired conditions we felt that it would be nice if we can get a natural class of loops built using integers.

Here we define the new class of loops of any even order, they are distinctly different from the loops constructed by other researchers. Here we enumerate several of the properties enjoyed by these loops.

**DEFINITION [32]:** *Let $L_n(m)$ = {e, 1, 2, …, n} be a set where n > 3, n is odd and m is a positive integer such that (m, n) = 1 and (m –1, n) = 1 with m < n.*

*Define on $L_n(m)$ a binary operation '$\bullet$' as follows:*

i. *$e \bullet i = i \bullet e = i$ for all $i \in L_n(m)$*
ii. *$i^2 = i \bullet i = e$ for all $i \in L_n(m)$*
iii. *$i \bullet j = t$ where $t = (mj – (m-1)i) \pmod n$*

*for all i, j $\in L_n(m)$; i $\neq$ j, i $\neq$ e and j $\neq$ e, then $L_n(m)$ is a loop under the binary operation '$\bullet$'.*

***Example 1.2.1:*** Consider the loop $L_5(2)$ = {e, 1, 2, 3, 4, 5}. The composition table for $L_5(2)$ is given below:

| • | e | 1 | 2 | 3 | 4 | 5 |
|---|---|---|---|---|---|---|
| e | e | 1 | 2 | 3 | 4 | 5 |
| 1 | 1 | e | 3 | 5 | 2 | 4 |
| 2 | 2 | 5 | e | 4 | 1 | 3 |
| 3 | 3 | 4 | 1 | e | 5 | 2 |
| 4 | 4 | 3 | 5 | 2 | e | 1 |
| 5 | 5 | 2 | 4 | 1 | 3 | e |

This loop is of order 6 which is both non-associative and non-commutative.

*Physical interpretation of the operation in the loop $L_n(m)$*:

We give a physical interpretation of this class of loops as follows: Let $L_n(m)$= {e, 1, 2, … , n} be a loop in this identity element of the loop are equidistantly placed on a circle with e as its centre.

We assume the elements to move always in the clockwise direction.

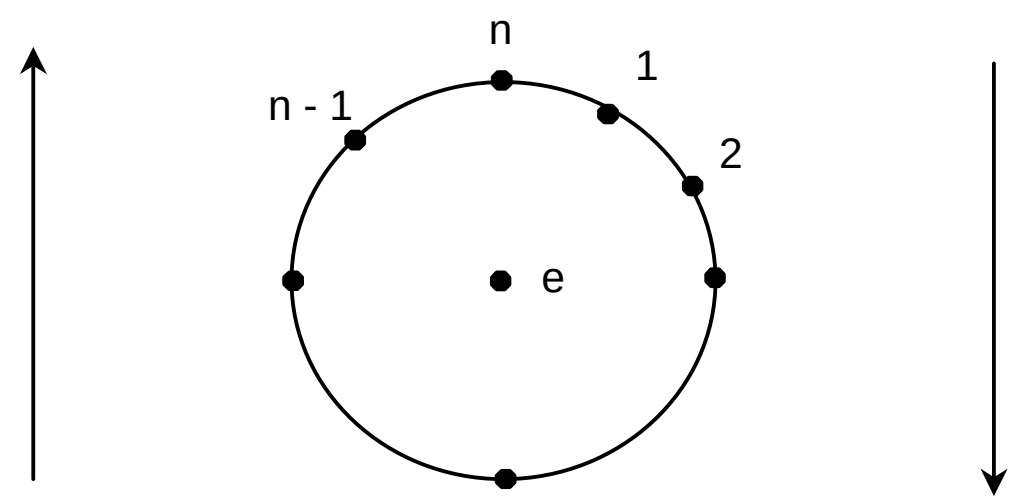

Let i, j ∈ $L_n(m)$ (i ≠ j, i ≠ e, j ≠ e). If j is the $r^{th}$ element from i counting in the clockwise direction the i • j will be the $t^{th}$ element from j in the clockwise direction where t = (m –1)r. We see that in general i • j need not be equal to j • i. When i = j we define $i^2$ = e and i • e = e • i = i for all i ∈ $L_n(m)$ and e acts as the identity in $L_n(m)$.

***Example 1.2.2***: Now the loop $L_7(4)$ is given by the following table:

| • | e | 1 | 2 | 3 | 4 | 5 | 6 | 7 |
|---|---|---|---|---|---|---|---|---|
| e | e | 1 | 2 | 3 | 4 | 5 | 6 | 7 |
| 1 | 1 | e | 5 | 2 | 6 | 3 | 7 | 4 |
| 2 | 2 | 5 | e | 6 | 3 | 7 | 4 | 1 |
| 3 | 3 | 2 | 6 | e | 7 | 4 | 1 | 5 |
| 4 | 4 | 6 | 3 | 7 | e | 1 | 5 | 2 |
| 5 | 5 | 3 | 7 | 4 | 1 | e | 2 | 6 |
| 6 | 6 | 7 | 4 | 1 | 5 | 2 | e | 3 |
| 7 | 7 | 4 | 1 | 5 | 2 | 6 | 3 | e |

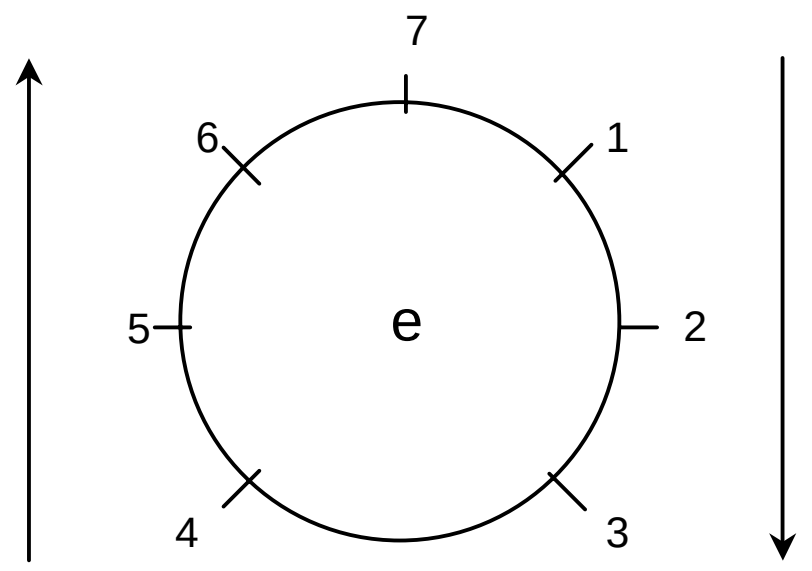

Let 2, 4 $\in$ $L_7(4)$. Now 4 is the $2^{nd}$ element from 2 in the clockwise direction. So 2.4 will be (4 – 1)2 that is the $6^{th}$ element from 4 in the clockwise direction which is 3.

Hence 2.4 = 3.

**Notation**: Let $L_n$ denote the class of loops. $L_n(m)$ for fixed n and various m's satisfying the conditions m < n, (m, n) = 1 and (m – 1, n) = 1, that is $L_n = \{L_n(m) \mid n > 3, n \text{ odd}, m < n, (m, n) = 1 \text{ and } (m-1, n) = 1\}$.

***Example 1.2.3:*** Let n = 5. The class $L_5$ contains three loops; viz. $L_5(2)$, $L_5(3)$ and $L_5(4)$ given by the following tables:

$L_5(2)$

| • | e | 1 | 2 | 3 | 4 | 5 |
|---|---|---|---|---|---|---|
| e | e | 1 | 2 | 3 | 4 | 5 |
| 1 | 1 | e | 3 | 5 | 2 | 4 |
| 2 | 2 | 5 | e | 4 | 1 | 3 |
| 3 | 3 | 4 | 1 | e | 5 | 2 |
| 4 | 4 | 3 | 5 | 2 | e | 1 |
| 5 | 5 | 2 | 4 | 1 | 3 | e |

$L_5(3)$

| • | e | 1 | 2 | 3 | 4 | 5 |
|---|---|---|---|---|---|---|
| e | e | 1 | 2 | 3 | 4 | 5 |
| 1 | 1 | e | 4 | 2 | 5 | 3 |
| 2 | 2 | 4 | e | 5 | 3 | 1 |
| 3 | 3 | 2 | 5 | e | 1 | 4 |
| 4 | 4 | 5 | 3 | 1 | e | 2 |
| 5 | 5 | 3 | 1 | 4 | 2 | e |

$L_5(4)$

| • | e | 1 | 2 | 3 | 4 | 5 |
|---|---|---|---|---|---|---|
| e | e | 1 | 2 | 3 | 4 | 5 |
| 1 | 1 | e | 5 | 4 | 3 | 2 |
| 2 | 2 | 3 | e | 1 | 5 | 4 |
| 3 | 3 | 5 | 4 | e | 2 | 1 |
| 4 | 4 | 2 | 1 | 5 | e | 3 |
| 5 | 5 | 4 | 3 | 2 | 1 | e |

**THEOREM [23]:** *Let $L_n$ be the class of loops for any $n > 3$, if $n = p_1^{\alpha_1} p_2^{\alpha_2} \ldots p_k^{\alpha_k}$ ($\alpha_i > 1$, for $i = 1, 2, \ldots, k$), then $|L_n| = \prod_{i=1}^{k} (p_i - 2)\, p_i^{\alpha_i - 1}$ where $|L_n|$ denotes the number of loops in $L_n$.*

The proof is left for the reader as an exercise.

**THEOREM [23]**: *$L_n$ contains one and only one commutative loop. This happens when $m = (n + 1) / 2$. Clearly for this $m$, we have $(m, n) = 1$ and $(m – 1, n) = 1$.*

It can be easily verified by using simple number theoretic techniques.

**THEOREM [23]**: *Let $L_n$ be the class of loops. If* $n = p_1^{\alpha_1} p_2^{\alpha_2} \ldots p_k^{\alpha_k}$ *, then $L_n$ contains exactly $F_n$ loops which are strictly non-commutative where* $F_n = \prod_{i=1}^{k} \left( p_i - 3 \right) p_i^{\alpha_i - 1}$.

The proof is left for the reader as an exercise.

Note: If n = p where p is a prime greater than or equal to 5 then in $L_n$ a loop is either commutative or strictly non-commutative. Further it is interesting to note if n = 3t then the class $L_n$ does not contain any strictly non-commutative loop.

**THEOREM [32]**: *The class of loops $L_n$ contains exactly one left alternative loop and one right alternative loop but does not contain any alternative loop.*

*Proof*: We see $L_n(2)$ is the only right alternative loop that is when m = 2 (Left for the reader to prove using simple number theoretic techniques). When m = n –1 that is $L_n(n – 1)$ is the only left alternative loop in the class of loops $L_n$.

From this it is impossible to find a loop in $L_n$, which is simultaneously right alternative and left alternative. Further it is clear from earlier result both the right alternative loop and the left alternative loop is not commutative.

**THEOREM [23]:** *Let $L_n$ be the class of loops:*

i. *$L_n$ does not contain any Moufang loop*
ii. *$L_n$ does not contain any Bol loop*
iii. *$L_n$ does not contain any Bruck loop.*

The reader is requested to prove these results using number theoretic techniques.

**THEOREM [32]**: *Let $L_n(m) \in L_n$. Then $L_n(m)$ is a weak inverse property (WIP) loop if and only if $(m^2 - m + 1) \equiv 0 (mod\ n)$.*

*Proof*: It is easily checked that for a loop to be a WIP loop we have "if (xy)z = e then x(yz) = e where x, y, z $\in$ L." Both way conditions can be derived using the defining operation on the loop $L_n(m)$.

***Example 1.2.4:*** L be the loop $L_7(3)$ = {e, 1, 2, 3, 4, 5, 6, 7} be in $L_7$ given by the following table:

| • | e | 1 | 2 | 3 | 4 | 5 | 6 | 7 |
|---|---|---|---|---|---|---|---|---|
| e | e | 1 | 2 | 3 | 4 | 5 | 6 | 7 |
| 1 | 1 | e | 4 | 7 | 3 | 6 | 2 | 5 |
| 2 | 2 | 6 | e | 5 | 1 | 4 | 7 | 3 |
| 3 | 3 | 4 | 7 | e | 6 | 2 | 5 | 1 |
| 4 | 4 | 2 | 5 | 1 | e | 7 | 3 | 6 |
| 5 | 5 | 7 | 3 | 6 | 2 | e | 1 | 4 |
| 6 | 6 | 5 | 1 | 4 | 7 | 3 | e | 2 |
| 7 | 7 | 3 | 6 | 2 | 5 | 1 | 4 | e |

It is easily verified $L_7(3)$ is a WIP loop. One way is easy for $(m^2 - m + 1) \equiv 0 \pmod 7$ that is $9 - 3 + 1 = 9 + 4 + 1 \equiv 0 (mod\ 7)$. It is interesting to note that no loop in the class $L_n$ contain any associative loop.

**THEOREM [23]**: *Let $L_n$ be the class of loops. The number of strictly non-right (left) alternative loops is $P_n$ where $P_n = \prod_{i=1}^{k} (p_i - 3) p_i^{\alpha_i - 1}$ and $n = \prod_{i=1}^{k} p_i^{\alpha_i}$.*

The proof is left for the reader to verify.

Now we proceed on to study the associator and the commutator of the loops in $L_n$.

**THEOREM [23]**: *Let $L_n(m) \in L_n$. The associator $A(L_n(m)) = L_n(m)$.*

**DEFINITION 1.2.20:** *The Smarandache loop (S-loop) is defined to be a loop L such that a proper subset A of L is a subgroup (with respect to the same induced operation) that is $\phi \neq A \subset L$.*

For more literature about the new class of loops refer [32-37].

## 1.3 Groupoids and Smarandache Groupoids

In this section we just recall the notion of groupoids and S-groupoids. We also give some new classes of groupoids constructed using the set of modulo integers. This book uses in several examples the groupoids from these new classes of groupoids and just recalls the notion of Smarandache semi automaton. For more about S-groupoids please refer [].

**DEFINITION 1.3.1:** *Given an arbitrary set P a mapping of $P \times P$ into P is called a binary operation on P. Given such a mapping $\sigma$: $P \times P \rightarrow P$ we use it to define a product $*$ in P by declaring $a * b = c$ if $\sigma(a, b) = c$.*

**DEFINITION 1.3.2:** *A non empty set of elements G is said to form a groupoid if in G is defined a binary operation called the product denoted by $*$ such that $a * b \in G$ for all a, b $\in$ G.*

**DEFINITION 1.3.3:** *A groupoid G is said to be a commutative groupoid if for every a, b $\in$ G we have $a * b = b * a$.*

**DEFINITION 1.3.4:** *A groupoid G is said to have an identity element e in G if $a * e = e * a = a$ for all a $\in$ G.*

**DEFINITION 1.3.5:** *Let (G, $*$) be a groupoid a proper subset H $\subset$ G is a subgroupoid if (H, $*$) is itself a groupoid.*

**DEFINITION 1.3.6:** *A groupoid G is said to be a Moufang groupoid if it satisfies the Moufang identity (xy) (zx) = (x(yz))x for all x, y, z in G.*

**DEFINITION 1.3.7:** *A groupoid G is said to be a Bol groupoid if G satisfies the Bol identity ((xy) z) y = x ((yz) y) for all x, y, z in G.*

**DEFINITION 1.3.8:** *A groupoid G is said to be a P-groupoid if (xy) x = x (yx) for all x, y* $\in$ *G.*

**DEFINITION 1.3.9:** *A groupoid G is said to be right alternative if it satisfies the identity (xy) y = x (yy) for all x, y* $\in$ *G. Similarly we define G to be left alternative if (xx) y = x (xy) for all x, y* $\in$ *G.*

**DEFINITION 1.3.10:** *A groupoid G is alternative if it is both right and left alternative, simultaneously.*

**DEFINITION 1.3.11:** *Let (G, *) be a groupoid. A proper subset H of G is said to be a subgroupoid of G if (H, *) is itself a groupoid.*

**DEFINITION 1.3.12:** *A groupoid G is said to be an idempotent groupoid if* $x^2 = x$ *for all x* $\in$ *G.*

**DEFINITION 1.3.13:** *Let G be a groupoid. P a non empty proper subset of G, P is said to be a left ideal of the groupoid G if 1) P is a subgroupoid of G and 2) For all x* $\in$ *G and a* $\in$ *P, xa* $\in$ *P. One can similarly define right ideal of the groupoid G. P is called an ideal if P is simultaneously a left and a right ideal of the groupoid G.*

**DEFINITION 1.3.14:** *Let G be a groupoid A subgroupoid V of G is said to be a normal subgroupoid of G if*

i. *aV = Va*
ii. *(Vx)y = V(xy)*
iii. *y(xV) = (yx)V*

*for all x, y, a* $\in$ *V.*

**DEFINITION 1.3.15:** *A groupoid G is said to be simple if it has no non trivial normal subgroupoids.*

**DEFINITION 1.3.16:** *A groupoid G is normal if*

i. $xG = Gx$
ii. $G(xy) = (Gx)y$
iii. $y(xG) = (yx)G$ *for all* $x, y \in G$.

**DEFINITION 1.3.17:** *Let G be a groupoid H and K be two proper subgroupoids of G, with* $H \cap K = \phi$. *We say H is conjugate with K if there exists a* $x \in H$ *such that* $H = x K$ *or* $Kx$ *('or' in the mutually exclusive sense).*

**DEFINITION 1.3.18:** *Let* $(G_1, \theta_1), (G_2, \theta_2), \ldots , (G_n, \theta_n)$ *be n groupoids with* $\theta_i$ *binary operations defined on each* $G_i$, $i = 1, 2, 3, \ldots , n$. *The direct product of* $G_1, \ldots , G_n$ *denoted by* $G = G_1 \times \ldots \times G_n = \{(g_1, \ldots , g_n) \mid g_i \in G_i\}$ *by component wise multiplication on G, G becomes a groupoid.*

*For if* $g = (g_1, \ldots , g_n)$ *and* $h = (h_1, \ldots , h_n)$ *then* $g \bullet h = \{(g_1\theta_1h_1, g_2\theta_2h_2, \ldots , g_n\theta_nh_n)\}$. *Clearly,* $gh \in G$. *Hence G is a groupoid.*

**DEFINITION 1.3.19:** *Let G be a groupoid we say an element* $e \in G$ *is a left identity if* $ea = a$ *for all* $a \in G$. *Similarly we can define right identity of the groupoid G, if* $e \in G$ *happens to be simultaneously both right and left identity we say the groupoid G has an identity.*

**DEFINITION 1.3.20:** *Let G be a groupoid. We say a in G has right zero divisor if* $a * b = 0$ *for some* $b \neq 0$ *in G and a in G has left zero divisor if* $b * a = 0$. *We say G has zero divisors if* $a \bullet b = 0$ *and* $b * a = 0$ *for* $a, b \in G \setminus \{0\}$ *A zero divisor in G can be left or right divisor.*

**DEFINITION 1.3.21:** *Let G be a groupoid. The center of the groupoid* $C (G) = \{x \in G \mid ax = xa \text{ for all } a \in G\}$.

**DEFINITION 1.3.22:** *Let G be a groupoid. We say a, b $\in$ G is a conjugate pair if a = bx (or xa for some x $\in$ G) and b = ay (or ya for some y $\in$ G).*

**DEFINITION 1.3.23:** *Let G be a groupoid of order n. An element a in G is said to be right conjugate with b in G if we can find x, y $\in$ G such that a • x = b and b • y = a (x $\ast$ a = b and y $\ast$ b = a). Similarly, we define left conjugate.*

**DEFINITION 1.3.24:** *Let $Z^+$ be the set of integers. Define an operation $\ast$ on $Z^+$ by x $\ast$ y = mx + ny where m, n $\in Z^+$, m < $\propto$ and n < $\propto$ (m, n) = 1 and m $\neq$ n. Clearly {$Z^+$, $\ast$, (m, n)} is a groupoid denoted by $Z^+$ (m, n). We have for varying m and n get infinite number of groupoids of infinite order denoted by $\mathbb{Z}^+$.*

Here we define a new class of groupoids denoted by Z(n) using $Z_n$ and study their various properties.

**DEFINITION 1.3.25:** *Let $Z_n$ = {0, 1, 2, ... , n – 1} n $\geq$ 3. For a, b $\in Z_n \setminus \{0\}$ define a binary operation $\ast$ on $Z_n$ as follows. a $\ast$ b = ta + ub (mod n) where t, u are 2 distinct elements in $Z_n \setminus \{0\}$ and (t, u) =1 here ' + ' is the usual addition of 2 integers and ' ta ' means the product of the two integers t and a. We denote this groupoid by {$Z_n$, (t, u), $\ast$} or in short by $Z_n$ (t, u).*

It is interesting to note that for varying t, u $\in Z_n \setminus \{0\}$ with (t, u) = 1 we get a collection of groupoids for a fixed integer n. This collection of groupoids is denoted by Z(n) that is Z(n) = {$Z_n$, (t, u), $\ast$ | for integers t, u $\in Z_n \setminus \{0\}$ such that (t, u) = 1}. Clearly every groupoid in this class is of order n.

***Example 1.3.1:*** Using $Z_3$ = {0, 1, 2}. The groupoid {$Z_3$, (1, 2), $\ast$} = ($Z_3$ (1, 2)) is given by the following table:

| $\ast$ | 0 | 1 | 2 |
|---|---|---|---|
| 0 | 0 | 2 | 1 |
| 1 | 1 | 0 | 2 |
| 2 | 2 | 1 | 0 |

Clearly this groupoid is non associative and non commutative and its order is 3.

**THEOREM 1.3.1:** *Let $Z_n = \{0, 1, 2, \ldots, n\}$. A groupoid in Z (n) is a semigroup if and only if $t^2 \equiv t \pmod n$ and $u^2 \equiv u \pmod n$ for $t, u \in Z_n \setminus \{0\}$ and $(t, u) = 1$.*

**THEOREM 1.3.2:** *The groupoid $Z_n$ (t, u) is an idempotent groupoid if and only if $t + u \equiv 1 \pmod n$.*

**THEOREM 1.3.3:** *No groupoid in Z (n) has {0} as an ideal.*

**THEOREM 1.3.4:** *P is a left ideal of $Z_n$ (t, u) if and only if P is a right ideal of $Z_n$ (u, t).*

**THEOREM 1.3.5:** *Let $Z_n$ (t, u) be a groupoid. If $n = t + u$ where both t and u are primes then $Z_n$ (t, u) is simple.*

**DEFINITION 1.3.26:** *Let $Z_n = \{0, 1, 2, \ldots, n - 1\}$ $n \geq 3$, $n < \infty$. Define $*$ a closed binary operation on $Z_n$ as follows. For any a, $b \in Z_n$ define $a * b = at + bu \pmod n$ where (t, u) need not always be relatively prime but $t \neq u$ and $t, u \in Z_n \setminus \{0\}$.*

**THEOREM 1.3.6:** *The number of groupoids in $Z^*(n)$ is $(n - 1)(n - 2)$.*

**THEOREM 1.3.7:** *The number of groupoids in the class Z (n) is bounded by $(n - 1)(n - 2)$.*

**THEOREM 1.3.8:** *Let $Z_n$ (t, u) be a groupoid in $Z^*(n)$ such that $(t, u) = t$, $n = 2m$, $t / 2m$ and $t + u = 2m$. Then $Z_n$ (t, u) has subgroupoids of order 2m / t or n / t.*

*Proof:* Given n is even and $t + u = n$ so that $u = n - t$.
Thus $Z_n (t, u) = Z_n (t, n - t)$. Now using the fact $t \cdot Z_n = \left\{0, t, 2t, 3t, \ldots, \left(\frac{n}{t} - 1\right)t\right\}$ that is $t \,.\, Z_n$ has only n / t elements and these n / t elements from a subgroupoid. Hence $Z_n$ (t, n – t) where (t, n – t) = t has only subgroupoids of order n / t.

**DEFINITION 1.3.27:** *Let $Z_n = \{0, 1, 2, \ldots , n – 1\}$ $n \geq 3$, $n < \propto$. Define $*$ on $Z_n$ as $a * b = ta + ub \pmod n$ where t and $u \in Z_n \setminus \{0\}$ and t can also equal u. For a fixed n and for varying t and u we get a class of groupoids of order n which we denote by $Z^{**}(n)$.*

**DEFINITION 1.3.28:** *Let $Z_n = \{0, 1, 2, \ldots , n – 1\}$ $n \geq 3$, $n < \propto$. Define $*$ on $Z_n$ as follows. $a * b = ta + ub \pmod n$ where $t, u \in Z_n$. Here t or u can also be zero.*

**DEFINITION 1.3.29:** *A Smarandache groupoid (S-groupoid) G is a groupoid which has a proper subset S, $S \subset G$ such that S under the operations of G is a semigroup.*

***Example 1.3.2:*** Let (G, ∗) be a groupoid given by the following table:

| ∗ | 0 | 1 | 2 | 3 | 4 | 5 |
|---|---|---|---|---|---|---|
| 0 | 0 | 3 | 0 | 3 | 0 | 3 |
| 1 | 1 | 4 | 1 | 4 | 1 | 4 |
| 2 | 2 | 5 | 2 | 5 | 2 | 5 |
| 3 | 3 | 0 | 3 | 0 | 3 | 0 |
| 4 | 4 | 1 | 4 | 1 | 4 | 1 |
| 5 | 5 | 2 | 5 | 2 | 5 | 2 |

Clearly, $S_1 = \{0, 3\}$, $S_2 = \{1, 4\}$ and $S_3 = \{2, 5\}$ are proper subsets of G which are semigroups of G.

So (G, ∗) is a S-groupoid.

**DEFINITION 1.3.30:** *Let G be a S-groupoid if the number of elements in G is finite we say G is a finite S-groupoid, otherwise G is said to be an infinite S-groupoid.*

**DEFINITION 1.3.31:** *Let G be a S-groupoid. G is said to be a Smarandache commutative groupoid (S-commutative groupoid)*

*if there is a proper subset, which is a semigroup, is a commutative semigroup.*

**THEOREM 1.3.9:** *Let G be a commutative groupoid if G is a S-groupoid then G is a S-commutative groupoid. Conversely, if G is a S-commutative groupoid G need not in general be a commutative groupoid.*

**DEFINITION 1.3.32:** *Let (G, $*$) be a S-groupoid. A non-empty subset H of G is said to be a Smarandache subgroupoid (S-subgroupoid) if H contains a proper subset K $\subset$ H such that K is a semigroup under the operation $*$.*

**DEFINITION 1.3.33:** *Let G be a SG and V be a S-subgroupoid of G. V is said to be a Smarandache normal groupoid (S-normal groupoid) if aV = X and Va = Y for all a $\in$ G where both X and Y are S-subgroupoids of G.*

**THEOREM 1.3.10:** *Every S-normal groupoid is a Smarandache semi normal groupoid and not conversely.*

**DEFINITION 1.3.34:** *Let G be a S-groupoid. H and P be any two subgroupoids of G. We say H and P are Smarandache (S-semi conjugate subgroupoids) of G if*

i. *H and P are S-subgroupoids of G*
ii. *H = xP or Px or*
iii. *P = xH or Hx for some x $\in$ G.*

**DEFINITION 1.3.35:** *Let G be a S-groupoid. H and P be subgroupoids of G. We say H and P are S-conjugate subgroupoids of G if*

i. *H and P are S-subgroupoids of G*
ii. *H = xP or Px and*
iii. *P = xH or Hx.*

**DEFINITION 1.3.36:** *Let S be non empty set. Generate a free groupoid using S and denote it by <S>. Clearly the free*

*semigroup generated by the set S is properly contained in <S>; as in <S> we may or may not have the associative law to be true.*

**DEFINITION 1.3.37:** *A Semi Automaton is a triple Y = (Z, A, δ) consisting of two non - empty sets Z and A and a function δ: Z × A → Z, Z is called the set of states, A the input alphabet and δ the "next state function" of Y.*

**DEFINITION 1.3.38:** *$Y_s = (Z, \overline{A}_s, \overline{\delta}_s)$ is said to be a Smarandache semi automaton (S-semi automaton) if $\overline{A}_s = \langle A \rangle$ is the free groupoid generated by A with Λ the empty element adjoined with it and $\overline{\delta}_s$ is the function from $Z \times \overline{A}_s \to Z$. Thus the S-semi automaton contains $Y = (Z, \overline{A}, \overline{\delta})$ as a new semi automaton which is a proper sub-structure of $Y_s$.*

*Or equivalently, we define a S-semi automaton as one, which has a new semi automaton as a sub-structure.*

**DEFINITION 1.3.39:** *$\overline{Y}_s' = (Z_1, \overline{A}_s, \overline{\delta}'_s)$ is called the Smarandache sub semi automaton (S-subsemi automaton) of $\overline{Y}_s = (Z_2, \overline{A}_s, \overline{\delta}'_s)$ denoted by $\overline{Y}_s' \leq \overline{Y}_s$ if $Z_1 \subset Z_2$ and $\overline{\delta}'_s$ is the restriction of $\overline{\delta}_s$ on $Z_1 \times \overline{A}_s$ and $\overline{Y}_s'$ has a proper subset $\overline{H} \subset \overline{Y}_s'$ such that $\overline{H}$ is a new semi automaton.*

For more literature refer [36, 40].

**Chapter Two**

# N-GROUPS AND SMARANDACHE N-GROUPS

In this chapter for the first time we give definitions, properties and examples of N-groups and Smarandache N-groups. However we have given in some of our books about the generalizations from bistructures to N-structures. The bistructures can be totally derived when N takes the value two. Throughout this book N is a finite positive number, $N \geq 2$. This chapter has two sections; in the first section we introduce the definition of N-groups and its properties. In section two we give the notions of Smarandache N-groups.

## 2.1 Basic Definition of N-groups and their properties

In this section we introduce the notion of N-groups and give their properties. Several of the classical results are extended by defining appropriate substructures on the N-groups. Examples are given to illustrate these structures.

**DEFINITION 2.1.1:** *Let $\{G, *_1, \ldots, *_N\}$ be a non empty set with N binary operations. $\{G, *_1, \ldots, *_N\}$ is called a N-group if there exists N proper subsets $G_1, \ldots, G_N$ of G such that*

i. $G = G_1 \cup G_2 \ldots \cup G_N$.
ii. $(G_i, *_i)$ *is a group for* $i = 1, 2, \ldots, N$.

*We say proper subset of G if $G_i \not\subseteq G_j$ and $G_j \not\subseteq G_i$ if $i \neq j$ for $1 \leq i, j \leq N$. When N = 2 this definition reduces to the definition of bigroup.*

***Example 2.1.1:*** Take $G = G_1 \cup G_2 \cup G_3 \cup G_4 \cup G_5$ where $G_1 = S_3$, the symmetric group of degree 3. $G_2 = D_{2n}$ the dihedral group of order 2n. $G_3 = \langle g \mid g^p = 1\rangle$, $(p \neq n, p \nmid n)$; $G_4 = A_4$ the alternating subgroup of the symmetric group $S_4$ and $G_5$ the group of quaternions.

Then $\{G, *_1, \ldots, *_5\}$ is a 5-group. Clearly $G = G_1 \cup G_2 \cup G_3 \cup G_4 \cup G_5$ and each $(G_i, *_i)$ is a group. Suppose we have $\{G = A_3 \cup S_3 \cup A_4 \cup S_4 \cup D_{25}, *_1, *_2, \ldots, *_5\}$. G is not a 5-gorup for $A_3 \underset{\neq}{\subset} S_3$ and $A_4 \underset{\neq}{\subset} S_4$ so G is not a 5-group as $A_3$, $S_3$, $A_4$, $S_4$ and $D_{25}$ are not proper subsets of G.

Now we proceed on to define sub N group of a N-group, $N \geq 2$. (N = 2 is the case of bigroup).

**DEFINITION 2.1.2:** *Let $\{G, *_1, \ldots, *_N\}$ be a N-group. A subset H $(\neq \phi)$ of a N-group $(G, *_1, \ldots, *_N)$ is called a sub N-group if H itself is a N-group under $*_1, *_2, \ldots, *_N$ , binary operations defined on G.*

**THEOREM 2.1.1:** *Let $(G, *_1, \ldots, *_N)$ be a N-group. The subset $H \neq \phi$ of a N-group G is a sub N-group then $(H, *_i)$ in general are not groups for i =1, 2, …, N.*

*Proof:* Given $(G, *_1, \ldots, *_N)$ is a N-group. $H \neq \phi$ of G is a sub N-group of G. To show $(H, *_i)$ are not groups for i = 1, 2, …, N.

We give an example to prove this.

Take $G = \{\pm i, \pm j, \pm k, \pm 1, 0, 2, 3, \ldots, 10, g, g^2, g^3 (g^4 = 1)\} = \{G_1 = \{\pm i, \pm j, \pm k, \pm 1\}, \times\} \cup \{G_2 = (0, 1, 2, \ldots, 9)$, addition modulo $10\} \cup \{G_3 = (1, g, g^2, g^3), \times\}$ is a 3-group, i.e., $G = G_1 \cup G_2 \cup G_3$.

Take $H = \{0, 5, g^2, \pm 1\}$; clearly H is a subset of G and H is a sub-3-group. But H is not a group under '×' or '+' the operations on $G_1$ or $G_2$ or $G_3$. Thus H = {±1; ×} ∪ {0, 5, +} ∪ {1, $g^2$, ×} is a 3-group. Hence the claim.

Next we give a characterization theorem about sub N-group.

**THEOREM 2.1.2:** *Let (G, $*_1$, $*_2$, …, $*_N$) be a N-group (N ≥ 2). The subset H (≠ ϕ) of G is a sub N group of G if and only if there exists N proper subsets $G_1$, …, $G_N$ of G such that*

i. *$G = G_1 \cup G_2 \cup \ldots \cup G_N$ where ($G_i$, $*_i$) are groups; i = 1, 2, ..., N.*
ii. *($H \cap G_i$, $*_i$ ) is a subgroup of ($G_i$, $*_i$) for i = 1, 2, …, N.*

*Proof:* Let H (≠ ϕ) be a sub N-group of G then (H, $*_1$, …, $*_N$) is a N-group. Therefore there exists proper subsets of $H_1$, $H_2$, …, $H_N$ of H such that

a. $H = H_1 \cup H_2 \cup \ldots \cup H_N$.
b. ($H_i$, $*_i$) is a group; for i = 1, 2, …, N.

Now we choose $H_i = H \cap G_i$ (for i = 1, 2, …, N) then we see that $H_i$ is a subset of $G_i$ and by (ii) ($H_i$, $*_i$) is itself a group. Hence $H_i$ = ($H \cap G_i$, $*_i$) is a subgroup of ($G_i$, $*_i$), for i = 1, 2, …, N.

Conversely let (i) and (ii) of the statement of the theorem be true. To prove (H, $*_1$, …, $*_N$) is a N-group it is seen the definition for N-group is satisfied by the set H so (H, $*_1$, …, $*_N$) is a N-group (ii) by note below! . Hence the claim.

*Note:*

1. It is important to note that in the above theorem just proved condition (i) can be removed. We include it only for simplicity of understanding.
2. As $(H \cap G_1) \cup (H \cap G_2) = H$ using this for finite unions we get $(H \cap G_1) \cup (H \cap G_2) \cup \ldots \cup (H \cap G_N) = H$ and each ($H \cap G_i$, $*_i$) is a subgroup of ($G_i$, $*_i$) for i = 1, 2, …, N given in Theorem 2.1.2.

We now define the notion of commutative N-group.

**DEFINITION 2.1.3:** *Let (G, $*_1, \ldots, *_N$) be a N-group where G = $G_1 \cup G_2 \cup \ldots \cup G_N$. G is said to be commutative if each ($G_i$, $*_i$) is commutative, $1 \leq i \leq N$.*

***Example 2.1.2:*** Let G = $G_1 \cup G_2 \cup \ldots . \cup G_7$ where $G_1 = \langle g \mid g^2 = e\rangle$, $G_2 = \langle g \mid g^7 = e'\rangle$, $G_3 = (Z_{10}, +)$, $G_4 = (Z^+, +)$, $G_5 = (Z_{11}, +)$, $G_6 = \langle h \mid h^{19} = e_1\rangle$ and $G_7 = \{Z_5 \setminus \{0\}, \times\}$. Clearly G is a commutative 7-group.

*Note:* If (G, $*_1, \ldots, *_N$) is a N-group with G = $G_1 \cup G_2 \cup \ldots \cup G_N$, the order of N-group G is the number of distinct elements in G. It is denoted by o(G) or |G| .
In general

$$o(G) \neq o(G_1) + \ldots + o(G_N).$$

If o(G) = n and n < $\propto$ we say the N-group G is finite. If n = $\infty$ then we say the N-group G is infinite.

It is interesting to note that the N-group given in example 2.1.2 is infinite. Thus we can say a N-group is infinite even if one of the group ($G_i$, $*_i$) is infinite. If all the groups ($G_i$, $*_i$), i = 1, 2, …, N are finite then we say the N-group G is finite.

**DEFINITION 2.1.4:** *Let G = (G, $*_1, \ldots, *_N$) be a N-group. We say the N-group G is cyclic if each of the group ($G_i$, $*_i$) is cyclic for i = 1, 2, …, N (Here G = $G_1 \cup G_2 \cup \ldots \cup G_N$ and ($G_i$, $*_i$) is a group for i = 1, 2,…, N).*

***Example 2.1.3:*** Let G = $G_1 \cup G_2 \cup G_3 \cup G_4$ where

$G_1 = \langle g \mid g^{10} = 1\rangle$, cyclic group of order 10,
$G_2 = (Z_{15}, +)$ group under '+' modulo 15,
$G_3 = \langle h \mid h^{19} = e\rangle$ cyclic group of order 19 and
$G_4 = \langle Z_{11}, +\rangle$ group under addition modulo 11.

Clearly G is 4-group which is cyclic infact G is also a finite cyclic 4-group.

It is easily verified as in case of groups even in case of N-groups every cyclic N-group is abelian but every abelian N-group in general need not be a cyclic N-group.

***Example 2.1.4:*** Let G = $\{G_1 \cup G_2 \cup G_3, *_1, *_2, *_3\}$ where $G_1 = (Z_{12}, +)$, $G_2 = (Z, +)$ and $G_3 = \{a, b \mid a^2 = b^2 = 1, ab = ba\}$. Clearly G is a abelian 3-group which is not a cyclic 3-group.

It is still important to see that as in case of finite groups where the order of the subgroup divides the order of the group, we see in case of finite N-groups the order of the sub N-group need not in general divide the order of the N-group.

***Result:*** Let $(G, *_1, \ldots, *_N)$ be a finite N-group. Suppose $(H, *_1, *_2, \ldots, *_N)$ be a sub N-group of G then, o(H) in general need not divide the order of G.

We prove this only by an example.

***Example 2.1.5:*** Let $(G, *_1, \ldots, *_4)$ be a 4-group where $G = G_1 \cup G_2 \cup G_3 \cup G_4$ with $G_1 = S_3$, $G_2 = (Z_{11}, +)$, $G_3 = A_4$ and $G_4 = \langle g \mid g^5 = 1\rangle$. Clearly order of the 4-group is 34.

Consider H, a sub 4- group of the 4-group G with $H = H_1 \cup H_2 \cup H_3 \cup H_4$; where

$$H_1 = \left\langle e, \begin{pmatrix} 1 & 2 & 3 \\ 2 & 3 & 1 \end{pmatrix}, \begin{pmatrix} 1 & 2 & 3 \\ 3 & 1 & 2 \end{pmatrix} \right\rangle, H_2 = \langle 0 \rangle,$$

$$H_3 = \left\langle e, \begin{pmatrix} 1 & 2 & 3 & 4 \\ 2 & 1 & 4 & 3 \end{pmatrix}, \begin{pmatrix} 1 & 2 & 3 & 4 \\ 3 & 4 & 1 & 2 \end{pmatrix}, \begin{pmatrix} 1 & 2 & 3 & 4 \\ 4 & 3 & 2 & 1 \end{pmatrix} \right\rangle \text{ and}$$

$$H_4 = \langle 1 \rangle.$$

Clearly o (H) = 9 and 9 does not divide 34.

Just we make a mention that the classical Lagrange’s theorem for finite groups cannot hold good in case of finite N-groups.

Now we proceed on to define the notion of normal sub N-group.

**DEFINITION 2.1.5:** *Let $(G, *_1, \ldots, *_N)$ be a N-group where $G = G_1 \cup G_2 \cup \ldots \cup G_N$. Let $(H, *_1, \ldots, *_N)$ be a sub N-group of $(G, *_1, \ldots, *_N)$ where $H = H_1 \cup H_2 \cup \ldots \cup H_N$ we say $(H, *_1, \ldots, *_N)$ is a normal sub N-group of $(G, *_1, \ldots, *_N)$ if each $H_i$ is a normal subgroup of $G_i$ for i = 1, 2,…, N.*

*Even if one of the subgroups $H_i$ happens to be non normal subgroup of $G_i$ still we do not call H a normal sub-N-group of the N-group G.*

***Example 2.1.6:*** Let $(G, *_1, *_2, *_3, *_4)$ be a 4-group where $G = G_1 \cup G_2 \cup G_3 \cup G_4$ with $G_1 = A_4$, $G_2 = S_3$, $G_3 = S_5$ and $G_4 = \langle g \mid g^{12} = 1\rangle$. Take $H = H_1 \cup H_2 \cup H_3 \cup H_4$ where

$$H_1 = \left\{\begin{pmatrix}1 & 2 & 3 & 4\\ 1 & 2 & 3 & 4\end{pmatrix}, \begin{pmatrix}1 & 2 & 3 & 4\\ 2 & 1 & 4 & 3\end{pmatrix}, \begin{pmatrix}1 & 2 & 3 & 4\\ 3 & 4 & 1 & 2\end{pmatrix}, \begin{pmatrix}1 & 2 & 3 & 4\\ 4 & 3 & 2 & 1\end{pmatrix}\right\}$$

$$H_2 = \left\{\begin{pmatrix}1 & 2 & 3\\ 1 & 2 & 3\end{pmatrix}, \begin{pmatrix}1 & 2 & 3\\ 2 & 3 & 1\end{pmatrix}, \begin{pmatrix}1 & 2 & 3\\ 3 & 1 & 2\end{pmatrix}\right\},$$

$$H_3 = A_5 \text{ and}$$
$$H_4 = \{g^3, g^6, g^9, 1\}.$$

Clearly H is a normal sub 4-group of the 4-group G for each $H_i$ is normal in $G_i$, i = 1, 2, 3, 4.

Now o(G) = 150 and o(H) = 71. We see even the order of the normal sub N-group in general need not divide the order of the N-group when G is a finite N-group.

From the above example we can note down the following important result.

**Result:** Let $(G, *_1, *_2, \ldots, *_N)$ be a finite N-group. If $(H, *_1, *_2, \ldots, *_N)$ is a normal sub N-group of G, in general o(H) does not divide the order of G.

Now we proceed onto define the notion of homomorphism of N-groups.

**DEFINITION 2.1.6:** *Let ($G = G_1 \cup G_2 \cup \ldots \cup G_N$, $*_1$, $*_2$,…, $*_N$) and ($K = K_1 \cup K_2 \cup \ldots \cup K_N$, $*_1$,…, $*_N$) be any two N- groups. We say a map $\phi : G \to K$ to be a N-group homomorphism if $\phi$ | $G_i$ is a group homomorphism from $G_i$ to $K_i$ for i = 1, 2,…, N. i.e. $\phi \big|_{G_i} : G_i \to K_i$ is a group homomorphism of the group $G_i$ to the group $K_i$; for i = 1, 2, …, N.*

Thus it is important and interesting to note that even if one of the restrictions say $\phi|_{G_j} : G_j \to K_j$ is not defined as a group homomorphism we cannot have N-group homomorphism.

Another important thing to be noted about N-groups are that in general the classical Cauchy theorem for finite groups is not true in case of finite N-groups. Suppose n is the order of the N-group G and n is a composite number and suppose p is a prime such that p divides n then in general the N-group G may not contain a sub N-group of order p.

First we illustrate this by the following example.

***Example 2.1.7:*** Let ($G = G_1 \cup G_2 \cup G_3 \cup G_4$, $*_1$, $*_2$, $*_3$, $*_4$) be a 4-group, where

$G_1$ = $\langle g \mid g^6 = 1\rangle$ ; cyclic group of order 6,
$G_2$ = {$S_3$}, $S_3$ the symmetric group of order 3;
$G_3$ = {$Z_{10}$, +} the group of integers modulo 10 and
$G_4$ = {$D_{2.3}$, the dihedral group having six elements}
= {a, b | $a^2 = b^3 = 1$, bab = a}
= {1, a, ab, $ab^2$, b, $b^2$}.

Clearly order of the 4-group G is 28. The primes which divide 28 are 2 and 7. By inspection one cannot have a sub 4- group of order four for its order can be only greater than 4.

H = ($H_1 \cup H_2 \cup H_3 \cup H_4$, $*_1$, $*_2$, $*_3$, $*_4$) be a sub 4-group of the 4-group G where

$$H_1 = \{1, g^3\},$$

$$H_2 = \left\{\begin{pmatrix} 1 & 2 & 3 \\ 1 & 2 & 3 \end{pmatrix}\right\},$$

$$H_3 = \{0, 5\} \text{ and } H_4 = \{1, b, b^3\}.$$

Order of the sub 4-group is 7. Thus we have 7/28, a sub-4 group of order 7 but no sub-4-group of order 2 or 4.

Thus we have the following nice result.

**Result:** Let G be a N-group of finite order say n, where n is a composite number. Let p be a prime such that p / n. The N-group G in general may not have a sub N-group of order p.

Another unique observation which differentiates a group from a N-group is that if G is a N-group of order n, where n is a prime number; then also the N-group G can have sub N-group of order m; m < n, m need not in general be a prime. Also m can be prime.

This is the marked difference between N-groups and groups. We illustrate our claim by the following example.

***Example 2.1.8:*** Let ($G = G_1 \cup G_2 \cup G_3 \cup G_4, \cup G_5, *_1, *_2, *_3, *_4\ *_5$) be a finite 5-group. Let $G_1 = S_3$, $G_2 = \{Z_9, +\}$, $G_3 = A_4$, $G_4 = \langle g \mid g^8 = e\rangle$ and $G_5 = \langle g \mid g^9 = 1\rangle$. Clearly order of the 5-group is 41. 41 is a prime. But the 5-group G has proper sub 5-groups. For consider

$$H = H_1 \cup H_2 \cup H_3 \cup H_4 \cup H_5$$

where

$$H_1 = \left\{\begin{pmatrix} 1 & 2 & 3 \\ 1 & 2 & 3 \end{pmatrix}, \begin{pmatrix} 1 & 2 & 3 \\ 1 & 3 & 2 \end{pmatrix}\right\},\ H_2 = \{0, 3, 6\},$$

$$H_3 = \left\{\begin{pmatrix} 1 & 2 & 3 & 4 \\ 1 & 2 & 3 & 4 \end{pmatrix}, \begin{pmatrix} 1 & 2 & 3 & 4 \\ 2 & 3 & 1 & 4 \end{pmatrix}, \begin{pmatrix} 1 & 2 & 3 & 4 \\ 3 & 2 & 1 & 4 \end{pmatrix}\right\},$$

$H_4 = \{1\}$ and $H_5 = \{e, g^4\}$.

Clearly H is a sub 5- group and order of H is 11. H is itself a sub 5-group of order a prime 11.

Clearly 11 $\nmid$ 41.

We have seen both the classical theorem for finite groups viz. Lagrange Theorem and Cauchy theorem do not in general hold good for finite N-groups.

**DEFINITION 2.1.7:** *Let G = ($G_1 \cup G_2 \cup \ldots \cup G_N$, $*_1, \ldots, *_N$) be a finite N-group. If H = ($H_1 \cup H_2 \cup \ldots \cup H_N$, $*_1, \ldots, *_N$) is a sub N group such that o(H) /o(G) then we call H a Lagrange sub N-group. If every sub N group of G is a Lagrange sub N-group then we call G a Lagrange N-group. If G has atleast one Lagrange sub N group then we call G a weakly Lagrange N group. If G has no Lagrange N-group then we call G a Lagrange free N-group.*

The following theorem can easily be proved.

**THEOREM 2.1.3:** *All N-groups of prime order are Lagrange free N-groups.*

Now we proceed on to define Sylow N group.

Now we try, can we ever think of classical Sylow theorems for finite N-groups. We have already shown if order a 4-group is 28, 2 / 28, 4 / 28, 8 $\nmid$ 28 but the 4-group G has no sub 4-group of order 2 or 4; but 7 / 28 and infact the 4-group G had a sub 4-group of order 7. So at the outset we can only say in general the verbatim extension of Sylow theorems for N-groups is not possible. To this end we define the following new concept namely p-Sylow sub- N group for any N-group G.

**DEFINITION 2.1.8:** *Let G = {$G_1 \cup G_2 \cup \ldots \cup G_N$, $*_1, \ldots, *_N$} be a N-group of finite order n and p be a prime. If G has a sub N-group H of order $p^r$, $r \geq 1$ such that p / n then we say G has a p-Sylow sub N-group.*

*Note:* We do not as in case of usual finite groups demand $p^r / n$. The only simple thing we demand is $p / n$; even it may so happen that the order of H is $p^t$ and $p^t \nmid n$ still we call H a p-Sylow sub N-group. It may also happen that $p^t / n$ and still order of H is just $p^t$, t finite.

We are forced to define yet another new notion called p-Sylow pseudo sub-N-group.

**DEFINITION 2.1.9:** *Let $G = (G = G_1 \cup G_2 \cup \ldots \cup G_N, *_1, *_2, \ldots, *_N)$ be a N-group. Let order of G be finite say n. Suppose n is not a composite number or say even a composite number such that there is a prime p such that $p \nmid n$. If the N-group G has a sub N- group H of order p or $p^t$ $(t \geq 1)$ then we call H the pseudo p-Sylow sub-N- group of the N-group G.*

Now we illustrate all our definitions by the following examples.

***Example 2.1.9:*** Consider the 3-group $(G = G_1 \cup G_2 \cup G_3, *_1, *_2, *_3)$ where $G_1 = S_3$, $G_2 = A_4$ and $G_3 = \langle g \mid g^{11} = 1\rangle$. Clearly G is a 3-group of finite order; and o(G) is 29 which is a prime.

Consider $H = (H_1 \cup H_2 \cup H_3, *_1, *_2, *_3)$ a sub 3- group given by

$$H_1 = \left\{\begin{pmatrix}1 & 2 & 3\\ 1 & 2 & 3\end{pmatrix}, \begin{pmatrix}1 & 2 & 3\\ 2 & 1 & 3\end{pmatrix}\right\},$$

$$H_2 = \left\{\begin{pmatrix}1 & 2 & 3 & 4\\ 1 & 2 & 3 & 4\end{pmatrix}, \begin{pmatrix}1 & 2 & 3 & 4\\ 2 & 1 & 4 & 3\end{pmatrix}, \begin{pmatrix}1 & 2 & 3 & 4\\ 3 & 4 & 1 & 2\end{pmatrix}, \begin{pmatrix}1 & 2 & 3 & 4\\ 4 & 3 & 2 & 1\end{pmatrix}\right\}$$

and $H_3 = \{1\}$.

Clearly $H = H_1 \cup H_2 \cup H_3$ is sub 3- group of G and it is a pseudo 7-Sylow sub 3- group of G.

Now consider $K = K_1 \cup K_2 \cup K_3$ a proper subset of G, where

$$K_1 = \left\{ \begin{pmatrix} 1 & 2 & 3 \\ 1 & 2 & 3 \end{pmatrix}, \begin{pmatrix} 1 & 2 & 3 \\ 1 & 3 & 2 \end{pmatrix} \right\},$$

$$K_2 = \left\{ \begin{pmatrix} 1 & 2 & 3 & 4 \\ 1 & 2 & 3 & 4 \end{pmatrix}, \begin{pmatrix} 1 & 2 & 3 & 4 \\ 2 & 1 & 4 & 3 \end{pmatrix} \right\} \text{ and } K_3 = \{1\}.$$

Clearly $K = K_1 \cup K_2 \cup K_3$ is a pseudo 5-Sylow sub 3 group of G for o (K) = 5.

Now consider $P = P_1 \cup P_2 \cup P_3$ where

$$P_1 = \left\{ \begin{pmatrix} 1 & 2 & 3 \\ 1 & 2 & 3 \end{pmatrix} \right\},$$

$$P_2 = \left\{ \begin{pmatrix} 1 & 2 & 3 & 4 \\ 1 & 2 & 3 & 4 \end{pmatrix} \right\} \text{ and } P_3 = \langle g \mid g^{11} = 1 \rangle.$$

Clearly P is a pseudo 13-Sylow sub 3- group of the 3-group G.

Next we give an example of a N-group of composite order and show it has pseudo p-Sylow sub N-groups.

***Example 2.1.10:*** Let $(G = G_1 \cup G_2 \cup G_3 \cup G_4, *_1, *_2, *_3, *_4)$ be a 4-group where

$G_1$ = $D_{2.6} = \{a, b \mid a^2 = b^6 = 1, bab = a\}$, so $|G_1| = 12$,
$G_2$ = $\{Z_{12}$, under addition modulo 12$\}$,
$G_3$ = $\{A_4\}$ the alternating group of order 12 and
$G_4$ = $\{g / g^{10} = 1\}$.

Clearly G is a N-group or order 46, $46 = 2 \times 23$ a composite member.

Consider a subset $H = H_1 \cup H_2 \cup H_3 \cup H_4$ where

$$H_1 = \{1, a\}, \; H_2 = \{0, 6\}$$

$$H_3 = \left\{ \begin{pmatrix} 1 & 2 & 3 & 4 \\ 1 & 2 & 3 & 4 \end{pmatrix}, \begin{pmatrix} 1 & 2 & 3 & 4 \\ 2 & 1 & 4 & 3 \end{pmatrix} \right\} \text{ and } H_4 = \{g^5, 1\}.$$

Clearly H is a sub 4- group of the 4-group G but H is a pseudo 7-Sylow sub 4-group of G, for $7 \nmid 46$.

Let $K = K_1 \cup K_2 \cup K_3, \cup K_4$ where

$$K_1 = \{b^2\ b^4\ 1\},\ K_2 = \{0, 3, 6, 9\},$$

$$K_3 = \left\{ \begin{pmatrix} 1 & 2 & 3 & 4 \\ 1 & 2 & 3 & 4 \end{pmatrix}, \begin{pmatrix} 1 & 2 & 3 & 4 \\ 4 & 3 & 2 & 1 \end{pmatrix} \right\} \text{ and}$$

$$K_4 = \{1, g^2\ g^4\ g^6\ g^8\}.$$

K is a pseudo 13 Sylow sub 4-group of G, for $13 \nmid 46$.

Now consider $P = P_1 \cup P_2 \cup P_3 \cup P_4$ a proper subset of G, where

$$P_1 = \{1, b^2, b^9\},\ P_2 = \{0, 6\},$$

$$P_3 = \left\{ \begin{pmatrix} 1 & 2 & 3 & 4 \\ 1 & 2 & 3 & 4 \end{pmatrix}, \begin{pmatrix} 1 & 2 & 3 & 4 \\ 3 & 4 & 1 & 2 \end{pmatrix} \right\} \text{ and}$$

$$P_4 = \{1, g^2, g^4, g^6, g^8\}.$$

Let $V = V_1 \cup V_2 \cup V_3 \cup V_4$ where

$$V_1 = A_4,$$
$$V_2 = \{0, 2, 4, 6, 8, 10\},$$

$$V_3 = \left\{ \begin{pmatrix} 1 & 2 & 3 & 4 \\ 1 & 2 & 3 & 4 \end{pmatrix}, \begin{pmatrix} 1 & 2 & 3 & 4 \\ 1 & 3 & 4 & 2 \end{pmatrix}, \begin{pmatrix} 1 & 2 & 3 & 4 \\ 1 & 4 & 2 & 3 \end{pmatrix} \right\} \text{ and}$$

$$V_4 = \{1, g_5\}.$$

o(V) = 23 and 23 /46.

Clearly order P is 11. P is a sub 4-group of the 4-group G. In fact P is a pseudo 11-Sylow sub-4 group of G for 11 $\not|$ 46. Thus a N-group can have both p-Sylow sub N-group as well as a pseudo $p_1$-Sylow sub N-group.

The only criteria for this is that order of the N-group must be a composite number still we have the following interesting result.

**THEOREM 2.1.4:** *Let $G = (G = G_1 \cup G_2 \cup \ldots \cup G_N, *_1, *_2, \ldots, *_N)$ be a N-group of finite order n. (n a composite number). If p is a prime such that p / n; then in general G need not have a p-Sylow sub N-group H in G.*

*Proof:* By an example. Let $G = \{G_1 \cup G_2 \cup G_3, *_1, *_2, *_3\}$ where

$G_1 = \langle g \mid g^7 = 1\rangle$,
$G_2 = \{Z_{11}$, addition modulo 11$\}$ and
$G_3 = \{Z_5 \setminus \{0\}$ under multiplication modulo 5$\}$.

Clearly order of the 3-group G is of order 20. 5/20 and 2/20 and 4/20. It is easily verified G has no 2- Sylow sub 3-group or 5-Sylow sub 3-group.

Now we proceed on to define quasi-symmetric N-groups.

**DEFINITION 2.1.10:** *Let $G = G_1 \cup G_2 \cup \ldots \cup G_N , *_1, *_2, \ldots, *_N)$ be a N-group. If each of $G_i$ is a symmetric group or an alternating group then we call G a quasi- symmetric N-group.*

***Example 2.1.11:*** Let $G = (G_1 \cup G_2 \cup G_3 \cup G_4, \cup G_5, *_1, *_2, *_3, *_4\ *_5)$ be a 5-group. Let $G_1 = A_3$, $G_2 = S_4$, $G_3 = S_5$, $G_4 = A_6$ and $G_5 = A_7$. Clearly G is a quasi symmetric 5-group.

Now we proceed on to define the concept of symmetric N-group.

**DEFINITION 2.1.11:** *Let G = ($G_1 \cup G_2 \cup \ldots \cup G_N$, $*_1$, $*_2$, …, $*_N$) where each $G_i$ = $S_{n_i}$ a symmetric group of degree $n_i$. Then we call G the symmetric N-group of degree ($n_1$, $n_2$,…, $n_N$).*

Now we can extend the notion of the classical Cayley theorem for groups.

**THEOREM 2.1.5:** *Every N-group is isomorphic to an appropriate sub N-group of the symmetric N-group.*

*Proof:* Let G = (G = $G_1 \cup G_2 \cup \ldots \cup G_N$, $*_1$, $*_2$, …, $*_N$) be a N-group and ($S = S_{N_1} \cup S_{N_2} \cup \ldots \cup S_{N_N}, *_1, \ldots, *_N$) be a suitably chosen symmetric N-group.

Then by usual Cayleys theorem we can define a homomorphism $\phi$: G $\to$ S such that $\phi$ / $G_i$ is isomorphic to a subgroup in $S_{N_i}$ for i = 1, 2, …, N. Hence the claim.

To over come the problems of obtaining analogous results in case of N-groups like LaGrange theorem and Cauchy theorem we define a new notion called N-order in a N-group G.

**DEFINITION 2.1.12:** *Let (G = $G_1 \cup G_2 \cup \ldots \cup G_N$, $*_1$, $*_2$,…, $*_N$) be a N-group. Let H = ($H_1 \cup H_2 \cup \ldots \cup H_N$, $*_1$,…, $*_N$) be a sub N-group of the N-group G. We define N-order of H to be equal to o($H_1$) + o($H_2$) + … + o($H_N$) where $H_i = H \cap G_i$ , i = 1, 2, …, N and we denote the N-order of H = ($H_1 \cup H_2 \cup \ldots \cup H_N$,, $*_1$,…, $*_N$) by N(H); when H = G we get N (H) = N (G).*

Now we need yet another concept for extending Cauchy theorem.

**DEFINITION 2.1.13:** *Let (G, $*_1$, …, $*_N$) be a finite N-group (i.e. G = $G_1 \cup G_2 \cup \ldots \cup G_N$) and H = ($H_1 \cup H_2 \cup \ldots \cup H_N$, $*_1$,…, $*_N$) be a sub N-group of the N-group G, we say the N-order of H pseudo divides o(G) and we denote it by N(H) $/_p$ o(G).*

We illustrate this by the following example.

***Example 2.1.12:*** Let G = ($G_1 \cup G_2 \cup \ldots \cup G_N$, $*_1$, $*_2$, …, $*_N$) be a N-group where N = 5 with $G_1 = S_3$, $G_2 = A_4$, $G_3$ = {$Z_{12}$, under + modulo 12}, $G_4 = D_{2.6}$ and $G_5 = \{g \mid g^8 = 1\}$. Consider $H = H_1 \cup H_2 \cup H_3 \cup H_4 \cup H_5$ where

$$H_1 = \left\{\begin{pmatrix}1 & 2 & 3\\1 & 2 & 3\end{pmatrix}, \begin{pmatrix}1 & 2 & 3\\2 & 1 & 3\end{pmatrix}\right\},$$

$$H_2 = \left\{\begin{pmatrix}1 & 2 & 3 & 4\\1 & 2 & 3 & 4\end{pmatrix}, \begin{pmatrix}1 & 2 & 3 & 4\\2 & 1 & 4 & 3\end{pmatrix}, \begin{pmatrix}1 & 2 & 3 & 4\\3 & 4 & 1 & 2\end{pmatrix}, \begin{pmatrix}1 & 2 & 3 & 4\\4 & 3 & 2 & 1\end{pmatrix}\right\}$$

$$H_3 = \{0, 6\},\ H_4 = \{1, b^2, b^4\} \text{ and } H_5 = \{1, g^4\}.$$

Clearly $o(H_i) / o(G_i)$ for i = 1, 2, ..., 5 as $H_i$ are subgroups of $G_i$, $N(H) /_p\, o(G)$.

o(G) = 48 and N(H) = 12. In truth 12/48, but all sub 5- groups of G need not divide the order of the N-group G. For take K = $K_1 \cup K_2 \cup K_3, \cup K_4, \cup K_5$ where

$$K_1 = \left\{\begin{pmatrix}1 & 2 & 3\\1 & 2 & 3\end{pmatrix}, \begin{pmatrix}1 & 2 & 3\\2 & 3 & 1\end{pmatrix}, \begin{pmatrix}1 & 2 & 3\\3 & 1 & 2\end{pmatrix}\right\},$$

$$K_2 = \left\{\begin{pmatrix}1 & 2 & 3 & 4\\1 & 2 & 3 & 4\end{pmatrix}, \begin{pmatrix}1 & 2 & 3 & 4\\2 & 1 & 4 & 3\end{pmatrix}\right\},$$

$$K_3 = \{0, 2, 4, 6, 8, 10\},$$
$$K_4 = \{1, a\} \text{ and } K_5 = \{1, g^2, g^4, g^6\}.$$

Now o(K) = 16 and $o(K) \mid_p o(G)$ in fact 16/48.
Take $P = P_1 \cup P_2 \cup P_3 \cup P_4 \cup P_5$ where

$$P_1 = K_1,\ P_2 = K_2,\ P_3 = \{0, 4, 8\},\ P_4 = \{1, b^2, b^4\}$$
$$\text{and } P_5 = \{1, g^2\, g^4\, g^6\}.$$

$o(P) = 14$, $o(P) /_{p} o(G)$ for $14 \nmid 48$.
Thus we have $N(H) \neq o(H)$ where H is a sub N-group.

**THEOREM 2.1.6:** *Let $(G, *_1, *_2, \ldots, *_N)$ be a N-group with $G = G_1 \cup G_2 \cup \ldots \cup G_N$. $(H_1 \cup H_2 \cup \ldots \cup H_N, *_1, \ldots, *_N)$ be a sub N-group of G. In general $o(H) \neq N(H)$.*

*Proof:* From the above example 2.1.12 we see in all the cases of the subgroup H, P and K we have $o(H) \neq N(H)$; $o(K) \neq N(K)$ and $o(P) \neq N(P)$.

**THEOREM 2.1.7:** *Let $(G_1 \cup G_2 \cup \ldots \cup G_N, *_1, \ldots, *_N)$ be a N-group. N(H) the N-order of H, H a subset of G then $N(H) /_p o(G)$.*

*Proof:* We prove this by an example consider $G = G_1 \cup G_2 \cup G_3$ be a 3-group

$$G_1 = S_3,\ G_2 = \langle g \mid g^6 = 1\rangle \text{ and } G_3 = D_{2.6}.$$

$$H = \left\{ \begin{pmatrix} 1 & 2 & 3 \\ 1 & 2 & 3 \end{pmatrix}, \begin{pmatrix} 1 & 2 & 3 \\ 2 & 1 & 3 \end{pmatrix}, \begin{pmatrix} 1 & 2 & 3 \\ 1 & 3 & 2 \end{pmatrix}, \begin{pmatrix} 1 & 2 & 3 \\ 3 & 2 & 1 \end{pmatrix}, 1, g^3, g^5, a, ab, ab^4, b^3, b^5 \right\}.$$

$$o(G_1) = 6,\ o(G_2) = 6 \text{ and } o(G_3) = 12.$$

Clearly $o(H_1) / o(G_1)$, $o(H_2) / o(G_2)$ and $o(H_3) / o(G_3)$, $N(H) = 12$, $o(G) = 24$ . $N(H) / o(G)$ but $N(H) /_{p} o(G)$.

**THEOREM 2.1.8:** *Let $(G, *_1, \ldots, *_N)$ be a N-group where $G = G_1 \cup G_2 \cup \ldots \cup G_N$ and $(H = H_1 \cup \ldots \cup H_N, *_1, \ldots, *_N)$ be a sub N-group of G then o(G) is pseudo divisible by N-order N(H) of H.*

*Proof:* Follows from the fact if $H = H_1 \cup \ldots \cup H_N$ is a sub N-group of the N-group G then we have $o(H_i) / o(G_i)$, $i = 1, 2, \ldots, N$. So $N(H) /_p o(G)$.

Now we proceed on to define a new notion called ($p_1, \ldots, p_N$) - Sylow sub N-group of a N-group G.

**DEFINITION 2.1.14:** *Let G = ($G_1 \cup G_2 \cup \ldots \cup G_N$, $*_1, \ldots, *_N$) be N-group. Let H = $H_1 \cup H_2 \cup \ldots \cup H_N$ be a sub N-group of G. We say H is a ($p_1, \ldots, p_N$)- Sylow sub N-group of G if $H_i$ is a $p_i$ - Sylow subgroup of $G_i$, i = 1, 2, …, N.*

We illustrate this by the following example.

***Example 2.1.13:*** Let G = $G_1 \cup G_2 \cup G_3 \cup G_4$ be a 4-group. H = $H_1 \cup H_2 \cup H_3 \cup H_4$ be a sub N-group of G. Let $G_1 = S_3$, $G_2 = A_4$, $G_3 = D_{2,7}$ and $G_4 = \langle g \mid g^{18} = 1\rangle$. Here

$$H_1 = \left\{\begin{pmatrix}1 & 2 & 3\\ 1 & 2 & 3\end{pmatrix}, \begin{pmatrix}1 & 2 & 3\\ 1 & 3 & 2\end{pmatrix}\right\},$$

$$H_2 = \left\{\begin{pmatrix}1 & 2 & 3 & 4\\ 1 & 2 & 3 & 4\end{pmatrix}, \begin{pmatrix}1 & 2 & 3 & 4\\ 2 & 3 & 1 & 4\end{pmatrix}, \begin{pmatrix}1 & 2 & 3 & 4\\ 3 & 1 & 2 & 4\end{pmatrix}\right\},$$

$$H_3 = \{1, b, b^2, \ldots, b^6\} \text{ and}$$
$$H_4 = \{g^2, g^4, g^6, g^8, g^{10}, g^{12}, g^{14}, g^{16}\ 1\}.$$

Clearly H is a sub 4-group and H is a (2, 3, 7, 3)-Sylow sub N-group of G.
The following theorem is immediate.

**THEOREM 2.1.9:** *Let (G, $*_1, \ldots, *_N$) be a finite N-group where G = $G_1 \cup G_2 \cup \ldots \cup G_N$. Let ($p_1, \ldots, p_N$) be a set of N-primes such that $p_i^{\alpha} \mid o(G_i)$ for i = 1, 2,…, N, then (G, $*_1, \ldots, *_N$) has a ($p_1, \ldots, p_N$) -Sylow sub N-group of N order $\sum_{i=1}^{N} p_i^{\alpha_i}$ .*

Now we proceed on to define the notion of conjugate sub N-groups of a N-group G.

**DEFINITION 2.1.15:** *Let $G = (G_1 \cup G_2 \cup \ldots \cup G_N, *_1, \ldots, *_N)$ be a N-group. We say two sub N-groups $H = H_1 \cup H_2 \cup \ldots \cup H_N$ and $K = K_1 \cup \ldots \cup K_N$ are conjugate sub N-groups of G if $H_i$ is conjugate to $K_i$, for i = 1, 2, …, N as subgroups of $G_i$.*

Now using this we can easily get the extension of second Sylow theorem.

**THEOREM 2.1.10:** *Let $G = (G_1 \cup G_2 \cup \ldots \cup G_N, *_1, \ldots, *_N)$ be a N-group. The number of $(p_1, \ldots, p_N)$ - Sylow sub N-groups of G for the N tuple of primes is of the form $\{(1 + k_1 p_1), (1 + k_2 p_2), \ldots, (1 + k_N p_N)\}$.*

Now we proceed on to define the notion of normalizer of a N-group G.

**DEFINITION 2.1.16:** *Let $G = (G_1 \cup G_2 \cup \ldots \cup G_N, *_1, *_2, \ldots, *_N)$ be a N-group, the normalizer of a in G is the set $N(a) = \{x \in G / xa = ax\}$. N(a) is a subgroup of $G_i$ depending on the fact whether $a \in G_i$. If a is in more than one $G_i$, then one is not in a position to say the structure of N(a).*

It can be set as an open problem if $G = G_1 \cup G_2 \cup \ldots \cup G_N$ is a N-group and $a \in G$ i.e. $a \in \cap G_i$, $i \in \{1, 2, \ldots, N\}$ then what is the structure of N(a) i.e. i does it not take all values from 1 to N, it takes some values whenever a is present in a particular $G_i$.

Only if we take $a \in G_i$ such that $a \notin G_j$ if $j \neq i$ we can define normalizer of a and can have

$$C_a^i = \frac{o(G_i)}{o(N(a))}.$$

It is still interesting to study

$$C_a = \frac{o(G_i)}{o(N(a))}$$

when a belongs to more than one $G_i$ .

Now we proceed on to define the concept of cosets of a N-group, where $G = G_1 \cup G_2 \cup \ldots \cup G_N$ is a N-group. $H = H_1 \cup H_2 \cup \ldots \cup H_N$ be the sub N-group of G. The right N coset of H in G for some $a \in G$ is defined to be $H_a^c = \{h_1 a \mid h_1 \in H_1$ and $a \in \cap G_i\} \cup \{h_2 a / h_2 \in H_2$ and $a \in \cap G_i\} \ldots \cup \{h_N a / h_N \in H_N$ and $a \in \cap G_i\}$. If $a \notin$ every $G_i$ and say only a is in $G_i$, $G_j$ and $G_k$, $i < j < k$, then,

$$H_a = H_1 \cup H_2 \cup \ldots \cup \{h_i a / h_i \in H_i\} \cup H_{i+1} \ldots \cup \{h_j a \mid h_j \in H_j\} \cup H_{j+1} \cup \ldots \cup \{h_k a / h_k \in H_k\} \cup H_{k+1} \cup \ldots \cup H_N.$$

The N-cosets largely depends on the ways we choose a. If $\bigcap_{i=1}^{k} G_i \neq \phi$, then the N coset is say for $a \in G_i$, $Ha = H_1 a \cup H_2 a \cup \ldots \cup \{h_i a / h_i \in G\} \cup H_{i+1} a \cup \ldots \cup H_N a$.

On similar lines one can define the left N-coset of a sub N-group H of G.

The next natural question is if $H = H_1 \cup H_2 \cup \ldots \cup H_N$ and $K = K_1 \cup K_2 \cup \ldots \cup K_N$ are two sub N-groups of the N-group $G = G_1 \cup G_2 \cup \ldots \cup G_N$ how to define HK we define,

$$HK = \left\{ h_1k_1, h_2k_2, \ldots, h_Nk_N \,\middle|\, \begin{matrix} h_i \in H_i \\ k_i \in K_i \end{matrix} \text{ with } i = 1, 2, \ldots, N \right\}.$$

The condition when is HK a sub N group is as follows. We make a definition for it.

**DEFINITION 2.1.17:** *Let $G = G_1 \cup G_2 \cup \ldots \cup G_N$ be a N-group. Let $H = H_1 \cup H_2 \cup \ldots \cup H_N$ and $K = K_1 \cup K_2 \cup \ldots \cup K_N$ be any two sub N-groups of G. HK will be a sub N-group of G if and only if $H_iK_i$ and $K_i H_i$ is a subgroup of $G_i$ and $H_iK_i = K_iH_i$ for $i = 1, 2, \ldots, N$.*

Now we proceed on to define the notion of quotient N-group of a N-group G.

**DEFINITION 2.1.18:** *Let $G = G_1 \cup G_2 \cup \ldots \cup G_N$ be a N-group we say $N = N_1 \cup N_2 \cup \ldots \cup N_N$ is a normal sub N-group of G if and only if each $N_i$ is normal in $G_i$, $i = 1, 2, \ldots, N$. We define the quotient N-group,*

$$G / N \text{ as } \left[\frac{G_1}{N_1} \cup \frac{G_2}{N_2} \cup \ldots \cup \frac{G_N}{N_N}\right]$$

*which is also a sub-N-group.*

We can define several other analogous properties about N-groups.

## 2.2 Smarandache N-groups and some of its properties

In this section we for the first time define the notion of Smarandache N-groups and when N = 2 we get the Smarandache bigroup. For more about Smarandache bigroup please refer [39].

**DEFINITION 2.2.1:** *Let $G = (G_1 \cup G_2 \cup \ldots \cup G_N, *_1, *_2, \ldots, *_N)$ be a non empty set such that $G_i$ are proper subsets of G. $G = (G_1 \cup G_2 \cup \ldots \cup G_N, *_1, *_2, \ldots, *_N)$ is called a Smarandache N-group (S-N-group) if the following conditions are satisfied.*

i. *$(G_i, *_i)$ is a group $(1 \leq i \leq N)$.*
ii. *$(G_j, *_j)$ is a S-semigroup for some $j \neq i$, $1 \leq j \leq N$.*

***Example 2.2.1:*** Let $G = \{G_1 \cup G_2 \cup G_3 \cup G_4 \cup G_5, *_1, *_2, *_3, *_4, *_5\}$ be a Smarandache 5-group where

$G_1 = S(3)$ the S-semigroup,
$G_2 = \{Z_{10}$, addition modulo 10} is a group,
$G_3 = \{Z_{12}$, multiplication modulo 12} is a S-semigroup,
$G_4 = \{g \mid g^5 = 1\}$ and
$G_5 = \{Z$; under multiplication} is a S-semigroup.

Clearly G is a Smarandache 5-group.

It is important to note that in general every S-N-group G need not in general be a N-group.

In the example 2.2.1 clearly G is not a 5-group it is only a S-5-group.

Now we proceed onto define the notion of Smarandache sub N-group.

**DEFINITION 2.2.2:** *Let* $G = (G_1 \cup G_2 \cup \ldots \cup G_N, *_1, \ldots, *_N)$ *be a S-N-group; a proper subset P of G is said to be a Smarandache sub N-group (S-sub N-group) of G if* $P = P_1 \cup P_2 \cup \ldots \cup P_N$ *where* $P_i \subseteq G_i$, $i = 1, 2, \ldots, N$, *and at least one of the* $P_i$ *is a S-semigroup under the operations of* $G_i$. *In short a proper subset P of G is a S-sub N-group if P itself a Smarandache-N-group under the operations of G.*

**THEOREM 2.2.1:** *Let* $G = \{G_1 \cup G_2 \cup G_3 \cup \ldots \cup G_N, *_1, *_2 \ldots, *_N\}$ *be a S-N-group. Then G has a proper subset H which is N-group.*

*Proof:* Given $G = G_1 \cup G_2 \cup \ldots \cup G_N$, is a S-N-group. Let $H = H_1 \cup H_2 \cup \ldots \cup H_N$ where each $H_i$ is a subgroup of $G_i$, $i = 1, 2, \ldots, N$ then H is a N-group. Thus we see every S-N-group contains a N-group.

Thus we can also make another definition, which is as follows.

Now we proceed on to define the notion of Smarandache commutative N-group.

**DEFINITION 2.2.3:** *Let* $(G = G_1 \cup G_2 \cup \ldots \cup G_N, *_1, \ldots, *_N)$ *be a S-N-group. We say G is a Smarandache commutative N-group (S-commutative N-group) if all* $G_i$ *which are not semigroups are commutative and when* $G_j$ *are S-semigroups every proper subset which is a group in* $G_j$ *is commutative.*

Thus if all the $G_i$ happen to be commutative we say G is a S-commutative N-group.

***Example 2.2.2:*** Let $G = G_1 \cup G_2 \cup G_3 \cup G_4 \cup G_5$ be a S-5-group. Let

$G_1$ = $\langle g \mid g^{12} = 1\rangle$,
$G_2$ = $\{Z_{18}$, multiplication modulo 18$\}$,
$G_3$ = $\{A_3\}$,
$G_4$ = $\{Z$, under multiplication$\}$ and
$G_5$ = $\{Z_7 \setminus \{0\}$, under multiplication modulo 7$\}$.

Clearly G is a S-commutative 5-group. We cannot say all S-N-groups are commutative for we have S-N-groups which are non commutative.

We illustrate this by an example.

***Example 2.2.3:*** Let $G = G_1 \cup G_2 \cup G_3 \cup G_4$ be a S-4 group, where

$G_1$ = $S(3)$,
$G_2$ = $A_4$,
$G_3$ = $\{Z_{20}$, under multiplication$\}$ and

$$G_4 = \left\{ \begin{pmatrix} a & b \\ c & d \end{pmatrix} \,\middle|\, a, b, c, d, \in Z_7 \right\} \text{ under matrix multiplication.}$$

Clearly G is a S-4-group which is non commutative.

Now we proceed on to define the notion of Smarandache weekly commutative N-group.

**DEFINITION 2.2.4:** *Let $\{G = G_1 \cup G_2 \cup \ldots \cup G_N, *_1, \ldots, *_N\}$be a -N-group. We say G is Smarandache weakly commutative N-group (S-weakly commutative N-group) if in all the S-semigroup $G_j$, they have a proper subset which is a commutative group and the rest of $G_i$'s which are groups are commutative.*

Now we illustrate this by the following example.

***Example 2.2.4:*** Let $G = G_1 \cup G_2 \cup G_3$ be a S-3-group where

$G_1$ = S(3),
$G_2$ = {$Z_{19}$, under modulo addition 19} and
$G_3$ = {$Z_{24}$, under multiplication modulo 24}.

Clearly G is a weakly commutative S-3-group for S (3) has

$$P = \left\{ \begin{pmatrix} 1 & 2 & 3 \\ 1 & 2 & 3 \end{pmatrix}, \begin{pmatrix} 1 & 2 & 3 \\ 2 & 1 & 3 \end{pmatrix} \right\},$$

proper subset which is a commutative subgroup. But $K = S_3 \underset{\neq}{\subset}$ S(3) is not a commutative group. So $G = G_1 \cup G_2 \cup G_3$ is only a S-weakly commutative 3-group.

In view of this we have the following theorem.

**THEOREM 2.2.2:** *Every S-commutative N-group is a S-weakly commutative N-group. But a S-weakly commutative N-group in general need not be a S-commutative N-group.*

*Proof:* From the very definition it implies that every S-commutative N-group is a S-weakly commutative N-group.

To show a S-weakly commutative N-group in general is not a S-commutative N-group. We prove by an example. Consider $G = G_1 \cup G_2 \cup G_3$ where $G_1$ = S (3), $G_2 = \langle Z_{12}$, under multiplication modulo 12$\rangle$ and $G_3 = \langle g \mid g^9 = 1\rangle$. Clearly G is not S-commutative 3-group for S(3) contains $S_3$ which is a proper subset of S(3) but is not a commutative group but is only a non commutative group. So G is only a S-weakly commutative 3-group and not a S-commutative 3-group. Hence the claim.

We proceed on to define the notion of Smarandache cyclic N-group and Smarandache weakly cyclic N-group.

**DEFINITION 2.2.5:** *Let $G = (G_1 \cup G_2 \cup \ldots \cup G_N, *_1, *_2, \ldots, *_N)$ be a S- N-group. We say G is a Smarandache cyclic N-group (S-cyclic N-group) if $G_i$ is a group it must be cyclic and if*

*$G_j$ is a S-semigroup it must be a S-cyclic semigroup for $1 \leq i, j \leq N$.*

Recall a S-semigroup, S is cyclic if each of the proper subset P of S which is a group must be such that P is a cyclic group. Then we say the S-semigroup is a S-cyclic semigroup. If the S-semigroup S has atleast one of the proper subset P which is a cyclic group then we call S a Smarandache weakly cyclic semigroup.

Now we proceed onto define S-weakly cyclic N-group.

**DEFINITION 2.2.6:** *Let $G = (G_1 \cup G_2 \cup \ldots \cup G_N, *_1, *_2, \ldots, *_N)$ be a S-N-group. We say G is a Smarandache weakly cyclic N-group (S-weakly cyclic N-group) if every group in the collection $\{G_i\}$ is a cyclic group and every S-semigroup in the collection $\{G_j\}$ is a S-weakly cyclic semigroup.*

Recall a S-semigroup S is a S-weakly cyclic semigroup if S has at least one proper subset which is a group, is cyclic.

Now we illustrate these by examples.

***Example 2.2.5:*** Let $G = G_1 \cup G_2 \cup G_3 \cup G_4$ be a S-4 group where

$G_1$ = $\{Z_{20}$, under multiplication modulo 20$\}$,
$G_2$ = $\{g \mid g^8 = 1\}$,
$G_3$ = $\{Z_7$, under '+' modulo 7$\}$ and
$G_4$ = $\{Z_{12}$, under multiplication modulo 12$\}$.

Clearly G is a S-cyclic 4-group.

Now we give an example of a S-weakly cyclic N-group.

***Example 2.2.6:*** Let $G = G_1 \cup G_2 \cup G_3 \cup G_4 \cup G_5$; be a S-5-group, where

$G_1$ = S(3), $G_2$ = S (5), $G_3 = \{g \mid g^7 = 1\}$,
$G_4$ = $\{Z_{15}$, multiplication modulo 15$\}$ and
$G_5$ = $\{Z_{11} \setminus \{0\}$ under multiplication modulo 11$\}$.

Clearly G is only a S-weakly cyclic 5-group and is not a S-cyclic 5- group for S(3) the S-semigroup contains $S_3$ ($\underset{\neq}{\subset}$ S(3)) as a proper subset) which is a group that is not cyclic similarly S(5) contains subsets $S_5$ and $A_5$ which are groups that are not cyclic. Thus G is only a S-weakly cyclic 5-group which is not a S-cyclic 5-group.

It is interesting to note that: If G = ($G_1 \cup G_2 \cup \ldots \cup G_N$, $*_1$, $*_2$, …, $*_N$) is a S-cyclic N-group then it is a S-weakly cyclic N-group but a S-weakly cyclic N-group in general is not a S-cyclic N-group. The example 2.2.6 just discussed will show how S(3) is a S-semigroup such that all its subgroups are not cyclic for $S_3$ is not cyclic.

Now we proceed on to define the notion of Smarandache Lagrange N-group.

**DEFINITION 2.2.7:** *Let G = {$G_1 \cup G_2 \cup \ldots \cup G_N$, $*_1$, …, $*_N$} be a S-N group. The number of distinct elements in G is called the order of G. If o(G) is finite we say G is a finite S-N-group otherwise it is infinite.*

**DEFINITION 2.2.8:** *Let G = $G_1 \cup G_2 \cup \ldots \cup G_N$ be a finite S-N-group. If the order of every S-sub N-group H of G divides the order of G then we call G a Smarandache Lagrange N-group (S-Lagrange N-group). If the order of at least one of the S-sub N-groups of G say H divides order G then we call G a Smarandache weakly Lagrange N-group (S-weakly Lagrange N-group). We call the S-sub N group H of G such that o(H) / o(G) to be a Lagrange Smarandache sub N group (Lagrange S-sub N group).*

The following remark is very important which is evident from the definitions.

**Remark:** Every S-Lagrange N-group is a S-weakly Lagrange N-group.

Now we proceed on to define the notion of Smarandache p-Sylow N-groups.

**DEFINITION 2.2.9:** *Let $G = (G_1 \cup G_2 \cup \ldots \cup G_N, *_1, *_2, \ldots, *_N)$ be a finite S-N-group. Let p be a prime such that p divides the order of G. If G has a S-sub N-group H of order p or $p^t$ ($t > 1$) then we say G has a Smarandache p-Sylow sub N-group (S-p-Sylow N-group).*

It is very important to note that $p / o(G)$ but $p^t \not| \; o(G)$ for any ($t \geq 1$); still we may have S-p-Sylow sub N-groups having $p^t$ elements in them. Even if $p / o(G)$ still in case of S-N-groups we may or may not have S-p-Sylow sub N-groups in them.

We illustrate our claim by the following example.

***Example 2.2.7:*** Let $G = G_1 \cup G_2 \cup G_3$ where

$G_1$ = S (3),
$G_2$ = $Z_{10}$ , {semigroup under multiplication modulo 10} and
$G_3$ = $\langle g \mid g^6 = 1\rangle$.

Clearly G is a S-3-group, o(G) = 43. Now order of G is 43. G has sub N-groups which are not S-sub N-groups. Take

$$H_1 = \left\{\begin{pmatrix}1 & 2 & 3\\1 & 2 & 3\end{pmatrix}, \begin{pmatrix}1 & 2 & 3\\2 & 1 & 3\end{pmatrix}\right\},$$

$H_2$ = {1, 9 multiplication modulo 10} and $H_3 = \{g \mid g^6 = 1\}$.

Clearly $H = H_1 \cup H_2 \cup H_3$ is a sub 3-group of order 10. This example is given to prove a S-N-group of even prime order can have sub N-groups.

Now we see this S-3-group also has S-sub 3-groups, for take $K = K_1 \cup K_2 \cup K_3$ where $K_1 = S_3$, $K_2 = \{0, 2, 4, 6, 8\}$ under multiplication modulo 10 and $K_3 = \{1\}$. Clearly K is a S-N-group in fact a S-sub N-group and order of K is 12.

Now we proceed on to give example of S-Sylow sub N-groups of a S-N-group.

***Example 2.2.8:*** $G = G_1 \cup G_2 \cup G_3$ be a 3-group where

$G_1 = S_3$,
$G_2 = \{Z_{14}$, semigroup under multiplication modulo 14$\}$ and
$G_3 = \{g \mid g = 1\}$.

Clearly o(G) = 28, and G is a S-3 group. Now 2 / 28 and 4 / 28 but G has no S- sub 3-group of order 4. Also 7/28 but $7^2 \nmid 28$. Consider $H = H_1 \cup H_2 \cup H_3$ where

$$H_1 = \left\{ \begin{pmatrix} 1 & 2 & 3 \\ 1 & 2 & 3 \end{pmatrix}, \begin{pmatrix} 1 & 2 & 3 \\ 2 & 1 & 3 \end{pmatrix} \right\},$$

$$H_2 = \{0, 2, 4, 8\} \text{ and } H_3 = \{1\}.$$

Clearly H is a S-sub 3-group in fact H is a 7-Sylow sub 3-group of G as order of H is 7.

Now we proceed on to see how far we can go for extending the notion of Cauchy element in case of S-N-groups.

**DEFINITION 2.2.10:** *Let $G = G_1 \cup G_2 \cup \ldots \cup G_N$ be a S-N-group of finite order. Let $a \in G$; a is said to be a Smarandache special Cauchy element (S-special Cauchy element) of G if n / o(G) when n > 1 such that $a^n = 1$ (n is the least number). We call $a \in G$ to be the (Smarandache Cauchy element) S-Cauchy element of G if $a \in G_i$ and $n \mid o(G_i)$. We see that by no means we can associate at all times the S-Cauchy elements and S-special Cauchy element of a S-N-group.*

Also it is important to note that in general in a S-N group every element need not be a S-special Cauchy element or a S-Cauchy element of the S-N-group.

Now we proceed on to define the notion of getting Smarandache right coset of a S-N-group G.

**DEFINITION 2.2.11:** *Let G = ($G_1 \cup G_2 \cup \ldots \cup G_N$, $*_1$, $*_2$, …, $*_N$ ) be a S-N-group; let H = ($H_1 \cup H_2 \cup \ldots \cup H_N$, $*_1$, $*_2$, …, $*_N$ ) be a S-sub N group of G.*

*For any g $\in$ G we define the Smarandache right coset (S-right coset) of the N-group G as Ha = $H_1 \cup H_2 \cup \ldots \cup H_i\, a \cup \ldots \cup H_N$ if a $\in H_i$ alone if a $\in \bigcap_{i=1}^{K} H_i$, then Ha = $H_1'a \cup H_2'\, a \cup \ldots \cup H_N'\, a$; here $H_i = H'_i$ if $H_i$ is a group; $H_i' \subset H_i$ if $H_i$ is a S-semigroup, $H'_i$ is a group. Similarly we can define Smarandache left coset (S-left coset) a H. We say H is a Smarandache coset if a H = Ha.*

The following observation is important and interesting.

*Observation*: Let G = ($G_1 \cup G_2 \cup \ldots \cup G_N$, $*_1$, $*_2$, …, $*_N$) be a finite S-N-group. If m is a positive integer such that m / o(G) then

1. G need not have a S-sub N group of order m
2. Even if H is a S-sub N-group of G

then o(H) $\not|$ o(G) in general. We have given several examples to this effect.

In view of this we define a new notion called Smarandache non Lagrange N-group.

**DEFINITION 2.2.12:** *Let G = ($G_1 \cup G_2 \cup \ldots \cup G_N$, $*_1$, $*_2$, …, $*_N$ ) be a finite S-N-group. If the order of none of the S-sub N-groups of G divides the order of G then we call G a Smarandache non Lagrange N-group (S- non Lagrange N-group).*

We illustrate this by the following example.

***Example 2.2.9:*** Let G = $G_1 \cup G_2 \cup G_3$ be a S-3 group where $G_1 = \langle g \mid g^2 = 1\rangle$, $G_2 = \{g \mid g^4 = 1\}$ under $\times$ and $G_3 = Z_6$ be the semigroup under multiplication. Clearly G is S-3-group of order 12.

Now consider the S-sub 3-group H = $G_1 \cup \{1, g^2\} \cup \{0, 2, 4\}$. Clearly H is a S-sub semigroup and o(H) = 7 and o(H) $\nmid$ o(G). So G is a S-non Lagrange 3-group.

***Example 2.2.10:*** Let G = $G_1 \cup G_2 \cup G_3$, where

$G_1$ = $\langle g \mid g^2 = 1\rangle$,
$G_2$ = $\langle h \mid h^3 = 1\rangle$ and
$G_3$ = $\{Z_4$, under multiplication modulo 4$\}$.

Clearly G is S-3-group and o(G) = 9. G has no S- sub -3 group which divides o(G). But G has S-sub-3 groups.

For take H = $H_1 \cup H_2 \cup H_3$ where $H_1 = G_1$, $H_2 = G_2$ and $H_3 = \{0, 2\}$; H is a S-sub 3-group of order 7 and all other S-sub 3-groups are of order greater than 6.

However one can have a S-sub 3-group of order seven by taking K = $K_1 \cup K_2 \cup K_3$ where $K_1 = G_1$, $K_2 = G_2$ and $K_3 = \{0, 1, 3\}$. Clearly order of K is 7, 7 $\nmid$ 9.

Now we proceed on to define a new notion called Smarandache symmetric N-groups.

**DEFINITION 2.2.13:** *G = ($G_1 \cup G_2 \cup \ldots \cup G_N$, $*_1, *_2, \ldots, *_N$ ) where each $G_i$ is either S ($n_i$) or $S_{n_i}$ i.e. each $G_i$ is either a symmetric semigroup or a symmetric group, for i =1,2, …, N. Then G is defined as the Smarandache Symmetric N-group (S-symmetric N-groups).*

**THEOREM 2.2.3: (SMARANDACHE CAYLEY'S THEOREM FOR S-N-GROUPS):** *Let G = $G_1 \cup \ldots \cup G_N$, $*_1, *_2, \ldots, *_N$ ) be a S-N-group. Every S-N-group is embeddable in a S-symmetric N-group for a suitable $n_i$ , i = 1, 2, …, N.*

*Proof:* We know from the classical Cayley's theorem every group is embeddable in a symmetric group for an appropriate n. Also we know from the Cayley's theorem proved for S-semigroup every S-semigroup is isomorphic to a S-semigroup in S(n) for some n.

Hence we have every S-N-group is embeddable in a S-symmetric N-group.

Now we proceed on the define Smarandache inverse pair in a S-N-group G.

**DEFINITION 2.2.14:** *Let $G = (G_1 \cup G_2 \cup G_3 \cup \ldots \cup G_N, *_1, *_2, \ldots, *_N)$ be a S-N-group. An element $x \in G \setminus \{1\}$ is said to have a Smarandache inverse (S-inverse) $y \in G$, if $xy = 1$ and for $a, b \in G \setminus \{1, x, y\}$ we have $xa = y$ and $y b = x$ with $ab = 1$. The pair (x, y) is called the Smarandache-inverse pair (S-inverse pair).*

*If x be a S-inverse pair of y and (x, y) is a inverse pair with the related pair (a, b). If the pair (a, b) happens to be a S-inverse pair not necessarily with (x, y) as related pair we say (a, b) is the Smarandache co-inverse pair (S-co-inverse pair).*

**THEOREM 2.2.4 :** *Let $G = (G_1 \cup G_2 \cup \ldots \cup G_N, *_1, *_2, \ldots, *_N)$ be a S-N-group. Every S-inverse in $G_i$ or $H \subset G_j$ (H a subgroup of the S-semigroup $G_j$) has the inverse in the subset P of G. $P = G_i \cup G_j$ but every inverse in G need not have an S inverse.*

The proof is left as an exercise for the reader.

*Note:* The S-inverse may be in more than one $G_i$ or $H \subset G_j$.

Now we define the notion of Smarandache conjugate in S-N-groups.

**DEFINITION 2.2.15:** *Let $G = (G_1 \cup G_2 \cup \ldots \cup G_N, *_1, *_2, \ldots, *_N)$ be a S-N group. We say an element $x \in G$ has a Smarandache conjugate (S-conjugate) y in G if.*

i. *$xa = ay$ for some $a \in G$*
ii. *$ab = bx$ and $ac = cy$ for some b, c in G.*

*It is easy to verify if G is a S-symmetric N-group then G has S-elements which are S-conjugate.*

Several interesting results as in case of groups can be derived with suitable modifications in case of N-groups and S-N-groups.

Now we just give the definition of N-semigroups and S-N-semigroups and just mention a very few of its important properties.

**DEFINITION 2.2.16:** *Let $S = (S_1 \cup \ldots \cup S_N, *_1, \ldots, *_N)$ where S is a non empty set on which is defined, N-binary operations $*_1, *_2, \ldots, *_N$. S is called a N-semigroup if the following condition are true*

i. *$S = S_1 \cup S_2 \cup \ldots \cup S_N$; $S_i$'s are proper subsets of S. $1 \leq i \leq N$.*
ii. *$(S_i, *i)$ are semigroups, $i = 1, 2, \ldots, N$.*

***Example 2.2.11:*** Let $\{S = S_1 \cup S_2 \cup S_3, *_1, *_2, *_3\}$ where

$S_1$ = $\{Z_{10}$ semigroup under $\times$ modulo 10$\}$,
$S_2$ = S(5) and
$S_3$ = {Z, semigroup under multiplication}.

S is a 3-semigroup.

Now we just give the definition of Smarandache N- semigroup.

**DEFINITION 2.2.17:** *Let $S = \{S_1 \cup S_2 \cup \ldots \cup S_N, *_1, *_2, \ldots, *_N\}$ where S is a non empty set and $*_1, \ldots, *_N$ are N-binary operations defined on S. S is said to be Smarandache N semigroup (S-N semigroup) if the following conditions are satisfied*

i. *$S = S_1 \cup S_2 \cup \ldots \cup S_N$ is such that $S_i$'s are proper subsets of S.*

ii. *Some of ($S_i$, $*i$) are groups and some of ($S_j$, $*j$) are S-semigroups, $1 \leq i, j \leq N$. ($i \neq j$)*

We illustrate this by the following example

***Example 2.2.12:*** Let S = $\{S_1 \cup S_2 \cup S_3 \cup S_4, *_1, *_2, *_3, *_4\}$ where

$S_1$ = S(5),
$S_2$ = $\{Z_{12}$; a semigroup under multiplication modulo 12$\}$,
$S_3$ = $\left\{\begin{pmatrix} a & b \\ c & d \end{pmatrix} \mid a, b, c, d \in Q\right\}$ Semigroup under matrix multiplication and
$S_4$ = $\{(a, b) \mid a \in Z_{10}$ and $b \in Z_{12}\}$ a semigroup under component wise multiplication.

Clearly S is a S-4 semigroup of infinite order.

All results derived for groups and groupoids can be derived in case of N-semigroups and S-N-semigroups with appropriate modifications. Problems are suggested in chapter six of this book.

**Chapter Three**

# N-LOOPS AND SMARANDACHE N-LOOPS

In this chapter we introduce the notion of N-loops and Smarandache N-loops. We give some of its properties and illustrate them with examples. This chapter has two sections in section one we introduce the notion of N-loops and in section two we define Smarandache N- loops and give some of its properties.

## 3.1 Definition of N-loops and their properties

In this section we define the notion of N-loops and give some of their properties. We mainly study Sylow and Lagrange properties for N- loops for they alone can give some special nature of these structures, which will be interesting and innovative.

**DEFINITION 3.1.1:** *Let $(L, *_1, \ldots, *_N)$ be a non empty set with N binary operations $*_i$. L is said to be a N loop if L satisfies the following conditions:*

i. *$L = L_1 \cup L_2 \cup \ldots \cup L_N$ where each $L_i$ is a proper subset of L; i.e., $L_i \not\subseteq L_j$ or $L_j \not\subseteq L_i$ if $i \neq j$ for $1 \leq i, j \leq N$.*
ii. *$(L_i, *_i)$ is a loop for some i, $1 \leq i \leq N$.*
iii. *$(L_j, *_j)$ is a loop or a group for some j, $1 \leq j \leq N$.*

*For a N-loop we demand atleast one ($L_j$, $*_j$ ) to be a loop.*

***Example 3.1.1:*** Let ($L = L_1 \cup L_2 \cup L_3$, $*_1$, $*_2$, $*_3$) be a 3-loop where ($L_1$, $*_1$) is a loop given by the following table:

| $*_1$ | e | $g_1$ | $g_2$ | $g_3$ | $g_4$ | $g_5$ |
|---|---|---|---|---|---|---|
| e | e | $g_1$ | $g_2$ | $g_3$ | $g_4$ | $g_5$ |
| $g_1$ | $g_1$ | e | $g_3$ | $g_5$ | $g_2$ | $g_4$ |
| $g_2$ | $g_2$ | $g_5$ | e | $g_4$ | $g_1$ | $g_3$ |
| $g_3$ | $g_3$ | $g_4$ | $g_1$ | e | $g_5$ | $g_2$ |
| $g_4$ | $g_4$ | $g_3$ | $g_5$ | $g_2$ | e | $g_1$ |
| $g_5$ | $g_5$ | $g_2$ | $g_4$ | $g_1$ | $g_3$ | e |

$L_2$ = {$Z_{10}$, group under '+' modulo 10} and
$L_3$ = {$S_3$} symmetric group of degree 3.

Clearly L is a 3-loop.

**DEFINITION 3.1.2:** *Let ($L = L_1 \cup L_2 \cup \ldots \cup L_N$, $*_1$, …, $*_N$) be a N-loop. L is said to be a commutative N-loop if each ($L_i$, $*_i$) is commutative, i = 1, 2, …, N. We say L is inner commutative if each of its proper subset which is N-loop under the binary operations of L are commutative.*

***Example 3.1.2:*** Let ($L = L_1 \cup L_2 \cup L_3 \cup L_4$, $*_1$, $*_2$, $*_3$, $*_4$) be a 4-loop, where $L_1 = \langle g \mid g^{12} = 1\rangle$, cyclic group of order 12, $L_2$ = {$Z_{15}$, group under '+' modulo 15}, $L_3$ is given by the following table:

| $*_3$ | e | $g_1$ | $g_2$ | $g_3$ | $g_4$ | $g_5$ |
|---|---|---|---|---|---|---|
| e | e | $g_1$ | $g_2$ | $g_3$ | $g_4$ | $g_5$ |
| $g_1$ | $g_1$ | e | $g_4$ | $g_2$ | $g_5$ | $g_3$ |
| $g_2$ | $g_2$ | $g_4$ | e | $g_5$ | $g_3$ | $g_1$ |
| $g_3$ | $g_3$ | $g_2$ | $g_5$ | e | $g_1$ | $g_4$ |
| $g_4$ | $g_4$ | $g_5$ | $g_3$ | $g_1$ | e | $g_2$ |
| $g_5$ | $g_5$ | $g_3$ | $g_1$ | $g_4$ | $g_2$ | e |

$L_3$ is a commutative loop. $L_4$ is given by the following table:

| $*_4$ | e | $a_1$ | $a_2$ | $a_3$ | $a_4$ | $a_5$ | $a_6$ | $a_7$ |
|---|---|---|---|---|---|---|---|---|
| e | e | $a_1$ | $a_2$ | $a_3$ | $a_4$ | $a_5$ | $a_6$ | $a_7$ |
| $a_1$ | $a_1$ | e | $a_5$ | $a_2$ | $a_6$ | $a_3$ | $a_7$ | $a_4$ |
| $a_2$ | $a_2$ | $a_5$ | e | $a_6$ | $a_3$ | $a_7$ | $a_4$ | $a_1$ |
| $a_3$ | $a_3$ | $a_2$ | $a_6$ | e | $a_7$ | $a_4$ | $a_1$ | $a_5$ |
| $a_4$ | $a_4$ | $a_6$ | $a_3$ | $a_7$ | e | $a_1$ | $a_5$ | $a_2$ |
| $a_5$ | $a_5$ | $a_3$ | $a_7$ | $a_4$ | $a_1$ | e | $a_2$ | $a_6$ |
| $a_6$ | $a_6$ | $a_7$ | $a_4$ | $a_1$ | $a_5$ | $a_2$ | e | $a_3$ |
| $a_7$ | $a_7$ | $a_4$ | $a_1$ | $a_5$ | $a_2$ | $a_6$ | $a_3$ | e |

Clearly $L_4$ is also a commutative loop so L is a 4-loop, which is commutative.

It is important to note that the number of distinct elements in a N-loop L gives the cardinality of L. If L has finite cardinality then L is a finite N-loop. If L has infinite cardinality then L is an infinite N-loop.

Now we proceed on to define certain types of N-loops.

**DEFINITION 3.1.3:** *Let $L = \{L_1 \cup L_2 \cup \ldots \cup L_N, *_1, *_2, \ldots, *_N\}$ be a N-loop. We say L is a Moufang N-loop if all the loops ($L_i$, $*_i$ ) satisfy the following identities.*

i. $(xy)(zx) = (x(yz))x$
ii. $((xy)z)y = x(y(zy))$
iii. $x(y(xz)) = ((xy)x)z$

*for all $x, y, z \in L_i$, $i = 1, 2, \ldots, N$.*

Now we proceed on to define a Bruck N-loop.

**DEFINITION 3.1.4:** *Let $L = (L_1 \cup L_2 \cup \ldots \cup L_N, *_1, \ldots, *_N)$ be a N-loop. We call L a Bruck N-loop if all the loops ($L_i$, $*_i$ ) satisfy the identities*

i. $(x(yz))z = x(y(xz))$
ii. $(xy)^{-1} = x^{-1}y^{-1}$

*for all* $x, y \in L_i$, $i = 1, 2, \ldots, N$.

Similarly we can define Bol N-loop, WIP-N-loop and alternative N-loop.

Now we proceed on to define sub N-loop of a N-loop more seriously.

**DEFINITION 3.1.5:** *Let* $L = (L_1 \cup L_2 \cup \ldots \cup L_N, *_1, \ldots, *_N)$ *be a N-loop. A non empty subset P of L is said to be a sub N-loop, if P is a N-loop under the operations of L i.e.,* $P = \{P_1 \cup P_2 \cup P_3 \cup \ldots \cup P_N, *_1, \ldots, *_N\}$ *with each* $\{P_i, *_i\}$ *is a loop or a group.*

We give the conditions under which a proper subset P of a N-loop L is a sub N loop of L.

**THEOREM 3.1.1:** *Let* $L = \{L_1 \cup L_2 \cup \ldots \cup L_N, *_1, \ldots, *_N\}$ *be a N-loop. A non empty subset* $P = \{P_1 \cup P_2 \cup \ldots \cup P_N, *_1, \ldots, *_N\}$ *of L is a sub N-loop of L if and only if*

i. $P_i = P \cap L_i$ *is a subloop or subgroup of* $L_i$, $i = 1, 2, \ldots, N$.

*Proof:* Clearly if each $P_i$ satisfies the condition $P_i = P \cap L_i$ is a subloop of $L_i$, i = 1, 2, …, N then P is a sub N-loop. Now suppose $P = \{P_1 \cup P_2 \cup \ldots \cup P_N, *_1, \ldots, *_N\}$ be a sub N-loop of L. Then clearly we have the two conditions to be true

(1) Each $P_i$ is a loop or a group under the operations of $*_i$, i = 1, 2, …, N.
(2) Also $P \cap L_i = P_i$ is true. Hence the claim.

Now we define normal sub N-loop of a N-loop L.

**DEFINITION 3.1.6:** *Let* $L = \{L_1 \cup L_2 \cup \ldots \cup L_N, *_1, \ldots, *_N\}$ *be a N-loop. A proper subset P of L is said to be a normal sub N-loop of L if*

i. *P is a sub N-loop of L*
ii. $x_i P_i = P_i x_i$ *for all* $x_i \in L_i$
iii. $y_i (x_i P_i) = (y_i x_i) P_i$ *for all* $x_i y_i \in L_i$.

*A N-loop is said to be a simple N-loop if L has no proper normal sub N-loop.*

Now we proceed on to define the notion of Moufang center.

**DEFINITION 3.1.7:** *Let* $L = \{L_1 \cup L_2 \cup \ldots \cup L_N, *_1, \ldots, *_N\}$ *be a N-loop. We say* $C_N(L)$ *is the Moufang N-centre of this N-loop if* $C_N(L) = C_1(L_1) \cup C_2(L_2) \cup \ldots \cup C_N(L_N)$ *where* $C_i(L_i) = \{x_i \in L_i / x_i y_i = y_i x_i$ *for all* $y_i \in L_i\}$, $i = 1, 2, \ldots, N$.

**DEFINITION 3.1.8:** *Let L and P to two N-loops i.e.* $L = \{L_1 \cup L_2 \cup \ldots \cup L_N, *_1, \ldots, *_N\}$ *and* $P = \{P_1 \cup P_2 \cup \ldots \cup P_N, o_1, \ldots, o_N\}$. *We say a map* $\theta : L \to P$ *is a N-loop homomorphism if* $\theta = \theta_1 \cup \theta_2 \cup \ldots \cup \theta_N$ *'$\cup$' is just a symbol and* $\theta_i$ *is a loop homomorphism from* $L_i$ *to* $P_i$ *for each* $i = 1, 2, \ldots, N$.

Here for a N-loop homomorphism to be defined we must have both L and P to be N-loops and if $L_i$ is a loop then $P_i$ is a loop if $L_j$ is a group then $P_j$ must be a group otherwise one cannot define N-loop homomorphism of N-loops.

We define the notion of weak Moufang N-loops, weak Bol N-loops and weak Bruck N-loops.

**DEFINITION 3.1.9:** *Let* $L = \{L_1 \cup L_2 \cup \ldots \cup L_N, *_1, \ldots, *_N\}$ *be a N-loop. We say L is weak Moufang N-loop if there exists atleast a loop* $(L_i, *_i)$ *such that* $L_i$ *is a Moufang loop.*

*Note:* $L_i$ should not be a group it should be only a loop.

Likewise we can define the notion of weak Bol N-loop, weak Bruck N-loop and weak WIP N-loop. It is very important to note that the Lagrange theorem for finite groups can by no means be true for finite N-loops. So we make a small modification and define a new notion called Lagrange N-loop.

**DEFINITION 3.1.10:** *Let L = {$L_1 \cup L_2 \cup \ldots \cup L_N$, $*_1, \ldots, *_N$} be a N-loop of finite order. We call L a Lagrange N-loop if the order of every sub N-loop divides the order of L, then we call L a Lagrange N-loop.*

It is interesting to note that in case of N-loops one cannot always easily find out whether a N-loop is a Lagrange loop or not. So we now proceed on to define a new concept called weak Lagrange N-loop.

**DEFINITION 3.1.11:** *Let L = {$L_1 \cup L_2 \cup \ldots \cup L_N$, $*_1, \ldots, *_N$} be a N-loop of finite order. Suppose L has atleast one sub N-loop P = {$P_1 \cup P_2 \cup \ldots \cup P_N$, $*_1, \ldots, *_N$} such that order of P divides order of L then we call L a weak Lagrange N-loop.*

It is easy to verify that if a finite N-loop is a Lagrange N-loop then L is a weak Lagrange N-loop. However a weak Lagrange N-loop in general is not a Lagrange N-loop.

We give some examples to illustrate this.

***Example 3.1.3:*** Let L = {$L_1 \cup L_2 \cup L_3$, $*_1$, $*_2$, $*_3$} be a 3-loop where $L_1 = \{g \mid g^{12} = 1\}$, $L_2 = S_3$ the symmetric group of degree 3 and $L_3$ is given by the following table which is a loop of order 8.

| $*_3$ | e | $a_1$ | $a_2$ | $a_3$ | $a_4$ | $a_5$ | $a_6$ | $a_7$ |
|---|---|---|---|---|---|---|---|---|
| e | e | $a_1$ | $a_2$ | $a_3$ | $a_4$ | $a_5$ | $a_8$ | $a_7$ |
| $a_1$ | $a_1$ | e | $a_5$ | $a_2$ | $a_6$ | $a_3$ | $a_7$ | $a_4$ |
| $a_2$ | $a_2$ | $a_5$ | e | $a_5$ | $a_3$ | $a_7$ | $a_4$ | $a_1$ |
| $a_3$ | $a_3$ | $a_2$ | $a_6$ | e | $a_7$ | $a_4$ | $a_1$ | $a_5$ |
| $a_4$ | $a_4$ | $a_6$ | $a_3$ | $a_7$ | e | $a_1$ | $a_5$ | $a_2$ |
| $a_5$ | $a_5$ | $a_3$ | $a_7$ | $a_4$ | $a_1$ | e | $a_2$ | $a_6$ |
| $a_6$ | $a_6$ | $a_7$ | $a_4$ | $a_1$ | $a_5$ | $a_2$ | e | $a_3$ |
| $a_7$ | $a_7$ | $a_4$ | $a_1$ | $a_5$ | $a_2$ | $a_6$ | $a_3$ | e |

Clearly L is a 3-loop of order 26.

Now consider the sub 3 loop given by $P_1 \cup P_2 \cup P_3$ where

$$P_1 = \{g^2, g^4, g^6, g^8, g^{10}, 1\},$$

$$P_2 = \left\{\begin{pmatrix}1 & 2 & 3\\1 & 2 & 3\end{pmatrix}, \begin{pmatrix}1 & 2 & 3\\2 & 3 & 1\end{pmatrix}, \begin{pmatrix}1 & 2 & 3\\3 & 1 & 2\end{pmatrix}\right\} \text{ and}$$

$$P_3 = \{e, a_1\}.$$

Now order of P = 6 + 3 + 2 = 11 and $11 \nmid 26$.
Now consider $K = K_1 \cup K_2 \cup K_3$ where

$$K_1 = \{1, g^2, g^4, g^6, g^8, g^{10}\}, K_2 = \{S_3\} \text{ and } K_3 = \{e\}.$$

Thus order of K is 6 + 6 + 1 = 13. Now 13/ 26 so L is a weak Lagrange 3 loop and not a Lagrange 3 loop.

We can construct several such examples. We proceed on to define the notion of p-Sylow N-loops.

**DEFINITION 3.1.12:** *Let $L = \{L_1 \cup L_2 \cup \ldots \cup L_N, *_1, \ldots, *_N\}$ be a N-loop of finite order. If p is a prime such that $p^\alpha / o(L)$ but $p^{\alpha+1} \nmid o(L)$ and if the N-loop L has a sub N-loop P of order $p^\alpha$ but not of order $p^{\alpha+1}$ then we call P the p-Sylow sub N-loop of the N-loop L.*

*If for every prime p we have a p-Sylow sub N-loop then we call L a Sylow N-loop.*

First we prove that the Sylow theorem cannot in general hold good for all N-loops L.

**THEOREM 3.1.2:** *Let $(L = L_1 \cup L_2 \cup \ldots \cup L_N, *_1, \ldots, *_N)$ be a N-loop of finite order. If p is a prime such that $p^\alpha / o(L)$ but $p^{\alpha+1} \nmid o(L)$ then the N-loop in general does not have a sub N-loop of order $p^\alpha$.*

*Proof:* We show that we can have N-loops of finite order and if $p^{\alpha} / o(L)$ and $p^{\alpha+1} \not| \ o(L)$ still we can have sub N-loops of order $p^{\alpha+1}$. This is true from example 3.1.3, 2 / 26, L has no sub 3-loop of order 2.

***Example 3.1.4:*** Consider the 3-loop L = $\{L_1 \cup L_2 \cup L_3, *_1, *_2, *_3\}$ where $L_1 = \{g \mid g^{14} = 1\}$ a cyclic group of order 14, $L_2 = \{S_3\}$ the symmetric group of degree 3 and $L_3$ the loop of order 8 given in page 22. The order of the N-loop L is given by 14 + 6 + 8 = 28, clearly $2^2 / 28$, $2^3 \not| \ 28$ but L has a sub 3-loop or order 8. Take $K = K_1 \cup K_2 \cup K_3$ where $K_1 = \{1\}$, $K_2 = S_3$ and $K_3 = \{e\}$. Clearly o(K) = 8, $8 \not| \ 28$ but L has a sub 3-loop of order 8. Hence the claim.

However the N-loop can have sub N-loops of order $2^2$. For consider $H = H_1 \cup H_2 \cup H_3$ where

$$H_1 = \{1\},\ \ H_2 = \left\{\begin{pmatrix} 1 & 2 & 3 \\ 1 & 2 & 3 \end{pmatrix}\right\} \text{ and } H_3 = \{e, a_3\}.$$

Clearly order of H is 4 and 4/28. Also it is easily verified that this 3-loop L has a sub 3-loop of order 7. Take $T = T_1 \cup T_2 \cup T_3$ where

$$T_1 = \{1, g^7\},\ T_2 = \left\{\begin{pmatrix} 1 & 2 & 3 \\ 1 & 2 & 2 \end{pmatrix}, \begin{pmatrix} 1 & 2 & 3 \\ 2 & 3 & 1 \end{pmatrix}, \begin{pmatrix} 1 & 2 & 3 \\ 3 & 1 & 2 \end{pmatrix}\right\}$$

$$\text{and } T_3 = \{e, a_5\}.$$

o(T) = 2 + 3 + 2 = 7, 7 / 28. L cannot have sub 3- loop of order 49. Now comes a confusion to one that should one call the 3-loop L a 2-Sylow sub 3-loop and 7-Sylow sub 3-loop.

So we make a new definition called the super p-Sylow N-loop.

**DEFINITION 3.1.13:** *Let L be a N-loop of finite order. Suppose p is a prime such that $p^{\alpha} / o(L)$ but $p^{\alpha+1} \not| \ o(L)$ but L has a sub N-* *loop of order $p^{\alpha+1}$ (or $p^{\alpha+r}$ $r \geq 1$) then we call L a super p-Sylow N-loop or just super Sylow N-loop provided L is also a Sylow N-loop.*

The 3-loop described in example is a super p-Sylow 3-loop.

*Note:* By the very definition we see a super Sylow N-loop is always a Sylow N-loop but in general a Sylow N-loop is not a super Sylow N-loop. Now one may have the following situation. Let L be a N-loop of finite order. Suppose p is a prime such that $p^{\alpha} / o(L)$ but $p^{\alpha+1} \not| \ o(L)$ and L has only sub-N-loops of order less than $p^{\alpha}$ say of order $p^r$ ($r < \alpha$) then we are forced to define the notion of weak Sylow N-loop.

**DEFINITION 3.1.14:** *Let $L = \{L_1 \cup L_2 \cup \ldots \cup L_N, *_1, \ldots, *_N\}$ be a N-loop of finite order. Suppose p is a prime such that $p^{\alpha} / o(L)$ but $p^{\alpha+1} \not| \ o(L)$ and L has only sub N-loop of order $p^r$ ($r < \alpha$) this is true for all primes p satisfying the condition $p^{\alpha} / o(L)$ but $p^{\alpha+1} \not| \ o(L)$ then we call the N-loop L to be a weak p-Sylow N-loop.*

We illustrate this by the following example:

***Example 3.1.5:*** Let $L = \{L_1 \cup L_2 \cup L_3, *_1, *_2, *_3\}$ be a 3-loop, with

$L_1$ = $\{g \mid g^7 = 1\}$ is a cyclic group of order 7,
$L_2$ = $S_3$ the symmetric group of order 3 and
$L_3$ = the loop of order 8 given in page 22.

$o(L) = 7 + 6 + 8 = 21$, $3 / 21$, $3^2 \not| \ 21$. But L has no non trivial sub 3-loop of order 3. Clearly L has a sub 3-loop of order 9 given by $K = K_1 \cup K_2 \cup K_3$ where

$$K_1 = \{1\}, K_2 = S_3 \text{ and } K_3 = \{e, a_2\}.$$

Order of K = 1 + 6 + 2 = 9. Thus we cannot claim L to be a super Sylow 3-loop.

***Example 3.1.6:*** Take L to be 4- loop where $L = L_1 \cup L_2 \cup L_3 \cup L_4$ where

$L_1$ = $\{g \mid g^7 = 1\}$
$L_2$ = $D_{2.5} = \{a, b \mid a^2 = b^3 = 1; bab = a\}$,
$L_3$ = {loop of order 8} and
$L_4$ = $S_3$ .

o(L) = 7 + 9 + 8 + 6 = 30, 5/30 and $5^2 \nmid 30$. But L has a sub 4-loop of order 25.

For consider $K = K_1 \cup K_2 \cup K_3 \cup K_4$ where

$$K_1 = L_1,\ K_2 = L_2,\ K_3 = L_3 \text{ and } K_4 = \left\{\begin{pmatrix} 1 & 2 & 3 \\ 1 & 2 & 3 \end{pmatrix}\right\}.$$

Clearly order of K is 25 still L is not a super Sylow N-loop as L has no sub 4-loop of order 3.

***Example 3.1.7:*** Consider the 3-loop $L = L_1 \cup L_2 \cup L_3$, where

$L_1$ = $\{g \mid g^{37} = 1\}$,
$L_2$ = loop of order 8 given in page 68 and
$L_3$ = $A_3$, the alternating subgroup of $S_3$

o(L) = 37 + 8 + 3 = 48, $2^3/ 48$, $2^4 \nmid 48$, L has no sub 3 loop of order 8 but L has a sub 3-loop of order $2^2 = 4$.

Take $K = K_1 \cup K_2 \cup K_3$ where

$$K_1 = \langle 1 \rangle,\ K_2 = \{e, a_2\} \text{ and } K_3 = \left\{\begin{pmatrix} 1 & 2 & 3 \\ 1 & 2 & 3 \end{pmatrix}\right\}.$$

o(K) = 4. Hence the claim. 3 / 48, $3^2 \nmid 48$, so L has no subloop of order less than 3. Thus L is a weak 2-Sylow N-loop.

**DEFINITION 3.1.15:** *Let $L = \{L_1 \cup L_2 \cup \ldots \cup L_N, *_1, \ldots, *_N\}$ be a N-loop. If x and y $\in$ L are elements of $L_i$ the N-commutator (x, y) is defined as xy = (yx) (x, y), $1 \leq i \leq N$.*

It is important note if $x \in L_i$ and $y \in L_j$, $i \neq j$ such that $x \notin L_j$ and $y \notin L_i$ then we cannot define N-commutator. So even the very existence of the definition of N-commutator is not always guaranteed.

Similarly the definition of associator also is defined only when the triple chosen is in one of the $L_i$'s.

**DEFINITION 3.1.16:** *Let $L = \{L_1 \cup L_2 \cup \ldots \cup L_N, *_1, \ldots, *_N\}$ be a N-loop. If x, y, z are elements of the N-loop L, an associator (x, y, z) is defined only if x, y, z $\in L_i$ for some i ($1 \leq i \leq N$) and is defined to be (xy) z = (x (y z)) (x, y, z).*

Likewise we can define the notion of left nucleus, middle nucleus and right nucleus of a N-loop, which is left as an exercise for the reader.

Now we proceed on to define the notion of right regular N-representation of a finite N-loop L.

**DEFINITION 3.1.17:** *Let $L = \{L_1 \cup L_2 \cup \ldots \cup L_N, *_1, *_2, \ldots, *_N\}$ be a N-loop of finite order. For $\alpha_i \in L_i$ define $R_{\alpha_i}$ as a permutation of the loop $L_i$, $R_{\alpha_i} : x_i \to x_i\alpha_i$. This is true for i = 1, 2,…, N we define the set*

$$\left\{R_{\alpha_1} \cup R_{\alpha_2} \cup \ldots \cup R_{\alpha_N} \,\middle|\, \alpha_i \in L_i; i = 1, 2, \ldots, N\right\}$$

*as the right regular N-representation of the N loop L.*

Now we proceed on to define the notion of principal isotope of a N-loop L.

**DEFINITION 3.1.18:** *Let $L = \{L_1 \cup L_2 \cup \ldots \cup L_N, *_1, \ldots, *_N\}$ be a N-loop. For any pre determined pair $a_i, b_i \in L_i$, $i \in \{1, 2, \ldots, N\}$ a principal isotope $(L, o_1, \ldots, o_N)$, of the N loop L is defined by $x_i\ o_i\ y_i = X_i *_i Y_i$ where $X_i + a_i = x_i$ and $b_i + Y_i = y_i$, $i = 1, 2,.., N$. L is called G-N-loop if it is isomorphic to all of its principal isotopes.*

Obtain some interesting results about principal isotopes. Can one get any analogue results for N-loops?

Is the principal isotope of a N-loop a N-loop?

## 3.2 Smarandache N-loops and their properties

Now we in this section proceed on to define the notion of Smarandache N-loop and list a few of its properties.

**DEFINITION 3.2.1:** *Let $L = \{L_1 \cup L_2 \cup \ldots \cup L_N, *_1, \ldots, *_N\}$ be a N-loop. We call L a Smarandache N-loop (S-N-loop) if L has a proper subset P; $(P = P_1 \cup P_2 \cup \ldots \cup P_N, *_1, \ldots, *_N)$ where P is a N-group.*

***Example 3.2.1:*** Let $L = \{L_1 \cup L_2 \cup L_3 \cup L_4, *_1, *_2, *_3, *_4\}$ be a 4-loop. $L_1$ is given by the table

| $*_1$ | e | $a_1$ | $a_2$ | $a_3$ | $a_4$ | $a_5$ |
|---|---|---|---|---|---|---|
| e | e | $a_1$ | $a_2$ | $a_3$ | $a_4$ | $a_5$ |
| $a_1$ | $a_1$ | e | $a_3$ | $a_5$ | $a_2$ | $a_4$ |
| $a_2$ | $a_2$ | $a_5$ | e | $a_4$ | $a_1$ | $a_3$ |
| $a_3$ | $a_3$ | $a_4$ | $a_1$ | e | $a_5$ | $a_2$ |
| $a_4$ | $a_4$ | $a_3$ | $a_5$ | $a_2$ | e | $a_1$ |
| $a_5$ | $a_5$ | $a_2$ | $a_4$ | $a_1$ | $a_3$ | e |

$L_2 = \{g \mid g^{10} = 1\}$, $L_3$ is the loop of order 8 given in page 68 of this book and $L_4 = S_3$. Take $P = P_1 \cup P_2 \cup P_3 \cup P_4$. $P_1 = \{e, a_3\}$, $P_2 = \{1, g^2, g^4, g^6, g^8\}$, $P_3 = \{1, g_3\}$ and $P_4 = A_3$.

Clearly P is a 4-group. So L is a Smarandache 4-loop.

Let L be a S-N-loop. L is said to be of finite order if L has finite number of distinct elements.

If L has infinite number of elements then L is called as a Smarandache N-loop of infinite order.

Now we proceed on to define the notion of Smarandache sub N-loop.

**DEFINITION 3.2.2:** *Let $L = \{L_1 \cup L_2 \cup \ldots \cup L_N, *_1, \ldots, *_N\}$ be a N-loop. We call a proper subset P of L where $P = \{P_1 \cup P_2 \cup \ldots \cup P_N, *_1, \ldots, *_N\}$ to be a Smarandache sub N-loop (S-sub N-loop) of L if P itself is a S-N-loop.*

The following theorem is straightforward.

**THEOREM 3.2.1:** *Let $L = \{L_1 \cup L_2 \cup \ldots \cup L_N, *_1, \ldots, *_N\}$ be a N-loop. Suppose L has a S-sub N-loop then L is a S-N-loop.*

*Proof:* Let L be a N-loop. Let P be a proper subset of L. We are given P is a S-sub N-loop of L. So P has a proper subset T such that T is a N-group.

Now $T \subset L$ and T is a proper subset of L and T is a N-group in L; hence L is a S- N-loop. Hence the claim.

Now we proceed on to define to notion of Smarandache commutative N-loop and Smarandache Weakly commutative N-loop.

**DEFINITION 3.2.3**: *Let $L = \{L_1 \cup L_2 \cup \ldots \cup L_N, *_1, \ldots, *_N\}$ be a N-loop. L is said to be a Smarandache commutative N-loop (S-commutative N-loop) if every proper subset P of L which are N-groups are commutative N-groups of L.*

*If the N-loop L has atleast one proper subset P which is a commutative N-group then we call L a Smarandache weakly commutative N-loop (S-weakly commutative loop).*

**DEFINITION 3.2.4:** *Let $L = \{L_1 \cup L_2 \cup \ldots \cup L_N, *_1, \ldots, *_N\}$ be a N-loop. We say the N-loop L is a Smarandache cyclic N-loop*

*(S-cyclic N-loop) if every proper subset which is a N-group is a cyclic N-group.*

If the N-loop L has atleast one proper subset which is a N-group which is cyclic then we call L a Smarandache weakly cyclic N-loop (S-weakly cyclic N-loop).

**DEFINITION 3.2.5:** *Let $L = \{L_1 \cup L_2 \cup \ldots \cup L_N, *_1, \ldots, *_N\}$ be a N-loop. A proper S-sub N-loop $(P = P_1 \cup P_2 \cup \ldots \cup P_N, *_1, \ldots, *_N)$ of L is said to be a Smarandache normal N-loop (S-normal N-loop) if*

i. $x_i P_i = P_i x_i$
ii. $P_i x_i(y_i) = P_i (x_i y_i)$
iii. $y_i (x_i P_i) = (y_i x_i) P_i$

*for all $x_i, y_i \in P_i$ for $i = 1, 2, ..., N$.*

*If the N-loop L has no proper S-normal sub N-loop then we call the N-loop to be Smarandache simple N-loop (S-simple N-loop).*

It is important to note the notions S-simple N-loop and simple N-loop may or may not have in general any proper relations.

**DEFINITION 3.2.6:** *Let $L = \{L_1 \cup L_2 \cup \ldots \cup L_N, *_1, \ldots, *_N\}$ be a N-loop and $P = (P_1 \cup P_2 \cup \ldots \cup P_N, *_1, \ldots, *_N)$ be a S-sub N-loop of L. $(P \subset L)$ for any N-pair of elements $x_i, y_i \in P_i$ $(i \in \{1, 2, \ldots, N\})$ the N-commutator $(x_i, y_i)$ is defined by $x_i y_i = (y_i x_i)(x_i, y_i)$.*

*The Smarandache commutator sub N-loop (S-commutator sub N-loop) of L denoted by $S(L^S)$ is the S-sub N-loop generated by all its commutators.*

If this set does not generate a S-sub N-loop we then say the S-commutator subloop generated relative to the S-sub N-loop P of L is empty.

Thus it is important to note that in case of S-commutators in N-loops we see there can be more than one S-commutator S-sub N-loop, which is a marked difference between the commutator in N-loops and S-commutator in N-loops.

Now we proceed onto define the S-associators of N-loops.

**DEFINITION 3.2.7:** *Let $L = \{L_1 \cup L_2 \cup \ldots \cup L_N, *_1, \ldots, *_N\}$ be a N-loop. Let $P \subset L$ be a S-sub N-loop of L. If x, y, z are elements of the S-sub N-loop P of L, an associator (x, y, z) is defined by $(xy)z = (x\,(y\,z))\,(x, y, z)$.*

*The associator S-sub N-loop of the S-sub N-loop P of L denoted by $A\left(L_N^S\right)$ is the S-sub N-loop generated by all the associators, that is $\langle \{x \in P \mid x = (a, b, c) \text{ for some } a, b, c \in P\}\rangle$. If $A\left(L_N^S\right)$ happens to be a S-sub N-loop then only we call $A\left(L_N^S\right)$ the Smarandache associator (S-associator) of L related to the S-sub N-loop P.*

Now we proceed on to define the notion of pseudo commutative N- loop.

**DEFINITION 3.2.8:** *Let $L = \{L_1 \cup L_2 \cup \ldots \cup L_N, *_1, \ldots, *_N\}$ be a N-loop. If for a, $b \in L$ with ab = ba in L, we have $(ax_i)\, b = (bx_i)\, a$ (or $b\,(x_i\, a)$) for all $x_i \in L_i$ (if a, $b \in L_i$ ), then we say the pair is pseudo commutative.*

*If every commutative pair in every $L_i$ is pseudo commutative then we call L a pseudo commutative N-loop. If we have for every pair of commuting elements a, b in a proper subset P of the N-loop L (which is a S-sub N-loop). If a, b is a pseudo commutative pair for every x in P then we call the pair a Smarandache pseudo commutative pair (S-pseudo commutative pair).*

*Note:* We see a pair is a pseudo commutative pair then it is obviously a S-pseudo commutative pair. But we see on the contrary that all S-pseudo commutative pairs in general need not be a pseudo commutative pairs. If in a S-sub N-loop P of L

every pair in P is a S-pseudo commutative pair then we call the N-loop L a Smarandache pseudo commutative N-loop (S-pseudo commutative N-loop).

Now we proceed on to define the notion of Smarandache pseudo commutator.

**DEFINITION 3.2.9:** *Let L = {$L_1 \cup L_2 \cup \ldots \cup L_N$, $*_1, \ldots, *_N$} be a N-loop. The pseudo commutator of L denoted by P(L) = $\langle$ {p $\in$ L / a(xb) = p ([bx] a)}$\rangle$; we define the Smarandache pseudo commutator (S-pseudo commutator) of P $\subset$ L (P a S-sub-N-loop of L) to be the S-sub N-loop generated by $\langle$ {p $\in$ P | a (xb) = p ( [b x] a), a, b $\in$ P}$\rangle$; denoted by $P_S^N$ (L). If $P_S^N$ (L) is not a S-sub N loop then we say $P_S^N$ (L) = $\phi$.*

It is important to note that we may have more than one S-sub N loop $P_S^N$ (L) for a given N-loop L.

Now we proceed on to define the notion of Smarandache pseudo associative triple in case of a N-loop L.

**DEFINITION 3.2.10:** *Let L = {$L_1 \cup L_2 \cup \ldots \cup L_N$, $*_1, \ldots, *_N$} be a N-loop. An associative triple a, b, c $\in$ P $\subset$ L where P is a S-sub N- loop of L is said to be Smarandache pseudo associative (S-pseudo associative) if (a b) (x c) = (a x) (b c) for all x $\in$ P. If (a b) (x c) = (a x) (b c) for some x $\in$ P we say the triple is Smarandache pseudo associative (S-pseudo associative) relative to those x in P.*

*If in particular every associative triple in P is S-pseudo associative then we say the N-loop is a Smarandache pseudo associative N-loop (S- pseudo associative N-loop).*

*Thus for a N-loop to be a S-pseudo associative N-loop it is sufficient that one of its S-sub N-loops are a S-pseudo associative N-loop.*

*The Smarandache pseudo associator of a N-loop (S-pseudo associator of a N-loop) L denoted by $S\left(AL_p^N\right)$ = $\langle$ {t $\in$ P / (a b) (t c) = (at) (bc) where a (bc) = (a b) c for a, b, c $\in$ P}$\rangle$ where P is a S-sub N-loop of L.*

**DEFINITION 3.2.11:** *Let $L = \{L_1 \cup L_2 \cup \ldots \cup L_N, *_1, \ldots, *_N\}$ be a N-loop. L is said to be a Smarandache Moufang N loop (S-Moufang N-loop) if there exists S-sub N-loop P of L which is a Moufang N-loop.*

*A N-loop L is said to be a Smarandache Bruck N-loop (S-Bruck N-loop) if L has a proper S-sub N-loop P, where P is a Bruck N-loop.*

*Similarly a N-loop L is said to be a Smarandache Bol N-loop (S-Bol N-loop) if L has a proper S-sub N loop) K, where K is a Bol N-loop.*

*We call a N-loop L to be a Smarandache Jordan N-loop (S-Jordan N-loop) if L has a S-sub N-loop T, where T is a Jordan N-loop.*

On similar lines we define Smarandache right (left) alternative N-loop.

Now proceed on to define the Smarandache N-loop homomorphism.

**DEFINITION 3.2.12:** *Let $L = \{L_1 \cup L_2 \cup \ldots \cup L_N, *_1, \ldots, *_N\}$ and $K = \{ K_1 \cup K_2 \cup \ldots \cup K_N, o_1, o_2, \ldots, o_N\}$ be two N-loops. We say a map $\phi = \phi_1 \cup \phi_2 \cup \ldots \cup \phi_N$ from L to K is a Smarandache N-loop homomorphism (S-N loop homomorphism) if $\phi$ is a N-group homomorphism from P to T where P and T are proper subsets of L and K respectively such that P and T are N-groups of L and K.*

Thus for Smarandache $\phi$ - homomorphism of N-loops we need not have $\phi$ to be defined on the whole of L. It is sufficient if $\phi$ is defined on a proper subset P of L which is a N-group.

Now we proceed on to define Smarandache Moufang center of a N loop.

**DEFINITION 3.2.13:** *Let $L = \{L_1 \cup L_2 \cup \ldots \cup L_N, *_1, \ldots, *_N\}$be a N-loop. Let P be a S-sub N-loop of L. The Smarandache Moufang Centre (S-Moufang Centre) of L is defined to be the Moufang centre of the S-sub N-loop; $P = (P_1 \cup P_2 \cup \ldots \cup P_N, *_1, \ldots, *_N)$.*

Thus for a given N-loop L we can have several S-Moufang centres or no Moufang centre if the N-loop L has no S-Sub N-loop.

Now we proceed on to define the notion of S-center.

**DEFINITION 3.2.14:** *Let L be a N-loop. P be a S-sub N-loop of L. The Smarandache N-centre (S N center) $SZ_N$ (L) is the center of the sub-N-loop $P \subset L$.*

Thus even in case of S-centre for a N-loop we may have several S-centers depending on the number of S-sub N-loops.

Now notion of Smarandache middle, left and right nucleus of a N-loop is defined.

**DEFINITION 3.2.15:** *Let L be a N-loop, P be a S-sub N-loop of L. To define Smarandache nucleus (S-nucleus) related to P we consider the Smarandache left nucleus (S-left nucleus) of L to be the left nucleus of the S-sub N-loop P denoted by $N_{\lambda}^{P_N}$.*

*Similarly the Smarandache right nucleus (S-right nucleus) of the N-loop L to be the right nucleus of the S-sub N-loop P denoted by $N_{P}^{P_N}$.*

*We can define the Smarandache middle nucleus (S-middle nucleus) of the N-loop L to be the middle nucleus of the S-sub N-loop P, denoted by $N_{m}^{P_N}$. Thus the Smarandache nucleus (S-nucleus) is $[SN(L)]_N = N_{\lambda}^{P_N} \cap N_{P}^{P_N} \cap N_{M}^{P_N}$.*

*Note:* We define left, right and middle S- nucleus only for the same S-sub N-loop P. As we vary P we can have several S-nucleus.

Now we proceed on to define the notion of Smarandache Lagrange N-loop, Smarandache Cauchy element in a N-loop and Smarandache p-Sylow N-loop.

**DEFINITION 3.2.16:** *Let $L = \{L_1 \cup L_2 \cup \ldots \cup L_N, *_1, \ldots, *_N\}$be a N-loop of finite order. A proper subset S of L where S is a S-sub-N-loop is said to be a Smarandache Lagrange sub N-loop (S-Lagrange sub N-loop) of L if o (S) / o (L).*

*If every S-sub N-loop of the N-loop L happens to be a S-Lagrange sub N-loop then we call L a Smarandache Lagrange N-loop (S-Lagrange N-loop). If no S-sub N-loop of L is a S-Lagrange sub N-loop then we call L a Smarandache non Lagrange N-loop (S-non Lagrange N-loop).*

Now we proceed on to define the notion of Smarandache Cauchy element of a N-loop L.

**DEFINITION 3.2.17:** *Let $L = \{L_1 \cup L_2 \cup \ldots \cup L_N, *_1, \ldots, *_N\}$be a N-loop of finite order. An element $x \in P$, where $P = \{P_1 \cup P_2 \cup \ldots \cup P_N, *_1, \ldots, *_N\}$ is a S-sub N-loop of L is said to be a Smarandache Cauchy element (S-Cauchy element) of the N-loop L if $x^n = e$, $(n > 1)$ and $n / o(P)$.*

*If every S-sub N-loop of L has a Cauchy element then we call the N-loop L to be a Smarandache Cauchy N-loop (S-Cauchy N-loop).*

We now proceed to define Smarandache p-Sylow N-loop.

**DEFINITION 3.2.18:** *Let $L = \{L_1 \cup L_2 \cup \ldots \cup L_N, *_1, \ldots, *_N\}$be a N-loop of finite order. Let p be a prime such that $p^\alpha / o(L)$ but $p^{\alpha+1} \not| \; o(L)$ and if L has a S-sub N-loop P of order $p^\alpha$ then we call P the Smarandache p-Sylow sub N-loop (S-p-Sylow sub N-loop). If for every prime p, $p^\alpha / o(L)$ and $p^{\alpha+1} \not| \; o(L)$; we have S-p-Sylow sub N-loop then we call the N-loop L to be a Smarandache Sylow N-loop (S-Sylow N-loop).*

We proceed on to define super Smarandache p-Sylow sub N-loop and the notion of super Smarandache Sylow N-loop.

**DEFINITION 3.2.19:** *Let $L = \{L_1 \cup L_2 \cup \ldots \cup L_N, *_1, \ldots, *_N\}$ be a N-loop of finite order. Suppose L is a S-Sylow N-loop. If for every prime p, $p^\alpha / o(L)$ but $p^{\alpha+1} \not| \; o(L)$ we have a S-sub N loop $P_1$ of order $p^{\alpha+r}$ $(r \geq 1)$ then we call L a (S-super Smarandache Sylow N-loop). $P_1$ is called as a super Smarandache p-Sylow S-sub N loop (super p-Sylow S-sub N loop) of L.*

*We define a N-loop L to be a Smarandache pseudo N-loop (S pseudo N-loop) if* $L = \{L_1 \cup L_2 \cup \ldots \cup L_N, *_1, \ldots , *_N\}$ *contains a non empty subset P where* $P = \{P_1 \cup P_2 \cup \ldots \cup P_N, *_1, \ldots, *_N\}$ *is a N-groupoid.*

Several results can be derived studied and analyzed using these definitions. This work is left for the researchers.

**Chapter Four**

# N-GROUPOIDS AND SMARANDACHE N-GROUPOIDS

In this chapter we for the first time introduce the notion of N-groupoids and Smarandache N-groupoids. The notion of bigroupoids was introduced in 2003 and also in the same year the notion of Smarandache bigroupoids was introduced.

We here introduce these N-groupoids (N > 2) for when N = 2 we get the usual bigroupoid and also give some of its applications.

This chapter has four sections. In section one we just introduce the concept of bigroupoids and S-bigroupoids. In section two we introduce the new notion of N-groupoid and give some of its properties. Section three gives the definition and properties of Smarandache-N-groupoids. In the final section we suggest some of the probable applications.

## 4.1 Introduction to bigroupoids and Smarandache bigroupoids

In this section we just recall the definition and some of the basic properties of bigroupoids and S-bigroupoids and illustrate them with examples.

**DEFINITION 4.1.1:** *Let G = (G, +, o) be a non empty set. We call G a bigroupoid if* $G = G_1 \cup G_2$ *where* $G_1$ *and* $G_2$ *are proper subsets of G and satisfy the following:*

i. *($G_1$, +) is a groupoid (i.e. '+' is non associative closed binary operation).*
ii. *($G_2$, o) is a semigroup.*

***Example 4.1.1:*** Let G = ($G_1 \cup G_2$, *, o) where $G_1$ is given by the table

| * | $h_1$ | $h_2$ | $h_3$ |
|---|---|---|---|
| $h_1$ | $h_1$ | $h_3$ | $h_2$ |
| $h_2$ | $h_2$ | $h_1$ | $h_3$ |
| $h_3$ | $h_3$ | $h_2$ | $h_1$ |

and $G_2$ = S(4), the semigroup of set of all mapping of (1, 2, 3, 4) to itself under composition of mappings. G is a bigroupoid.

The number of distinct elements in a bigroupoid G is called the order of the bigroupoid and is denoted by o (G). If o (G) = n and $n < \infty$ implies the bigroupoid G is of *finite order*. If $n = \infty$ we say G is a bigroupoid of *infinite order*.

Further it is important to note that in general $o(G) \neq o(G_1) + o(G_2)$ for we may have common elements in $G_1$ and $G_2$ but we can say

$$o(G) = o(G_1) + o(G_2) – o(G_1 \cap G_2).$$

Now we proceed on to define the notion of sub-bigroupoid of a bigroupoid G.

**DEFINITION 4.1.2:** *Let G = ($G_1 \cup G_2$, *, o) be a bigroupoid. We say a proper subset H of G is said to be a subbigroupoid if H under the operations of G is itself a bigroupoid.*

We illustrate this by the following example.

***Example 4.1.2:*** Let $G = G_1 \cup G_2$ be a bigroupoid where

$G_1$ = S(5) is a semigroup of mappings of (1 2 3 4 5) to itself under the composition of mappings and

$G_2$ is given by the following table:

| * | 0 | 1 | 2 | 3 | 4 | 5 | 6 | 7 |
|---|---|---|---|---|---|---|---|---|
| 0 | 0 | 6 | 4 | 2 | 0 | 6 | 4 | 2 |
| 1 | 2 | 0 | 6 | 4 | 2 | 0 | 6 | 4 |
| 2 | 4 | 2 | 0 | 6 | 4 | 2 | 0 | 6 |
| 3 | 6 | 4 | 2 | 0 | 6 | 4 | 2 | 0 |
| 4 | 0 | 6 | 4 | 2 | 0 | 6 | 4 | 2 |
| 5 | 2 | 0 | 6 | 4 | 2 | 0 | 6 | 4 |
| 6 | 4 | 2 | 0 | 6 | 4 | 2 | 0 | 6 |
| 7 | 6 | 4 | 2 | 0 | 6 | 4 | 2 | 0 |

Now clearly one can verify ($G_2$, *) is a groupoid and not a semigroup.

Now take $P = P_1 \cup P_2$ where $P_1 = \{0, 2, 4, 6\}$ and $P_2 = S'(3)$ i.e., we fix in the mappings of the set 1, 2, 3, 4, 5, 4 and 5 only 123 are mapped. This will form a semigroup isomorphic to S(3). Clearly P is a subbigroupoid of the bigroupoid G.

We now proceed on to define the notion of commutative bigroupoid.

**DEFINITION 4.1.3:** *A bigroupoid G = ($G_1 \cup G_2$, *, o) is said be a commutative bigroupoid if (G, *) is a commutative groupoid and ($G_2$, o) is a commutative semigroup.*

We can still define a new notion called idempotent bigroupoid.

**DEFINITION 4.1.4:** *A bigroupoid (G, *, o) is said to be an idempotent bigroupoid if ($G_1$, *) is an idempotent groupoid and ($G_2$, o) is an idempotent semigroup.*

We can also as in the case of semigroups define the notion of right and left ideals of a bigroupoid.

**DEFINITION 4.1.5:** *Let (G, *, o) be a bigroupoid where G = $G_1 \cup G_2$ with ($G_1$, *) a groupoid and ($G_2$, o) a semigroup. We say a nonempty subset P of G is said to be a right biideal of G if P =*

*$P_1 \cup P_2$ such that $P_1$ is a right ideal of the groupoid $G_1$ and $P_2$ is a right ideal of the semigroup $G_2$.*

*On similar lines we can define the notion of left biideal of the bigroupoid G. We say a nonempty subset P of the groupoid G to be a biideal of G if $P = P_1 \cup P_2$ where $P_1 \subset G_1$ and $P_2 \subset G_2$ with $P_1$ is an ideal of the groupoid $G_1$ and $P_2$ the ideal of the semigroup $G_2$.*

Now we proceed on to define the notion of normal subbigroupoid of a bigroupoid G.

**DEFINITION 4.1.6:** *Let $G = (G_1 \cup G_2, *, o)$ be a bigroupoid. A non empty subset P of G is said to be a normal subbigroupoid of the bigroupoid G if*

i. *$P = P_1 \cup P_2$ where P is a subbigroupoid of G.*
ii. *$a P_1 = P_1 a$ or $aP_2 = P_2 a$ according as $a \in P_1$ or $a \in P_2$.*
iii. *$P_1 (x y) = (P_1 x) y$ or $P_2 (xy) = (P_2 x) y$ according as $x, y \in P_1$ or $x, y \in P_2$.*
iv. *$y (x P_1) = (y x) P_1$ or $y (x P_2) = (x y) P_2$ according as $x, y \in P_1$, or $x, y \in P_2$ for all $x, y, a \in G$.*

*We call a bigroupoid simple if it has no nontrivial normal subbigroupoid.*

Next we define a new notion or condition under which a bigroupoid itself is normal.

**DEFINITION 4.1.7:** *Let (G, +, o) be a bigroupoid, $G = G_1 \cup G_2$ where $(G_1, +)$ is a groupoid and $(G_2, o)$ is a semigroup. We say G is a normal bigroupoid if*

i. *$xG = Gx$ (i.e. $x_1 G_1 = G_1 x_1$ if $x_1 \in G_1$ and $x_2 G_2 = G_2 x_2$ if $x_2 \in G_2$).*
ii. *$G (x y) = (G x) y$ (i.e. $G_1 (x y) = (G_1 x) y$ ; if $x, y \in G_1$ and $((G_2 x) y = G_2 (x y)$ if $x, y \in G_2$).*

*iii.* $y\ (x\ G) = (y\ x)\ G$ *i.e.,* $(y\ x)\ G_1 = y\ (x\ G_1)$, $x, y \in G_1$, *and* $(y\ x)\ G_2 = y\ (x\ G_2)$ *if* $(x, y \in G_2)$.

We proceed on to define some more concepts like biconjugate, bicenter, biconjugate pair, etc.

**DEFINITION 4.1.8:** *Let* $(G, *, o)$ *be a bigroupoid with* $G = (G_1 \cup G_2)$ *where* $(G_1, *)$ *is a groupoid and* $(G_2, o)$ *is a semigroup. Let H and K be two proper subbigroupoids of* $(G, *, o)$ *where* $H = H_1 \cup H_2$, $(H_1, *)$ *is a subgroupoid of* $(G_1, *)$ *and* $(H_2, o)$ *is a subsemigroup of* $(G_2, o)$.

*Similarly* $K = K_1 \cup K_2$ *with* $K \cap H = \phi$. *We say H is biconjugate with K if there exists* $x_1 \in H_1$ *and* $x_2 \in H_2$ *with* $x_1 K_1 = H_1$ *or* $(K_1 x_1 = H_1)$ *and* $x_2 K_2 = H_2$ *(or* $K_2 x_2 = H_2$*) 'or' in the mutually exclusive sense.*

**DEFINITION 4.1.9:** *Let* $(G, *, o)$ *be a bigroupoid. An element* $x \in G$ *is said to be a zero divisor if there exists* $y \in G$ *such that* $xy = yx = 0$.

**DEFINITION 4.1.10:** *Let* $(G = G_1 \cup G_2, *, o)$ *be a bigroupoid. The centre or bicentre of* $(G, *, o)$ *denoted by* $BC(G) = \{x \in G_1 \mid xa = ax \text{ for all } a \in G_1\} \cup \{y \in G_2 \mid yb = by \text{ for all } b \in G_2\}$ *i.e. why we call it as bicentre.*

**DEFINITION 4.1.11:** *Let* $(G = G_1 \cup G_2, *, o)$ *be a bigroupoid .* $G = G_1 \cup G_2$ *we say* $a, b \in G$ *is a biconjugate pair if* $a = b\ o\ x_1$, *(if* $a, b \in G_1$ *for some* $x_1 \in G_1$ *or* $a = x_1\ o\ b$*) and* $b = a\ o\ y$, *(or* $y_1\ o\ a$ *for some* $y_1 \in G_1$*) (if* $a, b \in G_2$ *we choose x and y from* $G_2$*). On similar lines we define right biconjugate only for* $a, b \in G_1$ *(or* $G_2$*) as* $a\ o\ x = b$ *and* $b\ o\ y = a$ *for* $x \in G_1$ $(x \in G_2)$.

**DEFINITION 4.1.12:** *Let* $(Y, *, o)$ *be a set with* $Y = Y_1 \cup Y_2$, $Y_1$ *and* $Y_2$ *are proper subsets of Y where* $(Y_1, *)$ *is a loop and* $(Y_2, o)$ *a groupoid, then we call* $(Y, *, o)$ *as a biquasi loop. A proper subset* $P \subset Y$ *is said to be a sub-biquasi loop if P itself under the operations of Y is a biquasi loop.*

Now we proceed on to define the notion of normal biquasi loop.

**DEFINITION 4.1.13:** *Let (X, +, o) be a biquasi loop $X = X_1 \cup X_2$. $(X_1, +)$ is a loop and $(X_2, o)$ is a groupoid. $P = P_1 \cup P_2$ a subbiquasi loop of X. We say P is a normal subbiquasi loop of X if $P_1$ is a normal subloop of $(X_1, +)$ and $P_2$ is an ideal of the groupoid $(X_2, o)$.*

Now we proceed on to define the notion of Lagrange biquasi loop.

**DEFINITION 4.1.14:** *Let (X, +, o) be a biquasi loop of finite order. If the order of every proper subbiquasi loop P of X divides the order of X, then we say X is a Lagrange biquasi loop. If the order of atleast one of the proper subbiquasi loop P of X divides the order of X then we say X is a weakly Lagrange biquasi loop. If the sub biquasi loop exists and none of its order divides the order X then X is said to be non-Lagrange biquasi loop.*

We just illustrate this by an example.

***Example 4.1.3:*** Consider the biquasi loop (Y, +, o) with $Y = Y_1 \cup Y_2$ where $Y_1$ a loop of order 8 given by the following table:

| * | e | $a_1$ | $a_2$ | $a_3$ | $a_4$ | $a_5$ | $a_6$ | $a_7$ |
|---|---|---|---|---|---|---|---|---|
| e | e | $a_1$ | $a_2$ | $a_3$ | $a_4$ | $a_5$ | $a_6$ | $a_7$ |
| $a_1$ | $a_1$ | e | $a_5$ | $a_2$ | $a_6$ | $a_3$ | $a_7$ | $a_4$ |
| $a_2$ | $a_2$ | $a_5$ | e | $a_6$ | $a_3$ | $a_7$ | $a_4$ | $a_1$ |
| $a_3$ | $a_3$ | $a_2$ | $a_6$ | e | $a_7$ | $a_4$ | $a_1$ | $a_5$ |
| $a_4$ | $a_4$ | $a_6$ | $a_3$ | $a_7$ | e | $a_1$ | $a_5$ | $a_2$ |
| $a_5$ | $a_5$ | $a_3$ | $a_7$ | $a_4$ | $a_1$ | e | $a_2$ | $a_6$ |
| $a_6$ | $a_6$ | $a_7$ | $a_4$ | $a_1$ | $a_5$ | $a_2$ | e | $a_3$ |
| $a_7$ | $a_7$ | $a_4$ | $a_1$ | $a_5$ | $a_2$ | $a_6$ | $a_3$ | e |

The groupoid $Y_2$ is given by the following table:

| + | $x_1$ | $x_2$ | $x_3$ |
|---|---|---|---|
| $x_1$ | $x_1$ | $x_3$ | $x_2$ |
| $x_2$ | $x_3$ | $x_2$ | $x_1$ |
| $x_3$ | $x_2$ | $x_1$ | $x_3$ |

Clearly o(Y) = 11. Y has no nontrivial subbiquasi loops but it has subbiquasi group. $P = P_1 \cup P_2$, where $P_1 = \{e, a_6\}$ is a group and $P_2 = \{x_1\}$ is a semigroup.

**DEFINITION 4.1.15:** *Let (Y, +, o) be a biquasi loop of finite order. If p is a prime such $p^\alpha$ / o(Y) and $p^{\alpha+1} \not| \; o(Y)$ we have a sub-biquasi loop P of order $p^\alpha$ then we say P is a p-Sylow sub-biquasi loop.*

*If for every prime p / o(Y) we have a p-Sylow subbiquasi loop then we call Y a Sylow biquasi loop.*

Just we proceed on to define the notion of Cauchy elements in a biquasi loop.

**DEFINITION 4.1.16:** *Let (Y, +, o) be a biquasi loop of finite order, an element $x \in Y$ such that $x^r = e$ (e-identity of Y) r> 1 with r / o(Y) is called the Cauchy element of Y.*

*We see in case of finite groups every element is a Cauchy element but in case of biquasi loop only a few of the elements may be Cauchy elements; if the order of a biquasi loop is a prime, p none of the elements in it would be a Cauchy element even though we will have $x \in Y$ with $x^t = e$, $t > 1$ $t \neq p$, as $t \not| \; p$.*

**DEFINITION 4.1.17:** *Let (P, +, o) be a non empty set. If (P, +) is a loop and (P, o) is a semigroup then we call P a biquasi semigroup.*

**THEOREM 4.1.1:** *All biquasi semigroups are biquasi loops.*

We now proceed on to define the notion of biquasi groupoid.

**DEFINITION 4.1.18:** *Let (T, +, o) be a non empty set. T is a biquasi groupoid if (T, +) is a group and (T, o) is a groupoid.*

**THEOREM 4.1.2:** *All biquasi groups are biquasi groupoids and not conversely.*

## 4.2 N-groupoids and their properties

In this section we proceed on to define the new notion of N-groupoids (N > 2) for when N = 2 we get the bigroupoids and give several of the properties and some other related structures.

**DEFINITION 4.2.1:** *Let* $G = (G_1 \cup G_2 \cup \ldots \cup G_N; *_1, \ldots, *_N)$ *be a non empty set with N-binary operations. G is called the N-groupoid if some of the* $G_i$*'s are groupoids and some of the* $G_j$*'s are semigroups,* $i \neq j$ *and* $G = G_1 \cup G_2 \cup \ldots \cup G_N$ *is the union of proper subsets of G.*

It is important to note that a $G_i$ is either a groupoid or a semigroup.

***Example 4.2.1:*** Let $G = (G_1 \cup G_2 \cup G_3 \cup G_4, *_1, *_2, *_3, *_4)$ where $G_1 = S(3)$, $G_2 = S(5)$, $G_3 = (Z_7, * \mid a * b = 3a + 2b \pmod 7)$ and $G_4 = (Z_4, a *_4 b = 3a + b \pmod 4$ for $a, b \in Z_4)$. Clearly G is 4-groupoid.

Now we proceed on to define the order of a N-groupoid G.

**DEFINITION 4.2.2:** *Let* $G = (G_1 \cup G_2 \cup \ldots \cup G_N; *_1, *_2, \ldots, *_N)$ *be a N-groupoid. The order of the N-groupoid G is the number of distinct elements in G. If the number of elements in G is finite we call G a finite N-groupoid. If the number of distinct elements in G is infinite we call G an infinite N-groupoid.*

The N-groupoid given in example 4.2.1 is a finite 4-gorupoid.

Now we proceed on to define the notion of commutative N-groupoid and idempotent N-groupoid.

**DEFINITION 4.2.3:** *Let $G = (G_1 \cup G_2 \cup \ldots \cup G_N; *_1, \ldots, *_N)$ be a N-groupoid. We say G is a commutative N-groupoid if each $(G_i, *_i)$ is either a commutative groupoid or a commutative semigroup.*

***Example 4.2.2:*** Let $G = \{G_1 \cup G_2 \cup G_3, *_1, *_2, *_3\}$ where

$G_1$ = $L_5(3)$ given by the table,
$G_2$ = $Z_{10}$ ={multiplicative semigroup under modulo 10} and
$G_3$ = {Z, under multiplication}.

G is a 3-groupoid which is commutative.

The table for the groupoid $L_5(3)$ is given which is easily checked to be commutative.

| * | e | 1 | 2 | 3 | 4 | 5 |
|---|---|---|---|---|---|---|
| e | e | 1 | 2 | 3 | 4 | 5 |
| 1 | 1 | e | 4 | 2 | 5 | 3 |
| 2 | 2 | 4 | e | 5 | 3 | 1 |
| 3 | 3 | 2 | 5 | e | 1 | 4 |
| 4 | 4 | 5 | 3 | 1 | e | 2 |
| 5 | 5 | 3 | 1 | 4 | 2 | e |

Next we proceed on to define the notion of sub N-groupoid of a N-groupoid G.

**DEFINITION 4.2.4:** *Let $G = \{G_1 \cup G_2 \cup \ldots \cup G_N, *_1, *_2, \ldots, *_N\}$ be a N-groupoid. We say a proper subset $H = \{H_1 \cup H_2 \cup \ldots \cup H_N, *_1, *_2, \ldots, *_N\}$ of G is said to be a sub N-groupoid of G if H itself is a N-groupoid of G.*

Now we illustrate this by the following example:

***Example 4.2.3:*** Let G = ($G_1 \cup G_2 \cup G_3$, $*_1$, $*_2$, $*_3$), be a 3-groupoid where

$G_1$ = {$Z_{12}$ / a * b = 3 a + 4 b (mod 12); a, b ∈ 12},
$G_2$ = {$Z_{16}$, semigroup under multiplication mod 16} and
$G_3$ = S(3).

Consider H = {$H_1 \cup H_2 \cup H_3$, $*_1$, $*_2$, $*_3$} where $H_1$ = {$G_1$}, $H_2$ = {0, 2, 4, 6, 8, 10, 12, 14, under multiplication modulo 16} and $H_3$ = {S(3)}. H is a sub N-groupoid of G.

It is interesting to note that as in case of all N-structures even in case of N-groupoids we may or may not have the order of the sub N-groupoid to divide the order of the N-groupoid G. For in the example 4.2.3 we see the order of the 3-groupoid is 55. The order of the sub 3-groupoid is 47, 47 $\nmid$ 55.

Thus we want conditions under which the sub N-groupoids of a finite N-groupoid divides the order of the N-groupoid. To this end define a new concept called Lagrange N- groupoid.

**DEFINITION 4.2.5:** *Let G = ($G_1 \cup G_2 \cup \ldots \cup G_N$, $*_1$, $*_2$, …, $*_N$) be a finite N-groupoid. If the order of every sub N-groupoid H divides the order of the N-groupoid, then we say G is a Lagrange N-groupoid.*

It is very important to mention here that in general every N-groupoid need not be a Lagrange N-groupoid.
Now we define still a weaker notion of a Lagrange N-groupoid.

**DEFINITION 4.2.6:** *Let G = ($G_1 \cup G_2 \cup \ldots \cup G_N$, $*_1$, …, $*_N$) be a finite N-groupoid. If G has atleast one nontrivial sub N-groupoid H such that o(H) / o(G) then we call G a weakly Lagrange N-groupoid.*

It is easily verified that every Lagrange N-groupoid is a weakly Lagrange N-groupoid. Also one cannot for certain say that all N-groupoids are weakly Lagrange N-groupoids. Thus we have the following containment relation.

{Lagrange N-groupoid} $\underset{\neq}{\subset}$ {Weakly Lagrange N-groupoid} $\underset{\neq}{\subset}$ {N-groupoid}.

Now we give an example of a weakly Lagrange N-groupoid.

***Example 4.2.4:*** Let G = ($G_1 \cup G_2 \cup G_3, \cup G_4, *_1, *_2, *_3, *_4$) be a 4-groupoid where

$G_1$ = ($Z_4$, semigroup under multiplication modulo 4),
$G_2$ = {$Z_6$ | a * b = 2a + b (mod 6); a, b $\in Z_6$},
$G_3$ = {$Z_7$, semigroup under multiplication modulo 7} and
$G_4$ = {$Z_{10}$ | a * b = 4a + 6b (mod 10) for a, b $\in Z_{10}$}.

Clearly G is a 4-groupoid and order of G is 27. Consider H = $H_1 \cup H_2 \cup H_3 \cup H_4$ where

$H_1$ = {0, 2, semigroup under multiplication modulo 4},
$H_2$ = {0, 2, 4},
$H_3$ = {1} and
$H_4$ = {0, 2, 4, 6, 8}.

Clearly o(H) = 11 and 11 $\not|$ 27.

Take K = $K_1 \cup K_2 \cup K_3 \cup K_4$ where $K_1$ = {0, 2}, $K_2$ = {3}, $K_3$ = {1} and $K_4$ = {0, 2, 4, 6, 8} where $K_i \subseteq G_i$, i = 1, 2, 3, 4. Clearly o(K) = 9 and 9 / 27. Thus G is a weakly Lagrange N-groupoid and not a Lagrange N-groupoid.

Now we proceed on to define the notion of normal sub N-groupoid.

**DEFINITION 4.2.7:** *Let G = ($G_1 \cup G_2 \cup \ldots \cup G_N, *_1, *_2, \ldots, *_N$) be a N groupoid. A sub N groupoid V = $V_1 \cup V_2 \cup \ldots \cup V_N$ of G is said to be a normal sub N-groupoid of G; if*

i. *a $V_i$ = $V_i$ a, i = 1, 2, …, N whenever a $\in V_i$*
ii. *$V_i$ (x y) = ($V_i$ x) y, i = 1, 2, …, N for x, y $\in V_i$*
iii. *y (x $V_i$ ) = (xy) $V_i$ , i = 1, 2, …, N for x, y $\in V_i$.*

*Now we say a N-groupoid G is simple if G has no nontrivial normal sub N-groupoids.*

Now we proceed on to define the notion of N-conjugate groupoids.

**DEFINITION 4.2.8:** *Let $G = (G_1 \cup G_2 \cup \ldots \cup G_N, *_1, \ldots, *_N)$ be a N-groupoid. Let $H = \{H_1 \cup \ldots \cup H_N; *_1, \ldots, *_N\}$ and $K = \{K_1 \cup K_2 \cup \ldots \cup K_N, *_1, \ldots, *_N\}$ be sub N-groupoids of $G = G_1 \cup G_2 \cup \ldots \cup G_N$; where $H_i$, $K_i$ are subgroupoids of $G_i$ ($i = 1, 2, \ldots, N$).*

*Let $K \cap H = \phi$. We say H is N-conjugate with K if there exists $x_i \in H_i$ such that $x_i K_i = H_i$ (or $K_i x_i = H_i$) for $i = 1, 2, \ldots, N$ 'or' in the mutually exclusive sense.*

Now we define the notion of inner commutative N-groupoid.

**DEFINITION 4.2.9:** *Let G be a N-groupoid. We say G is inner commutative if every sub N-groupoid of G is commutative.*

We now define the notion of weakly inner commutative N-groupoid.

**DEFINITION 4.2.10:** *Let $G = (G_1 \cup G_2 \cup \ldots \cup G_N, *_1, *_2, \ldots, *_N)$ be a N-groupoid if G has at least one non trivial N-sub groupoid which is commutative, then we call G a weakly inner commutative N-groupoid.*

We have got the following relations, which can be easily checked.

1. Every commutative N-groupoid is inner commutative.
2. Every inner commutative N-groupoid is weakly inner commutative.
3. A weakly inner commutative N-groupoid is not in general inner commutative.
4. All inner commutative N-groupoids need not in general be commutative.

Now we proceed on to define the notion of zero divisors.

**DEFINITION 4.2.11:** *Let $G = (G_1 \cup G_2 \cup \ldots \cup G_N, *_1, *_2, \ldots, *_N)$ be a N-groupoid. An element x in G is said to be a zero divisor if their exists a y in G such that $x *_i y = y *_i x = 0$ for some i in {1, 2, …, N}.*

We define N-centre of a N-groupoid G.

**DEFINITION 4.2.12:** *Let $G = (G_1 \cup G_2 \cup \ldots \cup G_N, *_1, *_2, \ldots, *_N)$ be a N-groupoid. The N-centre of $(G_1 \cup G_2 \cup \ldots \cup G_N, *_1, \ldots, *_N)$ denoted by NC (G) = $\{x_1 \in G_1 \mid x_1 a_1 = a_1 x_1$ for all $a_1 \in G_1\} \cup \{x_2 \in G_2 \mid x_2 a_2 = a_2 x_2$ for all $a_2 \in G_2\} \cup \ldots \cup \{x_N \in G_N \mid x_N a_N = a_N x_N$ for all $a_N \in G_N\}$ = $\bigcup_{i=1}^{N} \{x_i \in G \mid x_i a_i = a_i x_i \text{ for all } a_i \in G_i\}$ = NC (G).*

**DEFINITION 4.2.13**: *Let $G = (G_1 \cup G_2 \cup \ldots \cup G_N, *_1, *_2, \ldots, *_N)$ be a N-groupoid. We say a, b $\in$ G is a N-conjugate pair if a = b $*_i$ $x_i$ (if a, b $\in G_i$ for some $x_i \in G_i$ or a = $x_i *_i b$) and b = a $*_i$ $y_i$ (or $y_i *_i$ a for some $y_i \in G_i$) (if a, b $\in G_j$, $i \neq j$ we choose $x_j$ and $y_j$ from $G_j$ ). Similarly we define right N-conjugate only for a, b $\in G_i$ (or $G_j$ ) a $*_i x_i$ = b and b $*_i$ $y_i$ = a for $x_i \in G_i$.*

We can define left N-conjugate in an analogous way.

**DEFINITION 4.2.14:** *Let $G = (G_1 \cup \ldots \cup G_N, *_1, \ldots, *_N)$ be a non empty set with N-binary operations. If atleast one of $(G_i, *_i)$ is a loop and G a N-groupoid then we call G a N-quasi loop. The number of distinct elements in G is called the order of G. G can be a finite or an infinite N-quasi loop.*

We just give an illustration of a N quasi loop.

***Example 4.2.5:*** Let G = $\{G_1 \cup G_2 \cup G_3 \cup G_4, *_1, *_2, *_3, *_4)$ where $(G_1, *_1)$ = S(3) the semigroup under the operation composition of maps. $(G_2, *_2) = \{Z_{10} \mid a *_2 b = 3a + 2b$ (mod 10)} is a groupoid, $G_3$ = {$Z_9$, $*_3$ is usual multiplication modulo 9}. Clearly $G_3$ is a semigroup. $G_4$ is the loop under the operation $*_4$ given by the following table:

| $*_4$ | e | a | b | c | d |
|---|---|---|---|---|---|
| e | e | a | b | c | d |
| a | a | e | c | d | b |
| b | b | d | a | e | c |
| c | c | b | d | a | e |
| d | d | c | e | b | a |

Thus G is a 4-quasi loop.

We now proceed on to define sub N-quasi loop and normal sub N-quasi loop.

**DEFINITION 4.2.15:** *Let ($G_1 \cup G_2 \cup \ldots \cup G_N$, $*_1, \ldots, *_N$) be a N-quasi loop. Let K = ($K_1 \cup K_2 \cup \ldots \cup K_N$, $*_1, \ldots, *_N$ ) be a proper subset of G. K is said to be a normal sub N-quasi loop of G if K is itself a N-quasi loop such that for all loops $G_i$, $K_i$ is a normal subloop of $G_i$ and for all groupoids and semigroups $G_j$, $K_j$ are ideals of $G_j$.*

It is very important to note that in general in any finite N-quasi loop a proper sub N-quasi loop K of G does not in general divide order of G i.e. o(K) $\nmid$ o (G).

**DEFINITION 4.2.16:** *Let G = ($G_1 \cup G_2 \cup \ldots \cup G_N$, $*_1, \ldots, *_N$) be a N-quasi loop. A proper subset P = ($P_1 \cup P_2 \cup \ldots \cup P_N$, $*_1, \ldots, *_N$) of G is said to be sub N-quasi loop if P itself is a N-quasi loop under the operations of G.*

Now we define a special case of N-quasi loops.

**DEFINITION 4.2.17:** *Let G = ($G_1 \cup G_2 \cup \ldots \cup G_N$, $*_1, \ldots, *_N$) be a finite N-quasi loop. If the order of every proper sub N-quasi loop, P = {$P_1 \cup P_2 \cup \ldots \cup P_N$, $*_1, \ldots, *_N$} divides the order of*

*G, then we call G a Lagrange N-quasi loop. If the order of atleast one of the proper sub N-quasi loops P of G divides the order of G then G is said to be a weakly Lagrange N-quasi loop. If sub N-quasi loop exist and none of its order divides the order of G then G is said to be a non-Lagrange N-quasi loop.*

Thus given any finite N-quasi loop it falls only under one of the three types.

1. Non Lagrange N-quasi loop.
2. Lagrange N-quasi loop and
3. Weakly Lagrange N-quasi loop.

***Example 4.2.6:*** Let $\{G_1 \cup G_2 \cup G_3, *_1, *_2, *_3)$ where $G_1 = S(3)$, $G_2$ is given by the table

| $*_2$ | $x_1$ | $x_2$ | $x_3$ |
|---|---|---|---|
| $x_1$ | $x_1$ | $x_3$ | $x_2$ |
| $x_2$ | $x_2$ | $x_1$ | $x_3$ |
| $x_3$ | $x_3$ | $x_2$ | $x_1$ |

and $G_3$ is given by the table:

| $*_3$ | e | $a_1$ | $a_2$ | $a_3$ | $a_4$ | $a_5$ | $a_6$ | $a_7$ |
|---|---|---|---|---|---|---|---|---|
| e | e | $a_1$ | $a_2$ | $a_3$ | $a_4$ | $a_5$ | $a_6$ | $a_7$ |
| $a_1$ | $a_1$ | e | $a_5$ | $a_2$ | $a_6$ | $a_3$ | $a_7$ | $a_4$ |
| $a_2$ | $a_2$ | $a_5$ | e | $a_6$ | $a_3$ | $a_2$ | $a_4$ | $a_1$ |
| $a_3$ | $a_3$ | $a_2$ | $a_6$ | e | $a_7$ | $a_4$ | $a_1$ | $a_5$ |
| $a_4$ | $a_4$ | $a_6$ | $a_3$ | $a_7$ | e | $a_1$ | $a_5$ | $a_2$ |
| $a_5$ | $a_5$ | $a_3$ | $a_7$ | $a_4$ | $a_1$ | e | $a_2$ | $a_6$ |
| $a_6$ | $a_6$ | $a_7$ | $a_4$ | $a_1$ | $a_5$ | $a_2$ | e | $a_3$ |
| $a_7$ | $a_7$ | $a_4$ | $a_1$ | $a_5$ | $a_2$ | $a_6$ | $a_3$ | e |

$o(G) = 27 + 3 + 8 = 38.$

Now consider the sub N-quasi loop given by $P = P_1 \cup P_2 \cup P_3$ where

$P_1 = \{x_1, x_2, x_3\}$,

$P_2 =$

$$\left\{\begin{pmatrix}1 & 2 & 3\\1 & 2 & 3\end{pmatrix}, \begin{pmatrix}1 & 2 & 3\\1 & 3 & 2\end{pmatrix}, \begin{pmatrix}1 & 2 & 3\\2 & 1 & 3\end{pmatrix}, \begin{pmatrix}1 & 2 & 3\\3 & 2 & 1\end{pmatrix}, \begin{pmatrix}1 & 2 & 3\\2 & 3 & 1\end{pmatrix}, \begin{pmatrix}1 & 2 & 3\\3 & 1 & 2\end{pmatrix}\right\}$$

and $P_3 = \{e, a_1\}$.

Now o(P) = 3 + 6 + 2 = 11. Clearly o(P) $\not|$ o(G).

Now we proceed on to define the notion of p-Sylow sub N-quasi loop.

**DEFINITION 4.2.18:** *Let G = $(G_1 \cup G_2 \cup \ldots \cup G_N, *_1, *_2, \ldots, *_N)$ be a N-quasi loop of finite order. If p is a prime such that $p^\alpha$ / o(G) and $p^{\alpha+1} \not|$ o(G) and we have a sub-N-quasi loop P of order $p^\alpha$, then we say P is a p-Sylow sub N-quasi loop. If for every prime p / o(G) we have a p-Sylow sub N-quasi loop, then we call G a Sylow N-quasi loop.*

**DEFINITION 4.2.19:** *Let G = $(G_1 \cup G_2 \cup \ldots \cup G_N, *_1, \ldots, *_N)$ be a non empty set if each $G_i$ is either a loop or a semigroup, i = 1, 2, …, N then we call G a N-quasi semigroup. i.e. in case of N-quasi loops we can have $G_i$ to be loops and not groupoids and semigroups but in case of N-quasi semigroups we cannot have $G_i$ to be a groupoid it can only be a loop and a semigroup.*

It is interesting to note that All N-quasi semigroups are N-quasi loops. We still define a new notion called N-quasi groupoid.

**DEFINITION 4.2.20:** *Let G = $(G_1 \cup G_2 \cup \ldots \cup G_N, *_1, \ldots, *_N)$ be a non empty set with N-binary operations and each $G_i$ is a proper subset of G; if some of the $G_i$'s are groups and the rest of them are groupoids, then we call G a N-quasi groupoid.*

***Example 4.2.7:*** Let G = $\{G_1 \cup G_2 \cup G_3 \cup G_4, *_1, *_2, *_3, *_4\}$ where

$G_1 = \langle g \mid g^3 = 1\rangle$,

$G_2 = \{Z_{10}$, group under +},

$G_3$ = $\{Z_6 \mid a * b = 2a + 3b \pmod 6), a, b \in Z_6\}$ and $G_4$ is given by the following table:

| $*_4$ | $x_0$ | $x_1$ | $x_2$ | $x_3$ |
|---|---|---|---|---|
| $x_0$ | $x_0$ | $x_1$ | $x_2$ | $x_3$ |
| $x_1$ | $x_3$ | $x_0$ | $x_1$ | $x_2$ |
| $x_2$ | $x_2$ | $x_3$ | $x_0$ | $x_1$ |
| $x_3$ | $x_1$ | $x_2$ | $x_3$ | $x_0$ |

Clearly G is a 4 quasi groupoid.

If the number of distinct elements of a N-quasi groupoid G is *finite* we say G is a *finite N-quasi groupoid* if G has *infinite number* of elements then we say G is an *infinite N-quasi groupoid*.

We can define the notion of Lagrange N-quasi groupoid, non-Lagrange N-quasi groupoid and p-Sylow sub N-quasi groupoid. We define Cauchy elements of a N-quasi loop.

**DEFINITION 4.2.21:** *Let $G = \{G_1 \cup G_2 \cup \ldots \cup G_N, *_1, \ldots, *_N\}$ be a N-quasi loop of finite order. An element x of G such that $x^t = e$ (e the identity element of G), $t > 1$ with $t / o(G)$ is called the Cauchy element of G.*

We see in case of groups every element is a Cauchy element but in case of N-quasi loops only a few of the elements may be Cauchy elements, if the order of the N-quasi loop is a prime, p none of the elements in it would be a Cauchy element even though we will have $x \in G$ with $x^t = e$, $t > 1$, $t \neq p$ as $t \not| \; p$.

It is left as exercise for the reader to give examples and also derive some interesting results using these new concepts.

## 4.3 Smarandache N-groupoid

Now we proceed on to define the notion of Smarandache N-groupoids and give some of its properties.

**DEFINITION 4.3.1:** *Let $G = \{G_1 \cup G_2 \cup \ldots \cup G_N, *_1, \ldots, *_N\}$ be a non empty set with N-binary operations and $G_i$ are proper subsets of G. We call G a Smarandache N-groupoid (S-N-groupoid) if some of the $G_j$ are S-groupoids and the rest of them are S-semigroups; i.e. each $G_i$ is either a S-semigroup or a S-groupoid.*

We say $G = \{G_1 \cup G_2 \cup \ldots \cup G_N, *_1, *_2, \ldots, *_N\}$ is a finite S-N-groupoid if number of distinct elements in G is finite; otherwise we say G is an infinite S-N-groupoid.

**DEFINITION 4.3.2:** *Let $G = \{G_1 \cup G_2 \cup \ldots \cup G_N, *_1, \ldots, *_N\}$ be a N-groupoid. A non empty proper subset $K = (K_1 \cup K_2 \cup \ldots \cup K_N, *_1, \ldots, *_N)$ is said to be a Smarandache sub-N-groupoid (S-sub N-groupoid) if K itself is a S-N-groupoid.*

The reason for taking a N-groupoid and not a S-N-groupoid is evident from the following result.

***Result:*** Let $G = \{G_1 \cup G_2 \cup \ldots \cup G_N, *_1, \ldots, *_N\}$ be a N-groupoid. If G has a proper subset $K = \{K_1 \cup K_2 \cup \ldots \cup K_N, *_1, \ldots, *_N\}$ which is a S-sub N-groupoid then G is a S-N-groupoid.

Now we proceed on to define the notion of Smarandache-commutative N-groupoid.

**DEFINITION 4.3.3:** *Let $G = \{G_1 \cup G_2 \cup \ldots \cup G_N, *_1, \ldots, *_N\}$ be a N-groupoid. We say G is a Smarandache commutative N-groupoid (S-commutative N-groupoid) if every S-sub N-groupoid of G is commutative. If at least one S-sub-N-groupoid which is commutative then we call G a Smarandache weakly commutative N-groupoid (S-weakly commutative N-groupoid).*

We have the following result:

**Result:** Every S-commutative N-groupoid in a S-weakly commutative N-groupoid.

**DEFINITION 4.3.4:** *Let $\{G_1 \cup G_2 \cup \ldots \cup G_N, *_1, *_2 \ldots, *_N\}$ be a N-groupoid of finite order. We call a S-sub N-groupoid K of G to be a Smarandache Lagrange sub-N-groupoid (S-Lagrange sub-N-groupoid) of G if o(K) / o(G).*

*If every S-sub N-groupoid K of G is a S-Lagrange sub N-groupoid we call G to be a Smarandache Lagrange N-groupoid (S-Lagrange N-groupoid). If no S-sub N-groupoid K of G is a S-Lagrange sub N-groupoid then we call G a Smarandache non-Lagrange N-groupoid (S-non Lagrange N-groupoid).*

Now we proceed on to define the notion of Smarandache p-Sylow sub N-groupoid.

**DEFINITION 4.3.5:** *Let $G = \{G_1 \cup G_2 \cup \ldots \cup G_N, *_1, \ldots, *_N\}$ be a N-groupoid of finite order. If p is a prime such that $p^{\alpha} / o(G)$ and $p^{\alpha+1} \not/ o(G)$ further if G has a S-sub N-groupoid P of order $p^{\alpha}$ then we call P a Smarandache p-Sylow sub N-groupoid (S-p-Sylow sub N-groupoid). If for every prime we have a S-p-Sylow sub N-groupoid then we call the N-groupoid G to be a Smarandache Sylow N-groupoid (S-Sylow N-groupoid).*

Now we give some properties about the elements.

**DEFINITION 4.3.6:** *Let $G = \{G_1 \cup G_2 \cup \ldots \cup G_N, *_1, \ldots, *_N\}$ be a S-N-groupoid of finite order. An element $x \in G$ is said to be a Smarandache Cauchy element (S-Cauchy element) if $x^n = e$ and n / o(G). If the S-N-groupoid has atleast a Cauchy element then we call G a Smarandache Cauchy N-groupoid (S-Cauchy N-groupoid). An element $x \in P$ is a S-special Cauchy element if $x^t = e$ and t / o(P), P a S-sub N groupoid of G.*

**DEFINITION 4.3.7:** *Let $G = \{G_1 \cup G_2 \cup \ldots \cup G_N, *_1, *_2, \ldots, *_N\}$ be a N-groupoid. We call G an idempotent N-groupoid if every element of G is an idempotent. We call the S-N-groupoid to be a Smarandache idempotent N-groupoid (S-idempotent N-groupoid) if every element in every S- sub N-groupoid in G is a S-idempotent.*

**DEFINITION 4.3.8:** *Let $G = \{G_1 \cup G_2 \cup \ldots \cup G_N, *_1, *_2 \ldots, *_N\}$ be a N-groupoid. A non empty proper subset $P = (P_1 \cup P_2 \cup \ldots \cup P_N, *_1, *_2, \ldots, *_N)$ of G is said to be a Smarandache left N-ideal (S-left N ideal) of the N-groupoid G if*

i. *P is a S-sub N-groupoid.*
ii. *Each $(P_i, *_i)$ is a left ideal of $G_i$, $1 \leq i \leq N$.*

*On similar lines we can define S-right N-ideal. If P is both a S-right ideal and S-left ideal then we say P is a Smarandache N-ideal (S-N-ideal) of G.*

**DEFINITION 4.3.9:** *Let $G = \{G_1 \cup G_2 \cup \ldots \cup G_N , *_1, *_2, \ldots, *_N\}$ be a N-groupoid. A S-sub N-groupoid $K = (K_1 \cup K_2 \cup \ldots \cup K_N, *_1, \ldots, *_N)$ is said to be a Smarandache normal sub N-groupoid (S-normal sub N-groupoid) of G if*

i. $a_i P_i = P_i a_i$, $a_i \in P_i$ .
ii. $P_i (x_i y_i) = (P_i x_i) y_i$ . $x_i , y_i \in P_i$ .
iii. $x_i (y_i P_i) = (x_i y_i)P_i$ .

*for all $a_i$, $x_i y_i \in G$.*

*We call the Smarandache simple N-groupoid (S-simple N-groupoid) if it has no non trivial S- normal sub- N-groupoids.*

**DEFINITION 4.3.10:** *Let $G = \{G_1 \cup G_2 \cup \ldots \cup G_N , *_1, *_2, \ldots, *_N\}$ be a finite S- N-groupoid.*

*We say G to be a Smarandache normal N-groupoid (S-normal N-groupoid) if the largest S-sub N-groupoid of G is a normal N-groupoid. We call a S-sub N-groupoid to be the largest if the number of elements in it is the maximum.*

For more please refer [35-40]

**DEFINITION 4.3.11:** *Let $G = \{G_1 \cup G_2 \cup \ldots \cup G_N, *_1, *_2, \ldots, *_N\}$ be a N-groupoid. Let $H = \{H_1 \cup H_2 \cup \ldots \cup H_N, *_1, *_2, \ldots,$*

*$*_N$} and K = {$K_1 \cup K_2 \cup \ldots \cup K_N$, $*_1$, $*_2$, …, $*_N$} be any two S-sub N-groupoids of G. We say H is Smarandache N-conjugate (S-N-conjugate) with K if there exists $x_i \in H_i$ and $x_j \in H_j$ with $x_i K_i = H_i$ (or $K_i x_i$) and $x_j K_j = H_j$ (or $K_j x_j$) 'or' in the mutually exclusive sense.*

**DEFINITION 4.3.12:** *Let G be a N-groupoid. We say G is Smarandache inner commutative N-groupoid (S-inner commutative N-groupoid) if every S-sub N-groupoid of G is commutative.*

Now we proceed on to define Smarandache Moufang N-groupoid.

**DEFINITION 4.3.13:** *Let G = $G_1 \cup G_2 \cup \ldots \cup G_N$, $*_1$, $*_2$, …, $*_N$} be a N-groupoid. We say G is a Smarandache Moufang N-groupoid (S-Moufang N-groupoid) if G has a proper subset P where P is a S-sub N-groupoid of G and all elements of P satisfy the Moufang identity i.e. P is a Moufang N-groupoid; likewise we call a N-groupoid G = {$G_1 \cup G_2 \cup \ldots \cup G_N$, $*_1$, $*_2$, …, $*_N$} to be a Smarandache Bol N-groupoid (S-Bol N-groupoid) if G has a S-sub N-groupoid P where P satisfies the Bol identity, similarly if the N-groupoid G = { $G_1 \cup G_2 \cup \ldots \cup G_N$, $*_1$, $*_2$, …, $*_N$} has a S-sub N-groupoid K which is a Bruck N groupoid then we call G a Smarandache Bruck N-groupoid (S-Bruck N-groupoid).*

The notions of Smarandache right (left) alternative (S-right (left) alternative) or alternative N-groupoid is defined in a similar way. Thus for a N-groupoid to satisfy the Smarandache identity we demand a proper S-sub-N-groupoid of G must satisfy the identity.

**DEFINITION 4.3.14:** *Let G = { $G_1 \cup G_2 \cup \ldots \cup G_N$, $*_1$, $*_2$, …, $*_N$} be a N-groupoid. P be a S-sub-N-groupoid of G. We say G is a Smarandache P-N-groupoid (S-P-N-groupoid) if $(x * y) * x = x * (y * x)$ for all $x, y \in P$. If every S-sub N-groupoid of G* *satisfies the identity (x * y) * x = x * (y * x), we call G a Smarandache strong P-N-groupoid (S- strong P-N- groupoid).*

On similar lines we can define Smarandache strong Moufang N-groupoid, Smarandache strong Bol N-groupoid, Smarandache strong alternative N-groupoid and so on.

**DEFINITION 4.3.15:** *Let $G = \{G_1 \cup G_2 \cup \ldots \cup G_N, *_1, *_2, \ldots, *_N\}$ be a N quasi loop. We say G is a Smarandache N-quasi loop (S-N-quasi loop) if it has a proper subset $K = (K_1 \cup K_2 \cup \ldots \cup K_N, *_1, *_2, \ldots, *_N) \subset G$ such that K is a N-quasi loop with some of the $K_i$'s S-loops.*

Next we define Smarandache N-quasi semigroup.

**DEFINITION 4.3.16:** *Let $G = \{G_1 \cup G_2 \cup \ldots \cup G_N, *_1, *_2 \ldots, *_N\}$be a N-quasi semigroup. We say G is Smarandache N-quasi semigroup (S-N-quasi semigroup) if G has a proper subset $P = \{P_1 \cup P_2 \cup \ldots \cup P_N, *_1, *_2, \ldots, *_N)$ such that P is a N-quasi group under the operations of G.*

**DEFINITION 4.3.17:** *Let $G = \{G_1 \cup G_2 \cup \ldots \cup G_N, *_1, \ldots, *_N\}$ be a N-quasi groupoid. If $K = \{K_1 \cup K_2 \cup \ldots \cup K_N, *_1, \ldots, *_N\}$ is a proper subset of G and if $(K_i, *_i)$ is a group and $(K_j, *_j)$ is a S-groupoid $i \neq j$, $1 \leq i, j \leq N$ then we call G a Smarandache N-quasi groupoid (S-N-quasi groupoid).*

**DEFINITION 4.3.18:** *Let $G = \{G_1 \cup G_2 \cup \ldots \cup G_N, *_1, *_2, \ldots, *_N\}$ be a N groupoid we say G is a free N-groupoid; if $G = G_1 \cup G_2 \cup \ldots \cup G_N$ and each $G_i$ is freely generated by a subset $A_i$ of $G_i$.*

## 4.4 Application of N-groupoids and S-N-groupoids

Now we give applications of N-groupoids and S-N-groupoids. Now we define the notion of free S-N-groupoids and give their applications.

**DEFINITION 4.4.1:** *Let G = {$G_1 \cup G_2 \cup \ldots \cup G_N$, $*_1, \ldots, *_N$} be a N-groupoid. We say G is a Smarandache free N-groupoid (S-free N-groupoid) if each $G_i$ has a proper subset, which is a freely generated semigroup, $1 \le i \le N$.*

**DEFINITION 4.4.2:** *Let $A_s = (Z, \bar{A}_s\ \delta_s)$ where Z is a finite set of states, $\bar{A}_s = \bar{A}_1 \cup \bar{A}_2 \cup \ldots \cup \bar{A}_N$ and $\bar{\delta}_s = \bar{\delta}_1 \cup \delta_2 \cup \ldots \cup \bar{\delta}_N$ is said to be Smarandache N-semiautomaton (S-N-semiautomaton) if each $\bar{A}_i$ is a free groupoid (i.e. $\bar{A}$ is a free N-groupoid) and $\bar{\delta}_i$: $Z \times \bar{A}_i \to Z$ is a next state function. Thus the S-N-semiautomation contains A = (Z, $\bar{A}, \bar{\delta}$) as a new N-semi automaton which is a proper sub structure of $A_s$.*

*Or equivalently one can define S-N-semiautomaton as one which has a new N-semi automaton as a substructure.*

We define S-sub N-semiautomaton.

**DEFINITION 4.4.3:** $\bar{A}_s = (Z_1, \bar{A}_s, \bar{\delta}'_s)$ *is called the Smarandache N-sub semiautomaton of* $A_s = (Z_2, \bar{A}_s, \bar{\delta}_s)$ *denoted by* $\bar{A}_s \subset A_s$ *if $Z_1 \subset Z_2$ and $\bar{\delta}'_s$ is the restriction of $\bar{\delta}_s$ and is defined on $Z_1 \times \bar{A}_s$ and $\bar{A}_s$ has proper subset $\bar{H} \subset A'_s$ such that $\bar{H}$ is a new N-semi automaton.*

Now like wise we can define the notion of Smarandache N-automaton.

**DEFINITION 4.4.4:** $\bar{K}_s = (Z, \bar{A}_s, \bar{B}_s, \bar{\delta}_s, \bar{\lambda}_s)$ *is defined to be a Smarandache N-automaton (S-N-automaton) if* $\bar{K} = (Z, \bar{A}_s, \bar{B}_s, \bar{\delta}_s, \bar{\lambda}_s)$ *is the new N-automaton and $\bar{A}_s$ and $\bar{B}_s$ are the S-free N-groupoids so that* $\bar{K}_s = (Z, \bar{A}_s, \bar{B}_s, \bar{\delta}_s, \bar{\lambda}_s)$ *is the new automaton got from $\bar{K}_s$ and $\bar{K}$ is strictly contained in $\bar{K}_s$.*

A Smarandache automaton equipped with free N-groupoids can perform multi task unlike a S-automaton.

Now we proceed on to define the notion of Smarandache sub N-automaton.

**DEFINITION 4.4.5:** *Let* $\overline{K}_s = (Z_1, \overline{A}_s, \overline{B}_s, \overline{\delta}_s, \overline{\lambda}_s)$ *is called Smarandache sub N-automaton (S-sub N-automaton of* $\overline{K}'_s = (Z_2, \overline{A}_s, \overline{B}_s, \overline{\delta}_s, \overline{\lambda}_s)$ *denoted by* $\overline{K}'_s \le \overline{K}_s$ *if* $Z_1 \subseteq Z_2$ *and* $\overline{\delta}_s$ *and* $\overline{\lambda}_s$ *are restrictions of* $\overline{\delta}_s$ *and* $\overline{\lambda}_s$ *respectively on* $Z_1 \times \overline{A}_s$ *and has a proper subset* $\overline{H} \subset \overline{K}'_s$ *such that* $\overline{H}$ *is a new N-automaton.*

Also the notion of direct product of S-N-automatons can be defined using the notion of direct product of N-free groupoids. Now we proceed on to define the notion of N-semiautomaton and Smarandache N-semiautomaton.

**DEFINITION 4.4.6:** *Let* $S_N^S = (Z, \overline{A}_N, \overline{\delta}_N)$ *be a triple where* $\overline{A}_N$ *is a free N-groupoid and* $\overline{\delta}_N$ *is a function from Z ×* $\overline{A}_N$ *to Z i.e.* $\overline{\delta}_N = \overline{\delta}_N^1 \cup \overline{\delta}_N^2 \cup ... \cup \overline{\delta}_N^N$ *such that* $\overline{\delta}_N^i : Z \times \overline{A}_N^i \to Z$ *i = 1, 2, …, N.*

$$Z \times \overline{A}_N = (Z \times \overline{A}_N^1) \cup (Z \times \overline{A}_N^2) \cup ... \cup (Z \times \overline{A}_N^N).$$

**DEFINITION 4.4.7:** *Let* $S_N^s = (Z, \overline{A}_N^s, \overline{\delta}_N^s)$ *where* $\overline{A}_N^s$ *is a S-N-groupoid and* $\overline{\delta}_N = \overline{\delta}_N^1 \cup \overline{\delta}_N^2 \cup ... \cup \overline{\delta}_N^N$ *where* $\overline{\delta}_N^{is} : Z \times \overline{A}^{is} \to Z$ *is a function for i = 1, 2, …, N. Here* $P_N^S = (Z, \overline{A}_N^s, \overline{\delta}'_N)$ *is a new N-semi automation properly contained in* $S_N^s$*.*

The notion of sub N-semi automaton, homomorphism of N-semi automaton and direct product of N-semiautomaton can be defined.

**Chapter Five**

# MIXED N-ALGEBRAIC STRUCTURES

In this chapter for the first time we introduce the notion of mixed N-Algebraic structures. We will broadly classify this mixed N-Algebraic structures as

1. Associative N-algebraic structures.
2. Non-associative N-algebraic structures.
3. Associative and non-associative N-algebraic structures.

We also define their Smarandache analogue. For the reader should not be of the opinion we can have N-groupoids, N-groups, N-groupoids, N-semigroups or N-loops. One can have any N-algebraic structure got from the four structures.

This chapter has 3 sections. In section one we introduce the N-group semigroup algebraic structures and analyze them. In section two N-loop groupoids and their properties are studied. A new mixed algebraic structure called N-group loop semigroup groupoid (N-glsg) are introduced in section 3.

## 5.1 N-group semigroup algebraic structure

In this section we proceed onto define the notion of N-group-semigroup algebraic structure and enumerate some of its properties.

**DEFINITION 5.1.1**: *Let* $G = \{G_1 \cup G_2 \cup \ldots \cup G_N, *_1, \ldots, *_N\}$ *where some of the* $G_i$*'s are groups and the rest of them are* *semigroups. $G_i$'s are such that $G_i \not\subseteq G_j$ or $G_j \not\subseteq G_I$ if $i \neq j$, i.e. $G_i$'s are proper subsets of G. $*_1, \ldots, *_N$ are N binary operations which are obviously are associative then we call G a N-group semigroup.*

We can also redefine this concept as follows:

**DEFINITION 5.1.2:** *Let G be a non empty set with N-binary operations $*_1, \ldots, *_N$. We call G a N-group semigroup if G satisfies the following conditions:*

- i. *$G = G_1 \cup G_2 \cup \ldots \cup G_N$ such that each $G_i$ is a proper subset of G (By proper subset $G_i$ of G we mean $G_i \not\subseteq G_j$ or $G_j \not\subseteq G_i$ if $i \neq j$. This does not necessarily imply $G_i \cap G_j = \phi$).*
- ii. *$(G_i, *_i)$ is either a group or a semigroup, $i = 1, 2, \ldots, N$.*
- iii. *At least one of the $(G_i, *_i)$ is a group.*
- iv. *At least one of the $(G_j, *_j)$ is semigroup $i \neq j$.*

*Then we call $G = \{G_1 \cup G_2 \cup \ldots \cup G_N, *_1, \ldots, *_N\}$ to be a N-group semigroup ( $1 \leq i, j \leq N$).*

Now we first illustrate this algebraic structure with some examples.

***Example 5.1.1:*** Let $G = \{G_1 \cup G_2 \cup \ldots \cup G_5, *_1, *_2, \ldots, *_5\}$ where

$(G_1, *_1)$ = $\{Z_{10}$, under multiplication modulo 10} is a semigroup,

$\{G_2, *_2\}$ = $A_4$ is the alternating subgroup of $S_4$ is a group,

$G_3$ = S(3), the semigroup of all mapping from the set {1 2 3} to itself,

$$G_4 = \left\{ \begin{pmatrix} a_{11} & a_{12} & a_{13} \\ a_{21} & a_{22} & a_{23} \\ a_{31} & a_{32} & a_{33} \end{pmatrix} \middle| a_{ij} \in Q \right\}$$

semigroup under matrix multiplication

and

$G_5$ = $G' = \{g \mid g^{12} = 1\}$ the cyclic group of order 12.

Clearly G is a 5-group semigroup.

Now we as in case of any algebraic structure define the notion of order of that algebraic structure. Let $G = \{G_1 \cup G_2 \cup \ldots \cup G_N, *_1, \ldots, *_N\}$ be a N-group semigroup. The number of distinct elements in G is called the order of the N-group semigroup. If the order is finite then we call G a finite N-group semigroup. If the order is not finite then we call G an infinite N-group-semigroup.

**DEFINITION 5.1.3:** *Let $G = \{G_1 \cup G_2 \cup \ldots \cup G_N, *_1, \ldots, *_N\}$ be a N-group-semigroup. We say G is a commutative N-group semigroup if each $(G_i, *_i)$ is a commutative structure, i = 1, 2, …, N.*

We give an example of the commutative N-group-semigroup.

***Example 5.1.2:*** Let $G = \{G_1 \cup G_2 \cup G_3, *_1, *_2, *_3\}$ be a 3-group semigroup. Let

$\{G_1, *_1\}$ = $\langle g \mid g^7 = 1\rangle$ be a cyclic group of order 7,

$\{G_2, *_2\}$ = $\{Z_{12}$, semigroup under multiplication modulo 12}

and

$\{G_3, *_3\}$ = $\{Z_{11} \setminus \{0\}$, is a group under multiplication modulo 11}.

Clearly G is a commutative 3-group semigroup. Now the order of the 3-group semigroup is 27, i.e., o(G) = 7 + 11 + 9 = 27.

The 4-group semigroup given in example 5.1.1 is non commutative and is of infinite order and the 3-group semigroup given in example 5.1.2 is commutative and is of finite order.

Now we define N-subgroup of a N-group semigroup, N-subsemigroup of a N-group semigroup and sub N-group semigroup.

**DEFINITION 5.1.4:** *Let $G = \{G_1 \cup G_2 \cup \ldots \cup G_N, *_1, \ldots, *_N\}$ be a N-group. A proper subset P of G where $(P_1 \cup P_2 \cup \ldots \cup P_N, *_1, \ldots, *_N)$ is said to be a N-subgroup of the N-group semigroup G if each $(P_i, *_i)$ is a subgroup of $(G_i, *_i)$; $i = 1, 2, \ldots, N$.*

***Example 5.1.3:*** Let $G = \{G_1 \cup G_2 \cup G_3 \cup G_4, o_1, o_2, o_3, o_4\}$ be a 4-group semigroup where

$G_1$ = {Z, the set of positive and negative integers under '+'}.
$G_2$ = S(3),
$G_3$ = $\left\{\begin{pmatrix} a_{11} & a_{12} \\ a_{21} & a_{22} \end{pmatrix} = A \,\middle|\, a_{ij} \in Q\right\}$
the semigroup under multiplication and
$G_4$ = $\{g \mid g^{12} = 1\}$.

Clearly G is a 4-group semigroup.

The set $P = P_1 \cup P_2 \cup P_3 \cup P_4$ is a N-subgroup of the N-group semigroup G; where

$P_1$ = $\{2n / n \in Z\}$,
$P_2$ = $S_3$ the symmetric group of degree 3,
$P_3$ = $\left\{\begin{pmatrix} a_{11} & a_{12} \\ a_{21} & a_{22} \end{pmatrix} = A \,\middle|\, a_{11}a_{12} - a_{12}a_{21} \neq 0\right\} \subset G_3$
is a group of $2 \times 2$ matrices under multiplication and
$P_4$ = $\{g^3, g^6, g^9, 1\}$.

P is a 4-subgroup of the 4-group semigroup G.

Now we proceed on to illustrate by an example of the N-subsemigroup.

***Example 5.1.4:*** Let $G = \{G_1 \cup G_2 \cup G_3 \cup G_4 \cup G_5, *_1, *_2, *_3, *_4, *_5\}$ be a 5-group semigroup where

$G_1$ = {Z; the set of positive and negative integers with 0, a group under '+'},

$G_2$ = $A_{n\times n}$ = {($a_{ij}$), semigroup under 'matrix multiplication' $a_{ij} \in Q$},

$G_3$ = {$Z_{16}$; semigroup under multiplication modulo 16},

$G_4$ = {3 × 2 matrices with entries from Q}, $G_4$ is a group under matrix addition and

$G_5$ = S(4), semigroup of mapping of the set {1 2 3 4} to itself, under composition of mappings.

It is easy to see that G is a 5-group semigroup.
Now consider the set K = $K_1 \cup K_2 \cup K_3 \cup K_4 \cup K_5$ where

$K_1$ = {$Z^+$, set of all positive integers under '+'}, $K_1$ is a semigroup under '+'.

$K_2$ = {$A_{n \times n} = (a_{ij}) \mid a_{ij} \in Z \subseteq Q$}, $K_2$ is semigroup under matrix multiplication.

$K_3$ = {0, 2, 4, 6, 8, 10, 12, 14} is a semigroup under multiplication modulo 16,

$K_4$ = {set of all 3 × 2 matrices with entries from $Q^+$, the set of positive rationals}, a semigroup under matrix addition and

$K_5$ = {Set of all maps of the set (1 2 3 4) fixing 4 in the same place}.

Clearly $K_5$ is a semigroup isomorphic to S(3), the semigroup of all maps of the set (1 2 3) to itself.

Hence K = {$K_1 \cup K_2 \cup K_3 \cup K_4 \cup K_5$, $*_1$, $*_2$, $*_3$, $*_4$, $*_5$} is a 5- subsemigroup of the 5-group semigroup G.

Now we proceed on to give an example of a N-subgroup semigroup of a N-group semigroup, then we define both the notions of N- subsemigroup and N-group semigroup G.

***Example 5.1.5:*** Let G = {$G_1 \cup G_2 \cup G_3 \cup G_4$, $*_1$, … ,$*_4$}be a 4-group semigroup, where

$G_1$ = $A_4$, the alternating group of degree four.

$G_2$ = S(5) the semigroup of (1 2 3 4 5) to itself.

$G_3$ = $D_{2.7}$ = {a, b | $a^2 = b^7 = 1$, b a b = a}. The dihedral group of order 14.

$G_4$ = {set of all 3 × 7 matrices with entries from Q under, '+'} is a semigroup.

Now we denote by $P = P_1 \cup P_2 \cup P_3 \cup P_4$ where

$$P_1 = \left\{ \begin{pmatrix} 1 & 2 & 3 & 4 \\ 1 & 2 & 3 & 4 \end{pmatrix}, \begin{pmatrix} 1 & 2 & 3 & 4 \\ 2 & 1 & 4 & 3 \end{pmatrix} \right\},$$

$P_2$ = S(3) the semigroup got by fixing the elements 4 and 5 in the set (1 2 3 4 5) and mapping only (1 2 3) to itself,

$P_3$ = $\{1, b, \ldots, b^6\}$ and

$P_4$ = {set of all 3 × 7 matrices with entries from Z} $\subset G_4$.

P is a 4-subgroup semigroup of the 4-group semigroup G.

**DEFINITION 5.1.5:** *Let $G = \{G_1 \cup G_2 \cup \ldots \cup G_N, *_1, \ldots, *_N\}$ be a N-group semigroup where some of the $(G_i, *_i)$ are groups and rest are $(G_j, *_j)$ are semigroups, $1 \le i, j \le N$. A proper subset P of G is said to be a N-subsemigroup if $P = \{P_1 \cup P_2 \cup \ldots \cup P_N, *_1, \ldots, *_N\}$ where each $(P_i, *_i)$ is only a semigroup under the operation $*_i$.*

Now we proceed on to define the notion of N-subgroup semigroup of a N-group semigroup.

**DEFINITION 5.1.6:** *Let $G = \{G_1 \cup G_2 \cup \ldots \cup G_N, *_1, \ldots, *_N\}$ be a N-group semigroup. Let P be a proper subset of G. We say P is a N-subgroup semigroup of G if $P = \{P_1 \cup P_2 \cup \ldots \cup P_N, *_1, \ldots, *_N\}$ and each $(P_i, *_i)$ is a group or a semigroup.*

Now we are not always guaranteed of a N-subgroup or N-subsemigroup in a N-group-semigroup. At times even it may be impossible to find N-subgroup semigroup of a N-group semigroup.

So the study of the existence theorem for substructures can be taken up any person who finds interest in it.

Now we proceed on to define the notion of normal N subgroup semigroup of a N-group semigroup.

**DEFINITION 5.1.7**: *Let $G = \{G_1 \cup G_2 \cup \ldots \cup G_N, *_1, \ldots, *_N\}$ be a N-group semigroup. We call a proper subset P of G where $P = \{P_1 \cup P_2 \cup \ldots \cup P_N, *_1, \ldots, *_N\}$ to be a normal N-subgroup semigroup of G if $(G_i, *_i)$ is a group then $(P_i, *_i)$ is a normal subgroup of $G_i$ and if $(G_j, *_j)$ is a semigroup then $(P_j, *_j)$ is an ideal of the semigroup $G_j$. If G has no normal N-subgroup semigroup then we say N-group semigroup is simple.*

Now we provide some examples for the same.

***Example 5.1.6:*** Let $G = \{G_1 \cup G_2 \cup G_3, *_1, *_2, *_3\}$ be a 3-group semigroup where

$G_1$ = $\{g \mid g^{19} = 1\}$ be the cyclic group of order 19,
$G_2$ = $\{A_5$, the alternating subgroup of the symmetric group $S_5\}$ and
$G_3$ = S(3) be the symmetric semigroup.

Clearly G is simple 3-group semigroup.

Now we give an example of a N-group semigroup, which has non trivial normal N-subgroup semigroup i.e. N-group semigroup which is not simple.

***Example 5.1.7:*** Let $G = \{G_1 \cup G_2 \cup G_3 \cup G_4, *_1, \ldots, *_4\}$ be a 4-group semigroup where

$G_1$ = $S_3$, symmetric group of degree 3,
$G_2$ = {Z; the set of integers under multiplication is a semigroup},
$G_3$ = $\{Z_{12}$, under multiplication modulo 12}, $G_3$ under multiplication is a semigroup and
$G_4$ = $\{Z_{15}$, group under '+' modulo 15}.

Consider the set $P = P_1 \cup P_2 \cup P_3 \cup P_4$, where

$P_1$ = $A_3$; clearly $A_3$ is a normal subgroup of $S_3$,
$P_2$ = {2Z, 2Z is an ideal of the semigroup Z under multiplication},

$P_3$ = {0, 6} is an ideal of the semigroup $G_3$ and
$P_4$ = {0, 5, 10} is a normal subgroup of $G_4$.

Thus P is a normal 4-subgroup semigroup of the 4-group semigroup G.

Having defined the order of the N-group semigroup G the next natural question would be suppose the N-group semigroup G is of finite order and if P is a non trivial N-subgroup semigroup (or N-subgroup or N-subsemigroup) will order of P divide the order of G. To this end we give an example of a finite N-group semigroup.

***Example 5.1.8:*** Let G = {$G_1 \cup G_2 \cup G_3 \cup G_4$, $*_1$, $*_2$, $*_3$, $*_4$} be a 4-group semigroup of finite order, where

$G_1$ = $A_4$, the alternating subgroup of $S_4$,
$G_2$ = {$Z_{10}$, the semigroup under multiplication modulo 10},
$G_3$ = {S(3), the semigroup of mapping from the set (1 2 3) to itself} and
$G_4$ = {g | $g^5$ = e}.

Clearly order of G is 12 + 10 + 27 + 5 = 54. Now we will find N-subgroup semigroup of G.
Take P = $P_1 \cup P_2 \cup P_3 \cup P_4$ where

$$P_1 = \left\{ \begin{pmatrix} 1 & 2 & 3 & 4 \\ 1 & 2 & 3 & 4 \end{pmatrix} \begin{pmatrix} 1 & 2 & 3 & 4 \\ 3 & 4 & 1 & 2 \end{pmatrix} \right\},$$

$P_2$ = {0, 5} subsemigroup under multiplication modulo 10,
$P_3$ = $S_3$ and
$P_4$ = {e}.
o (P) = 2 + 2 + 6 + 1 = 11.

Clearly $11 \nmid 54$, so we see in general the order of a N-subgroup semigroup of a finite N-group semigroup need not divide the order of the N-group semigroup.

Now consider the subset K of G. Let K = $K_1 \cup K_2 \cup K_3 \cup K_4$ where

$K_1 = \left\{\begin{pmatrix} 1 & 2 & 3 & 4 \\ 1 & 2 & 3 & 4 \end{pmatrix}\right\}$,

$K_2 = Z_{10}$,

$K_3 = S_3$ and

$K_4 = \{e\}$.

Clearly order of K = 1 + 10 + 6 + 1 = 18, 18 /54.
Now consider the subset $T = T_1 \cup T_2 \cup T_3 \cup T_4$ where

$T_1 = \left\{\begin{pmatrix} 1 & 2 & 3 & 4 \\ 1 & 2 & 3 & 4 \end{pmatrix}\right\}$,

$T_2 = \{0, 2, 4, 6, 8\}$ semigroup under multiplication modulo 10.

$T_3 = \left\{\begin{pmatrix} 1 & 2 & 3 \\ 1 & 2 & 3 \end{pmatrix}, \begin{pmatrix} 1 & 2 & 3 \\ 2 & 1 & 3 \end{pmatrix}\right\}$ and

$T_4 = \{e\}$.

The order of T = 1 + 5 + 2 + 1 = 9, 9 / 54. Thus the order of this N-subgroup semigroup divides the order of G. Thus we see in the same finite N-group semigroup we have certain N-subgroup semigroup whose order divides the order of the N-group semigroup and some N-subgroup semigroup whose order does not divide the order of the N-group semigroup. So we define the following:

**DEFINITION 5.1.8:** *Let $G = \{G_1 \cup G_2 \cup \ldots \cup G_N, *_1, \ldots, *_N\}$ be a finite N-group semigroup. If the order of every N-subgroup semigroup divides order the N-group semigroup G then we call G a Lagrange N-group semigroup. If there is atleast one N-subgroup semigroup P such that o(P) / o(G) then we call G to be a weakly Lagrange N-group semigroup. If their does not exist any N-subgroup semigroup P of the N-group semigroup G such that o(P) / o(G) then we call G the free Lagrange N-group semigroup.*

We have all three types of N-group semigroups of finite order.

Now we proceed on to define the notion of Cauchy element of a finite N-group semigroup G.

**DEFINITION 5.1.9:** *Let $G = \{G_1 \cup G_2 \cup \ldots \cup G_N, *_1, \ldots, *_N\}$ be a finite N-group semigroup. Let $x \in G$; we say x is a Cauchy element of G if $x^t = 1$ ($t > 1$ and t is the smallest integer) then t / o(G). We say the N-group semigroup to be a Cauchy N-group semigroup if every element x in G such that $x^t = 1$ is a Cauchy element of G. It is very important to note that in general every element x in G need not be such that $x^t = 1$ for we may have x in G such that $x^2 = x$ and so on or $x^n = 0$. Thus if every element x which is such that $x^t = 1$ ($t > 1$) satisfies t / o(G) then only we call the N-group semigroup to be a Cauchy N-group semigroup.*

Now we illustrate this with examples.

***Example 5.1.9:*** Let $G = \{G_1 \cup G_2 \cup G_3 \cup G_4, *_1, *_2, *_3, *_4\}$ be a 4-group semigroup of finite order where

$G_1$ = $\{g \mid g^7 = e\}$ the cyclic group of order 7,
$G_2$ = S(3), semigroup of mappings of the set (1 2 3) to itself,
$G_3$ = $\{Z_{11}$, group under '+' modulo 11$\}$ and
$G_4$ = $\{Z_8$, semigroup under multiplication modulo 8$\}$.

Clearly o (G) = 7 + 27 + 11 + 8 = 53. Now since order of the 4-group semigroup is a prime we see G has no Cauchy elements. In truth we see $g^7 = e$ is an element of finite order

$$\left\{\begin{pmatrix} 1 & 2 & 3 \\ 2 & 3 & 1 \end{pmatrix}\right\}^3 = \begin{pmatrix} 1 & 2 & 3 \\ 1 & 2 & 3 \end{pmatrix}$$

is an element of order 3 and so on. As o (G) = 53 a prime this 4-group semigroup has no Cauchy element.

Next we go to another example.

***Example 5.1.10:*** Let G = {$G_1 \cup G_2 \cup G_3 \cup G_4 \cup G_5$, $*_1$, $*_2$, $*_3$, $*_4$, $*_5$} be a N-group semigroup of finite order where

$G_1$ = $\langle g \mid g^8 = e\rangle$, cyclic group of order 8.
$G_2$ = {$Z_{12}$, semigroup under multiplication modulo 12}.
$G_3$ = $S_3$, symmetric group of degree 3.
$G_4$ = {$Z_{15}$, semigroup under multiplication 14} and
$G_5$ = $D_{2.5}$ = {a, b | $a^2 = b^5 = 1$, b a b = a} dihedral group of order 10.

Clearly o (G) = 8 + 12 + 6 + 14 + 10 =50.

Take

$$x = \begin{pmatrix} 1 & 2 & 3 \\ 2 & 3 & 1 \end{pmatrix}$$

in G;

$$x^3 = \begin{pmatrix} 1 & 2 & 3 \\ 1 & 2 & 3 \end{pmatrix}$$

but $3 \nmid 50$ so x is not a Cauchy element of G.

Take $g^2 \in G$, $(g^2)^4 = e$ but $4 \nmid 50$ so $g^2$ is not a Cauchy element of G.

Consider

$$y = \begin{pmatrix} 1 & 2 & 3 \\ 1 & 3 & 2 \end{pmatrix} \in G.$$

We have

$$y^2 = \begin{pmatrix} 1 & 2 & 3 \\ 1 & 2 & 3 \end{pmatrix}$$

and 2/ 50 so y is a Cauchy element of G. Consider b ∈ G, we have $b^5 = 1$ thus b is a Cauchy element of G. So we see in this finite N-group semigroup G some elements are Cauchy elements and some of them are not Cauchy elements though the elements considered are of finite order.

Now we see the element 3 ∈ G, $3^6 \equiv 729 \equiv 1 \pmod{14}$, 3 is an element of finite order but $6 \nmid 50$.

Next we consider the following finite N-group semigroup G.

***Example 5.1.11:*** Let G = {$G_1 \cup G_2 \cup G_3 \cup G_4$, $*_1, \ldots, *_4$} be a 4-group semigroup of finite order where

$G_1$ = {$Z_6$, semigroup under multiplication modulo 6},
$G_2$ = $\langle g \mid g^3 = 1\rangle$, cyclic group of order 3,
$G_3$ = {$S_3$, the symmetric group of degree 3} and
$G_4$ = {$Z_3$, under modulo addition}.

Now order of G = 6 + 3 + 6 + 3 = 18.

It is easily verified every element of finite order is a Cauchy element and G is a Cauchy 4-group semigroup. In fact G has element x such that $x^2 = x$ for x = 4 in $G_1$ gives $4^2 = 4$, similarly $y \in G$ such that $y^3 = y$ for take $2 \in G$, $2^3 \equiv 2 \pmod 6$. Thus this 4-group semigroup also has elements which are such that $x^n = 1$.

Now we make the following definition.

**DEFINITION 5.1.10:** *Let G = {$G_1 \cup G_2 \cup \ldots \cup G_N$, $*_1, \ldots, *_N$} be a N-group semigroup of finite order. We say G is a weakly Cauchy N-group semigroup if G has atleast one element x such that $x^t = 1$, (t > 1) and t / o(G).*

We have also seen examples of N-group semigroup which has no Cauchy elements. However we can give a one way condition for the existence of a Cauchy N-group semigroup.

**THEOREM 5.1.1:** *Let G be a N-group semigroup of finite order n where n is a prime. Then G has no Cauchy elements.*

*Proof:* Given G is a N-group semigroup of order n where n is a prime. So even if G has elements say x such that $x^t = 1$ (t > 1) still $t \not| p$. Hence the claim.

So we define a finite N-group semigroup G to be a Cauchy free N-group semigroup if no element of finite order in G divides the order of G.

All finite N-group semigroup of order n, n any prime is a Cauchy free N-group semigroup.

Now we proceed on to define the notion of p-Sylow sub N group semigroup of a N-group semigroup G.

**DEFINITION 5.1.11:** *Let $G = \{G_1 \cup G_2 \cup \ldots \cup G_N, *_1, \ldots, *_N\}$ be a N-group semigroup of finite order. Let p be a prime such that $p^{\alpha} / o(G)$ but $p^{\alpha+1} \not| \; o(G)$ and if G has a N-subgroup semigroup P of order $p^{\alpha}$ then we call P the p-Sylow N-subgroup semigroup of the N-group semigroup G. If for every prime p such that $p^{\alpha} / o(G)$ but $p^{\alpha+1} \not| \; o(G)$ we have a p-Sylow N-subgroup semigroup then we call the N-group semigroup G to be a Sylow N-group semigroup.*

Now we proceed on to illustrate this concept by examples.

***Example 5.1.12:*** Let $G = \{G_1 \cup G_2 \cup G_3, *_1, *_2, *_3\}$ be a 3-group semigroup, where

$G_1 = S_3$,
$G_2 = \{Z_{24}$, is a semigroup under multiplication modulo 24$\}$
and
$G_3 = \{g \mid g^6 = e\}$.

The order of G = 6 + 24 + 6 = 36, 3 / 36, $3^2$ / 36 and $3^3 \not| \; 36$, 2 / 36, $2^2$ / 36 and $2^3 \not| \; 36$.

Take $P = P_1 \cup P_2 \cup P_3$ where

$$P_1 = \left\{\begin{pmatrix}1 & 2 & 3\\ 1 & 2 & 3\end{pmatrix}, \begin{pmatrix}1 & 2 & 3\\ 2 & 3 & 1\end{pmatrix}, \begin{pmatrix}1 & 2 & 3\\ 3 & 2 & 1\end{pmatrix}\right\},$$

$P_2 = \{0, 6, 12, 18\}$ and $P_3 = \{e, g^3\}$.

Clearly o(P) = 3 + 4 + 2 = 9, 9 / 36. Hence P is a 3-Sylow N-subgroup semigroup. Consider $K = K_1 \cup K_2 \cup K_3$ where

$$K_1 = \left\{\begin{pmatrix}1 & 2 & 3\\ 1 & 2 & 3\end{pmatrix}\right\},$$

$K_2 = \{0, 12\}$ and $K_3 = \{e\}$.

Clearly order of K = 1 + 2 + 1 = 4 and K is a 2-Sylow N-subgroup semigroup of G. Thus we can call G the Sylow N-group semigroup. Here N = 3.

Now consider the subset $T = T_1 \cup T_2 \cup T_3$ where

$$T_1 = \left\{ \begin{pmatrix} 1 & 2 & 3 \\ 1 & 2 & 3 \end{pmatrix}, \begin{pmatrix} 1 & 2 & 3 \\ 1 & 3 & 2 \end{pmatrix} \right\},$$

$T_2 = Z_{24}$ and $T_3 = \{e\}$.

The order of T is $2 + 24 + 1 = 27 = 3^3$, $3^3 \nmid 36$. But G has a N-subgroup semigroup of order 27. This property forces us to define some new structure called super Sylow N-group semigroup.

**DEFINITION 5.1.12:** *Let $G = \{G_1 \cup G_2 \cup \ldots \cup G_N, *_1, *_2, \ldots, *_N\}$ be a finite N-group semigroup. Suppose G is a Sylow N-group semigroup and if in addition for every prime p, with $p^{\alpha} / o(G)$ but $p^{\alpha+1} \nmid o(G)$ their exists N-subgroup semigroup of order $p^{\alpha+t}$, $t \geq 1$ then we call the N-group semigroup to be the super Sylow N-group semigroup.*

Now it is important to note that every super Sylow N-group semigroup is obviously a Sylow N-group semigroup, but in general every Sylow N-group semigroup need not be a super Sylow N-group semigroup. Now we proceed on to define still a weaker condition than the Sylow N-group semigroup.

**DEFINITION 5.1.13:** *Let $G = \{G_1 \cup G_2 \cup \ldots \cup G_N, *_1, \ldots, *_N\}$ be a N-group semigroup of finite order. Suppose p is a prime such that $p^{\alpha} / o(G)$ but $p^{\alpha+1} \nmid o(G)$ and G has a N subgroup semigroup $P = \{P_1 \cup P_2 \cup \ldots \cup P_N, *_1, \ldots, *_N\}$ of order $p^t$ ($t < \alpha$) then we call P the weak p-Sylow N-subgroup semigroup of G. If for every prime p such that $p^{\alpha} / o(G)$ but $p^{\alpha+1} \nmid o(G)$ we have a N-subgroup semigroup of order $p^t$ ($t < \alpha$) then we call* *the N-group semigroup to be a weakly Sylow N-group semigroup*.

We now illustrate this by the following example:

***Example 5.1.13:*** Let G = $\{G_1 \cup G_2 \cup G_3, *_1, *_2, *_3\}$ be a 3-group semigroup where

$G_1$ = $\langle g \mid g^{37} = e\rangle$ the cyclic group of order 37,
$G_2$ = $\{Z_5$, the semigroup under multiplication modulo 5$\}$ and
$G_3$ = $\{Z_3$, semigroup under +$\}$.

Now order of G = 37 + 5 + 3 = 45, 3/ 45, $3^2$/ 45, 5/45, $5^2 \not| \ 45$. Clearly G has no p-Sylow N-subgroup semigroups (p = 3).

If we take T = $T_1 \cup T_2 \cup T_3$ where T = {e}, $T_2$ = {0} and $T_3 = Z_3$ we get a 5- N-subgroup semigroup. This is another extreme case where the order of the N-group semigroup is not a prime yet we are not a position to get or make the N-group semigroup into a Sylow N-group semigroup.

***Example 5.1.14:*** Let G = $\{G_1 \cup G_2 \cup G_3, *_1, *_2, *_3\}$ be a 3 group semigroup where

$G_1$ = $\langle g \mid g^{19} = 1\rangle$ the cyclic group of order 19,
$G_2$ = $\{Z_3$, the semigroup under multiplication modulo 3$\}$ and
$G_3$ = {1, –1 group under multiplication},

o(G) = 24, $2^3$/24 and $2^4 \not| \ 24$, 3/24, $3^2 \not| \ 24$.

Take T = $T_1 \cup T_2 \cup T_3$ where $T_1$ = {e}, $T_2$ = {0} and $T_3$ = $G_3$. Order of T is 4 and G has no sub N-group semigroup of order 8. But it has a N-subgroup semigroup of order 4 viz. T.

Clearly G has no N-subgroup semigroup of order 3. So G is a weakly Sylow 3-group semigroup. We can as in case of any N-algebraic structure define homomorphism.

**DEFINITION 5.1.14:** *Let G = $\{G_1 \cup G_2 \cup \ldots \cup G_N, *_1, \ldots, *_N\}$ and K = $\{K_1 \cup \ldots \cup K_N, o_1, \ldots, o_N\}$ be any two N-group*

*semigroups. We can define a map $\phi$ from G to K to be a N-group semigroup homomorphism if and only if $\{G_i, *_i\}$ is a group then $\{K_i, o_i\}$ is also a group. Likewise if $\{G_t, *_t\}$ is a semigroup then $\{K_t, o_t\}$ is a semigroup. We have $\phi = \phi_1 \cup \ldots \cup \phi_N$ where $\phi_i$ $G_i \to K_i$ is either a group homomorphism or a semigroup homomorphism according as $G_i$ and $K_i$ are groups or semigroups.*

If in two N-group semigroups $G = \{G_1 \cup G_2 \cup \ldots \cup G_N, *_1, \ldots, *_N\}$ and $K = \{K_1 \cup K_2 \cup \ldots \cup K_N, o_1, \ldots, o_N\}$ we have $G_i$ happens to be a group and $K_i$ is only a semigroup then we cannot define homomorphism of N-group semigroups.

Now we proceed on to define the notion of conjugate sub N-group semigroups.

**DEFINITION 5.1.15:** *Let $G = \{G_1 \cup G_2 \cup \ldots \cup G_N, *_1, *_2, \ldots, *_N\}$ be a N-group semigroup. Let H and K be any two N-subgroups of the N-group semigroup G. i.e. $H = \{H_1 \cup H_2 \cup \ldots \cup H_N, *_1, \ldots, *_N\}$ where $(H_i, *_i)$ is a subgroup of $(G_i, *_i)$ and $K = \{K_1 \cup K_2 \cup \ldots \cup K_N, *_1, \ldots, *_N\}$ is such that each $(K_i, *_i)$ is again a subgroup of $(G_i, *_i)$. We say H is a N-conjugate subgroup of K if each $H_i$ is conjugate with $K_i$, i = 1, 2, …, N.*

We can also prove it in case of groups by the following result.

**THEOREM 5.1.2:** *Let $G = \{G_1 \cup G_2 \cup \ldots \cup G_N, *_1, \ldots, *_N\}$ be a N-group semigroup. Suppose $H = \{H_1 \cup \ldots \cup H_N, *_1, \ldots, *_N\}$ and $K = \{K_1 \cup K_2 \cup \ldots \cup K_N, *_1, \ldots, *_N\}$ be two sub N-groups of the N-group semigroup G. $HK = \{H_1 K_1 \cup \ldots \cup H_N K_N, *_1, \ldots, *_N\}$ will be a sub N-group if and only if $K_i H_i = H_i K_i$ for i = 1, 2, …, N.*

*Proof:* We know if G is any group, H and K are subgroups of G then we have HK to be a subgroup if and only if KH = HK. Thus we see in case of any two N-subgroup H and K of G, G a N-group semigroup we have HK = KH i.e. if and only if $H_i K_i = K_i H_i$, i = 1, 2, …, N. Hence the claim.

Now we proceed on to define the notion of Smarandache N-group semigroups.

**DEFINITION 5.1.16:** *Let* $(G = G_1 \cup G_2 \cup \ldots \cup G_N, *_1, \ldots, *_N)$ *be a non empty set on which is defined N-binary operations.* $(G, *_1, \ldots, *_N)$ *is defined as a Smarandache N-group semigroup (S-N-group semigroup) if*

i. $G = G_1 \cup G_2 \cup \ldots \cup G_N$ *where each* $G_i$ *is a proper subset of G* $(G_i \not\subseteq G_j, G_j \not\subseteq G_j, i \neq j)$.
ii. *Some of the* $(G_i, *_i)$ *are groups.*
iii. *Some of the* $(G_j, *_j)$ *are S-semigroups.*
iv. *Some of the* $(G_k, *_k)$ *are just semigroups* $(1 \leq i, j, k \leq N)$.

We now illustrate this by the following example:

***Example 5.1.15:*** Let $G = \{G_1 \cup G_2 \cup G_3 \cup G_4 \cup G_5, *_1, \ldots, *_5\}$ be a 5-group semigroup where

$G_1$ = {2Z, the set of integers under multiplication is a semigroup},
$G_2$ = S(3), the semigroup of the set of all maps from the set (1 2 3) to itself,
$G_3$ = $A_4$, the alternating group,
$G_4$ = {$Z_{12}$, the semigroup under multiplication modulo 12} and
$G_5$ = $\{g \mid g^{10} = 1\}$ the cyclic group of order 10.

Clearly G is a S-N-group semigroup for $Z_{12}$ is a S-semigroup but 2 Z is just a semigroup.

Now we can define the notion of commutativity of a S-N-group semigroup.

**DEFINITION 5.1.17:** *Let* $G = \{G_1 \cup G_2 \cup \ldots \cup G_N, *_1, \ldots, *_N\}$ *be a N-group semigroup. G be a S-N-group semigroup. We call G a Smarandache commutative N-group semigroup (S-commutative N-group semigroup) if* $G_i$ *are commutative groups,* *every proper subgroup of the S-semigroup is commutative and every semigroup which are not S-semigroups in G is also commutative.*

It is to be noted if $G = \{G_1 \cup G_2 \cup \ldots \cup G_N, *_1, \ldots, *_N\}$ is any S-N-group semigroup then G need not in general be a S-commutative N-group semigroup.

First we illustrate this by the following examples:

***Example 5.1.16:*** Let $G = \{G_1 \cup G_2 \cup G_3, *_1, *_2, *_3\}$ be a N-group semigroup, where

$G_1$ = S(3), the S-semigroup,
$G_2$ = $\{g \mid g^7 = 1\}$ and
$G_3$ = $\{Z_{15}$, the semigroup under multiplication modulo 15$\}$.

Clearly G is a S-N-group semigroup, but G is not a S-commutative N-group semigroup; for S(3) contains the subgroup $S_3$ which is non commutative.

We define Smarandache subcommutative N-group semigroup in what follows:

**DEFINITION 5.1.18:** *Let $G = \{G_1 \cup G_2 \cup \ldots \cup G_N, *_1, \ldots, *_N\}$ be a N-group semigroup. We call G a Smarandache subcommutative N-group semigroup (S-subcommutative N-group semigroup) if the following conditions are satisfied.*

- *i. G has non trivial sub N-groups.*
- *ii. Every sub N-group of G is commutative.*

We say a N-group semigroup G is Smarandache weakly subcommutative (S-weakly subcommutative) if atleast one of the sub N-group is commutative.

*Note:* Every S-subcommutative N-group semigroup is evidently S-weakly subcommutative.

***Example 5.1.17:*** Let $G = \{G_1 \cup G_2 \cup G_3 \cup G_4, *_1, *_2, *_3, *_4\}$ be a 4-group semigroup where

$G_1$ = $A_4$ the alternating group,
$G_2$ = {$Z_{12}$, semigroup under multiplication modulo 12},
$G_3$ = {Z, the semigroup under multiplication} and
$G_4$ = {$D_{2.3}$ / the dihedral group of order 6 in which $a^2 = b^3 = 1$, bab = a}.

Clearly G is a S-sub-commutative group.

*Note:* While defining the S-subcommutative group we demand if $P = \{P_1 \cup P_2 \cup \ldots \cup P_N, *_1, \ldots, *_N\}$ is a sub N-group then each $P_i$ is taken only as a proper subgroup of $G_i$. For we may have groups which are non commutative but every subgroup is commutative.

Examples of such groups are $A_4$ and $S_3$.

**DEFINITION 5.1.19:** *Let $G = \{G_1 \cup G_2 \cup \ldots \cup G_N, *_1, \ldots, *_N\}$ be a N-group semigroup. We say G is Smarandache cyclic N-group (S-cyclic N-group) if $\{G_i, *_i\}$ are cyclic groups and $\{G_j, *_j\}$ is a S-cyclic semigroup. We say G to be a Smarandache weakly cyclic N-group (S-weakly cyclic N-group) if every subgroup of $(G_i, *_i)$ are cyclic and $(G_j, *_j)$ is a S-weakly cyclic semigroup.*

Now we can define the Smarandache hyper N-group semigroup.

**DEFINITION 5.1.20:** *Suppose $G = \{G_1 \cup G_2 \cup \ldots \cup G_N, *_1, \ldots, *_N\}$ be a S-N-group. The Smarandache hyper N-group semigroup (S- hyper N-group semigroup) P is defined as follows.*

i. *If $(G_i, *_i)$ are groups then they have no subgroups.*
ii. *If $(G_j, *_j)$ is a S-semigroup then $(A_j, *_j)$ is the largest subgroup of $(G_j, *_j)$.*

iii. *If $(G_k, *_k)$ are semigroups then $(T_k, *_k)$ are the largest ideals of $(G_k, *_k)$, where*

$$P = \bigcup_{\substack{\text{over}\\ \text{relevant}\\ i}} (G_i, *_i) \bigcup_{\substack{\text{over}\\ \text{relevant}\\ j}} (A_j, *_j) \bigcup_{\substack{\text{over}\\ \text{relevant}\\ K}} (T_k, *_k).$$

*Thus $P = P_1 \cup P_2 \cup \ldots \cup P_N$ where $P_i = (G_i, *_i)$ or $P_j = (A_j, *_j)$ or $P_k = (T_k, *_k)$. $1 \leq i, j, k \leq N$.*

*We call a S-N group semigroup G to be Smarandache simple N-group semigroup (S-simple N-group semigroup) if G has no S-hyper N-group semigroup.*

Now we proceed on to study finite S-N-group semigroup.

**DEFINITION 5.1.21:** *Let $G = \{G_1 \cup G_2 \cup \ldots \cup G_N, *_1, \ldots, *_N\}$ be a S-N-group semigroup. Let $P = \{P_1 \cup P_2 \cup \ldots \cup P_N, *_1, \ldots, *_N\}$ be a subset of G. If P itself is a S-N-group semigroup then we call P a Smarandache sub N-group semigroup (S- sub N-group semigroup) of G. Now it may happen in case of finite N-group semigroup G, $o(P) / o(G)$ or $o(P) \nmid o(G)$.*

**DEFINITION 5.1.22:** *Let $G = \{G_1 \cup G_2 \cup \ldots \cup G_N, *_1, \ldots, *_N\}$ be a S-N-group semigroup of finite order. If the order of every S-sub N-group semigroup divides the order of G then we call G a Smarandache Lagrange N-group semigroup (S-Lagrange N-group semigroup). If the order of atleast one of the S-sub N-group semigroup of G say H divides the order of G then we call G a Smarandache weakly Lagrange N-group semigroup (S-weakly Lagrange N-group semigroup).*

It is easy to verify that every S-Lagrange N-group semigroup is a S-weakly Lagrange N-group semigroup.

We now proceed on to define the notion of Smarandache p-Sylow sub N-group semigroup.

**DEFINITION 5.1.23:** *Let $G = \{G_1 \cup G_2 \cup \ldots \cup G_N, *_1, \ldots, *_N\}$ be a finite S-N-group semigroup. Let p be a prime such that p divides the order of G. If G has a S-sub N-group semigroup H of* *order p or $p^t$ ($t \geq 1$), then we say G has a Smarandache p-Sylow sub N-group semigroup (S-p-Sylow sub N-group semigroup). It is very important to note that p / o(G) but $p^t \nmid$ o(G) for any t > 1, still we may have S-p Sylow sub N-group semigroup having $p^t$ elements in them. Even if p / o(G) still in case of S-N-group semigroups we may or may not have S-p Sylow sub N-group semigroups.*

We give some illustrative examples.

***Example 5.1.18:*** Let $G = \{G_1 \cup G_2 \cup G_3, *_1, *_2, *_3\}$ be a finite S-N group semigroup, where

$G_1$ = $\{Z_{24}$, is a semigroup under multiplication modulo 24$\}$,
$G_2$ = S(3) and
$G_3$ = $\langle g \mid g^{12} = 1\rangle$ be the cyclic group of order 12.

o (G) = 24 + 27 + 12 = 63.

Let $K = \{K_1 \cup K_2 \cup K_3, *_1, *_2, *_3\}$ where

$K_1$ = $\{0, 2, 4, \ldots, 22\} \subset Z_{24}$,
$K_2$ = $S_3$ and
$K_3$ = $\{g^3, g^6, g^9, 1\}$.

K is a S-N-group semigroup of G. o(K) = 12 + 6 + 4 = 22 clearly $22 \nmid 63$.

Take $T = \{T = T_1 \cup T_2 \cup T_3, *_1, *_2, *_3\}$, where

$T_1$ = $\{0, 8, 16\}$,

$$T_2 = \left\{\begin{pmatrix}1 & 2 & 3\\ 1 & 2 & 3\end{pmatrix}, \begin{pmatrix}1 & 2 & 3\\ 2 & 3 & 1\end{pmatrix}, \begin{pmatrix}1 & 2 & 3\\ 3 & 1 & 2\end{pmatrix}\right\} \text{ and}$$

$T_3$ = $\{1, g^4, g^8\}$.

o (T) = 9 and 9 / 63.

Consider $U = \{U_1 \cup U_2 \cup U_3, *_1, *_2, *_3\}$ where

$U_1$ = $\{0, 16\}$,

$U_2 = \left\{\begin{pmatrix} 1 & 2 & 3 \\ 1 & 2 & 3 \end{pmatrix}, \begin{pmatrix} 1 & 2 & 3 \\ 2 & 3 & 1 \end{pmatrix}, \begin{pmatrix} 1 & 2 & 3 \\ 3 & 1 & 2 \end{pmatrix}\right\}$ and

$U_3 = \{1, g^6\}$.

Clearly o (U) = 7 and 7 / 63.

We give yet another example.

***Example 5.1.19:*** Let G = $\{G_1 \cup G_2 \cup G_3 \cup G_4, *_1, *_2, *_3, *_4\}$ be a 4-group semigroup where

$G_1$ = $\{Z_{24}$, semigroup under multiplication modulo 24},

$G_2$ = $\{g \mid g^8 = 1\}$,

$G_3$ = $\{Z_{15}$, semigroup under multiplication modulo 15} and

$G_4$ = $A_3$, the alternating group of $S_3$.

o (G) = 24 + 8 + 15 + 3 = 50, 5 / 50 and $5^2$ / 50, 2 / 50 and $2^2 \not|$ 50 .

Take T = $\{T_1 \cup T_2 \cup T_3 \cup T_4, *_1, *_2, *_3, *_4\}$ where

$T_1$ = $\{0, 6, 12, 18\} \subset Z_4$,

$T_2$ = $\{g^2, g^4, g^6, 1\} \subset G_2$,

$T_3$ = $\{0, 3, 9, 12, 6\} \subset Z_{15}$ and

$T_4$ = $G_4 = A_3$.

o (T) = 4 + 4 + 5 + 3 = 16 = $2^4$. $2^4 \not|$ 50. But G has a S-2 Sylow sub N-group semigroup.

Take P = $P_1 \cup P_2 \cup P_3 \cup P_4$;

$P_1$ = $\{0, 2, 4, 6, 8, 10, 12, 14, 16, 18, 20, 22\} \subseteq Z_{24}$,

$P_2$ = $\{1, g^2, g^4, g^6\} \subseteq G_2$,

$P_3$ = $\{3, 9, 12, 6, 0\} \subseteq Z_{15}$ and

$P_4$ = $\{A_3\}$.

o (P) = 12 + 4 + 5 + 3 = 24.

Consider V = $\{V_1 \cup V_2 \cup V_3 \cup V_4, *_1, *_2, *_3, *_4\}$

$V_1$ = $\{0, 8, 16\} \subset Z_{24}$,
$V_2$ = $\{1, g^2, g^4, g^6\}$,
$V_3$ = $G_3 = Z_{15}$ and
$V_4$ = $A_3$.

o (V) = 3 + 4 + 15 + 3. Clearly 25 / 50 i.e. $5^2$ / 50 and G has a S-5-Sylow N-subgroup semigroup.

Now we know if S is any S-semigroup the largest proper subset which is a group is called a Smarandache hyper subgroup. For example in S (4) the S-semigroup has $S_4$ to be its S-hyper subgroup. We will call a subgroup H of a group G to be maximal if $H \subset P \subset G$ then P = H or P = G.

Now using these two notions we define S-hyper N-group semigroup.

**DEFINITION 5.1.24:** *Let $G = \{G_1 \cup G_2 \cup \ldots \cup G_N, *_1, \ldots, *_N\}$ be a S-N-group semigroup. A proper subset $P = \{P_1 \cup P_2 \cup P_3 \cup \ldots \cup P_N, *_1, \ldots, *_N\}$ is called the Smarandache P-hyper sub-N-group semigroup (S-P-hyper sub-N-group semigroup) if each $P_i \subset G_i$ is a maximal subgroup of $G_i$ and $P_i$ is a maximal subgroup of $G_i$, if $G_i$ is a group; or if $P_j \subset G_j$ is the largest group if $G_j$ is a S-semigroup.*

We illustrate this by the following example.

***Example 5.1.20:*** Let $G = \{G_1 \cup G_2 \cup G_3 \cup G_4, *_1, *_2, *_3, *_4\}$ be a S-4-group semigroup where

$G_1$ = S(3), semigroup,
$G_2$ = $\{g \mid g^{12} = 1\}$, cyclic group,
$G_3$ = {0, 1, …, 14 semigroup under multiplication modulo 15} and
$G_4$ = $D_{2.7} = \{a, b \mid a^2 = b^7 = 1; b\,a\,b = a\}$, the dihedral group of order 14.

o(G) = 27 + 12 + 15 + 14 = 68.

Take $P = \{P_1 \cup P_2 \cup P_3 \cup P_4, *_1, *_2, *_3, *_4\}$ where

$P_1$ = $S_3$,
$P_2$ = $\{g^2, g^4, g^6, g^8, g^{10}, 1\}$,
$P_3$ = {3, 6, 9, 12; this is the largest proper subset which is a group in this the element 6 acts as the multiplicative identity modulo 15} and
$P_4$ = $\{1, b, b^2, b^3, b^4, b^5, b^6\}$.

P is the S-hyper sub N-group semigroup of G.
o (P) = 6 + 6 + 4 + 7 = 23. o (G) = 27 + 12 + 15 + 14 = 68, 23 $\nmid$ 68. Thus we see in general order of a S-hyper sub N-group semigroup need not divide the order of the S-N-group semigroup.

Now we proceed on to define the notion of Smarandache Cauchy element in a finite order N-group semigroup G.

**DEFINITION 5.1.25:** *Let $G = \{G_1 \cup G_2 \cup \ldots \cup G_N, *_1, \ldots, *_N\}$be any S-N-group semigroup of finite order. Let $a \in G$; a is said to be a Smarandache special Cauchy element (S- special Cauchy element) of G if $t / o(G)$ where $t > 1$ and $a^t = 1$ (t is the least integer. It is important to note $a \in G_i$ for some $i \in \{1, 2, \ldots, N\}$).*

We can proceed on to define Smarandache right coset of a N-group semigroup.

**DEFINITION 5.1.26:** *Let $G = \{G_1 \cup G_2 \cup G_3 \cup \ldots \cup G_N, *_1, \ldots, *_N\}$ be a S-N-group semigroup. Let $H = \{H_1 \cup H_2 \cup \ldots \cup H_N, *_1, \ldots, *_N\}$ be a S-sub N group semigroup of G, we define Smarandache right coset (S-right coset) $Ha = \{H_i a \cup H_1 \cup \ldots \cup \hat{H}_i \cup \ldots \cup H_N \mid$ if $a \in G_i$ and $\hat{H}_i$ denotes $H_i$ is not present in the union; this is true for any i, i.e. whenever $a \in G_i\}$. Suppose a is in some r number of $G_i$'s (say) $G_i, G_j, \ldots, G_t, G_K$ then we have $Ha = \{H_1 \cup \ldots \cup a \cup H_j a \cup \ldots \cup H_i a \cup H_K a \cup \ldots \cup H_N\}$ if we assume $H_i$ are subgroups. If $H_i$'s are just S-semigroups then we take $H'_i$ to be the subgroup of the S-semigroup $H_i$.*

*Then $Ha = \{H'_i\, a \cup H_1 \cup \ldots \cup \hat{H}_i \cup \ldots \cup H_N$ if $a \in G_i$ and $H_i$ is the S-semigroup\}. If a is in some r, S-semigroups say $H_i$, $H_j$, … , $H_t$, $Ha = \{H_1 \cup H_2 \cup \ldots \cup H'_i\, a \cup \ldots \cup H'_j a \cup \ldots \cup H'_t\, a \cup \ldots \cup H_N \mid a \in G_i, a \in G_j, \ldots, a \in G_t\}$. If $a \in G_i$, $G_i$ a group and $a \in G_j$ a S-semigroup then $Ha = \{H_1 \cup \ldots \cup H_i\, a \cup \ldots \cup H'_j\, a \cup \ldots \cup H_N \mid H'_j$ is the subgroup of the S-semigroup $G_j\}$.*

*On similar lines we can define Smarandache left coset (S-left coset) of a S-N-group semigroup. We define Smarandache coset (S- coset) of H if Ha = a H.*

Now we illustrate this by a very simple example.

***Example 5.1.21:*** Let $G = \{G_1 \cup G_2 \cup G_3, *_1, *_2, *_3\}$ where

$G_1$ = $S_3$,
$G_2$ = $\{Z_{10}$, semigroup under multiplication modulo 10\} and
$G_3$ = $\{Z_{24}$, the semigroup under multiplication\}.

G is a 3-group semigroup.

Let $H = \{H_1 \cup H_2 \cup H_3, *_1, *_2, *_3\}$ where

$H_1$ = $A_3$,
$H_2$ = \{0, 2, 4, 6, 8\} and
$H_3$ = \{0, 8, 16\}.

Let

$$x = \begin{pmatrix} 1 & 2 & 3 \\ 1 & 3 & 2 \end{pmatrix} \in G.$$

Clearly $x \in G_1$ so $Hx = \{H_1\, x \cup H_2 \cup H_3\}$. Also x H = H x.

Now let $y = 18 \in G$, $y \in G_3$ so

$$Hy = \{H_1 \cup H_2 \cup \{8, 16\}.\, 18\}$$
$$= \{H_1 \cup H_2 \cup \{16, 0\}\} =$$
$$yH = \{H_1 \cup H_2 \cup \{16, 0\}\};$$

i.e. Hy = yH. We know given any finite N-group semigroup G and H is any sub N-group semigroup, in general $o(H) \not| \, o(G)$.

**DEFINITION 5.1.27:** *Let* $G = \{G_1 \cup G_2 \cup \ldots \cup G_N, *_1, \ldots, *_N\}$ *be a finite N-group semigroup. If the order of none of the S-sub N-group semigroups divides the order of G then we call G a Smarandache non Lagrange N-group semigroup (S-non Lagrange N-group semigroup).*

We now illustrate this by the following examples:

***Example 5.1.22:*** Let $G = \{G_1 \cup G_2 \cup G_3 \cup G_4, *_1, \ldots, *_4\}$ be a finite S-N-group semigroup where

$G_1$ = $\{Z_{12}$, semigroup under multiplication modulo 12$\}$,
$G_2$ = $S_3$ the symmetric group of degree 3,
$G_3$ = $\{Z_{10}$, semigroup under multiplication modulo 10$\}$ and
$G_4$ = $\{g \mid g^9 = 1\}$.

Clearly o (G) = 12 + 6 + 10 + 9 = 37.

The order of G is a prime so no number other than 1 and 37 will divide o (G).

Consider $H = \{H_1 \cup H_2 \cup H_3 \cup H_4, *_1, *_2, *_3, *_4\}$, the S-sub N-group semigroup of G where

$H_1$ = $\{0, 2, 4, 6, 8, 10\} \subset Z_{12}$,
$H_2$ = $A_3$,
$H_3$ = $\{0\ 5\}$, and
$H_4$ = $\{g^3, g^6, 1\}$.

Clearly o (H) = 6 + 3 + 2 + 3 = 14; $14 \nmid 37$. Thus G is a S-non Lagrange 4-group semigroup.

Now we proceed on to define the notion of Smarandache N-symmetric group semigroup.

**DEFINITION 5.1.28:** *Let* $G = \{G_1 \cup G_2 \cup \ldots \cup G_N, *_1, \ldots, *_N\}$ *be a N-group semigroup. We call G a Smarandache N-symmetric group semigroup (S-N-symmetric group semigroup) if each* $(G_i, *_i)$ *is either the symmetric group* $S_n$ *or the symmetric semigroup S (m).*

So based on this definition we can give an analogous of the Smarandache Cayley's theorem for N-group semigroup.

**THEOREM (SMARANDACHE CAYLEY'S THEOREM FOR N-GROUP SEMIGROUPS):** *Let* $G = \{G_1 \cup G_2 \cup \ldots \cup G_N, *_1, \ldots, *_N\}$ *be a S-N-group semigroup. Every S-N-group semigroup is embeddable in a suitable symmetric S- N-group semigroup.*

*Proof:* When we say suitable we mean for each group and S-semigroup we choose the suitable symmetric group and the symmetric semigroup respectively.

Let $G = \{G_1 \cup G_2 \cup \ldots \cup G_N, *_1, \ldots, *_N\}$ be any S-N-group semigroup G. We know every semigroup in G is a S-semigroup. We know a S-semigroup S is embeddable in a symmetric semigroup S(n).

So if G is a S-N group semigroup say $G = \{G_1 \cup G_2 \cup \ldots \cup G_N, *_1, \ldots, *_N\}$ then each group $G_i$ is embeddable in a suitable symmetric group $S_n$ and each S-semigroup is also embeddable in a symmetric semigroup S(n) for suitable n. Hence the claim.

Now we define the notion of double coset in a S-N-group semigroup G.

**DEFINITIONS 5.1.29:** *Let* $G = \{G_1 \cup G_2 \cup \ldots \cup G_N, *_1, \ldots, *_N\}$ *be a S-N- group semigroup. Suppose A and B be sub N groups of G. Then we define the double coset of x in G with respect to A, B if* $A x B = \{A_1 x_1 B_1 \cup \ldots \cup A_N x_N B_N \mid x = x_1 \cup x_2 \cup \ldots \cup x_N\}$.

We proceed onto define S-normal sub N-group semigroups of a S-N-group semigroup.

**DEFINITION 5.1.30:** *Let* $G = \{G_1 \cup G_2 \cup \ldots \cup G_N, *_1, \ldots, *_N\}$ *be a S-group semigroup. We call a S-sub N-group H of S-N-group semigroup G where* $A = \{A_1 \cup A_2 \cup \ldots \cup A_N, *_1, \ldots, *_N\}$ *to be a Smarandache normal sub N-group semigroup (S-normal* *sub N-group semigroup) of G, if $x \in G$ and $x = (x_1 \cup x_2 \cup \ldots \cup x_N) \in G$ then $x A = (x_1 A_1 \cup \ldots \cup x_N A_N) \subseteq A$ and $Ax = Ax_1 \cup \ldots \cup A_N x_N \subseteq A$ for all $x \in G$.*

We call a S-group semigroup to be Smarandache pseudo simple N-group if G has no S-normal sub N-group.

Now we proceed on to define the notion of the non associative N-algebraic structure viz loops and groupoids.

## 5.2 N-loop-groupoids and their properties

In this section we define a new mixed N-algebraic structure called N-loop groupoid. Clearly this algebraic structure is always non associative as both loops and groupoids are non-associative. We give several interesting properties about them. Further we give examples so that the interested reader can easily follow them.

**DEFINITION 5.2.1:** *Let $L = \{L_1 \cup L_2 \cup \ldots \cup L_N, *_1, \ldots, *_N\}$ be a non empty set with N-binary operations defined on it. We call L a N-loop groupoid if the following conditions are satisfied:*

*i. $L = L_1 \cup L_2 \cup \ldots \cup L_N$ where each $L_i$ is a proper subset of L i.e. $L_i \not\subseteq L_j$ or $L_j \not\subseteq L_i$ if $i \neq j$, for $1 \leq i, j \leq N$.*

*ii. $(L_i, *_i)$ is either a loop or a groupoid.*

*iii. There are some loops and some groupoids in the collection $\{L_1, \ldots, L_N\}$.*

*Clearly L is a non associative mixed N-structure.*

We illustrate them by the following examples.

***Example 5.2.1:*** Let $L = \{L_1 \cup L_2 \cup L_3, *_1, *_2, *_3\}$ where $L_1 = \{Z_{10}$, with $*_1$ defined on $Z_{10}$ as $a *_1 b = 2a + 3b \pmod{10}$ where $a, b \in Z_{10}\}$, $L_1$ is a groupoid,

$L_2$ is the loop given by the following table:

| $*_2$ | e | $a_1$ | $a_2$ | $a_3$ | $a_4$ | $a_5$ |
|---|---|---|---|---|---|---|
| e | e | $a_1$ | $a_2$ | $a_3$ | $a_4$ | $a_5$ |
| $a_1$ | $a_1$ | e | $a_3$ | $a_5$ | $a_2$ | $a_4$ |
| $a_2$ | $a_2$ | $a_5$ | e | $a_4$ | $a_1$ | $a_3$ |
| $a_3$ | $a_3$ | $a_1$ | $a_1$ | e | $a_5$ | $a_2$ |
| $a_4$ | $a_4$ | $a_3$ | $a_5$ | $a_2$ | e | $a_1$ |
| $a_5$ | $a_5$ | $a_2$ | $a_4$ | $a_1$ | $a_3$ | e |

$L_3$ is the loop given by the following table:

| $x_3$ | e | a | b | c | d |
|---|---|---|---|---|---|
| e | e | a | b | c | d |
| a | a | e | c | d | b |
| b | b | d | a | e | c |
| c | c | b | d | a | e |
| d | d | c | e | b | a |

Now L is a 3-loop groupoid.

We say the N-loop groupoid $L = \{L_1 \cup L_2 \cup \ldots \cup L_N, *_1, *_2, \ldots, *_N\}$ is finite if the number of distinct elements in L is finite and it is called as the order of L that is the N-loop groupoid is of finite order. If the number of elements in L is infinite we say L is of infinite order and L is an infinite N-loop groupoid.

**DEFINITION 5.2.2:** *Let $L = \{L_1 \cup L_2 \cup \ldots \cup L_N, *_1, \ldots, *_N\}$ be a N-loop groupoid. L is said to be a commutative N-loop groupoid if each of $\{L_i, *_i\}$ is commutative.*

Now we give an example of a commutative N-loop groupoid.

***Example 5.2.2:*** Let $L = \{L_1 \cup L_2 \cup L_3, *_1, *_2, *_3\}$ be 3-loop groupoid where

$L_1$ is given by the following table:

| $*_1$ | $g_1$ | $g_2$ | $g_3$ | $g_4$ | $g_5$ | $g_6$ | $g_7$ |
|---|---|---|---|---|---|---|---|
| $g_1$ | $g_1$ | $g_4$ | $g_7$ | $g_3$ | $g_6$ | $g_2$ | $g_5$ |
| $g_2$ | $g_6$ | $g_2$ | $g_5$ | $g_1$ | $g_4$ | $g_7$ | $g_3$ |
| $g_3$ | $g_4$ | $g_7$ | $g_3$ | $g_6$ | $g_2$ | $g_5$ | $g_1$ |
| $g_4$ | $g_2$ | $g_5$ | $g_1$ | $g_4$ | $g_7$ | $g_3$ | $g_6$ |
| $g_5$ | $g_7$ | $g_3$ | $g_6$ | $g_2$ | $g_5$ | $g_1$ | $g_4$ |
| $g_6$ | $g_5$ | $g_1$ | $g_4$ | $g_7$ | $g_3$ | $g_6$ | $g_2$ |
| $g_7$ | $g_3$ | $g_6$ | $g_2$ | $g_5$ | $g_1$ | $g_4$ | $g_7$ |

It is easily verified that $L_1$is a groupoid of order 7. Now we give the table of the loop ($L_2$, $*_2$):

| $*_2$ | e | $a_1$ | $a_2$ | $a_3$ | $a_4$ | $a_5$ | $a_6$ | $a_7$ |
|---|---|---|---|---|---|---|---|---|
| e | e | $a_1$ | $a_2$ | $a_3$ | $a_4$ | $a_5$ | e | $a_7$ |
| $a_1$ | $a_1$ | e | $a_4$ | $a_7$ | $a_3$ | $a_6$ | $a_2$ | $a_5$ |
| $a_2$ | $a_2$ | $a_6$ | e | $a_5$ | $a_1$ | $a_4$ | $a_7$ | $a_3$ |
| $a_3$ | $a_3$ | $a_4$ | $a_7$ | e | $a_6$ | $a_2$ | $a_5$ | $a_1$ |
| $a_4$ | $a_4$ | $a_2$ | $a_5$ | $a_1$ | e | $a_7$ | $a_3$ | $a_6$ |
| $a_5$ | $a_5$ | $a_7$ | $a_3$ | $a_6$ | $a_2$ | e | $a_1$ | $a_4$ |
| $a_6$ | $a_6$ | $a_5$ | $a_1$ | $a_4$ | $a_7$ | $a_3$ | e | $a_2$ |
| $a_7$ | $a_7$ | $a_3$ | $a_6$ | $a_2$ | $a_5$ | $a_1$ | $a_4$ | e |

Clearly ($L_2$, $*_2$) is a commutative loop of order 8. Now we give the table for the commutative loop of order 6.

| $*_3$ | e | 1 | 2 | 3 | 4 | 5 |
|---|---|---|---|---|---|---|
| e | e | 1 | 2 | 3 | 4 | 5 |
| 1 | 1 | e | 4 | 2 | 5 | 3 |
| 2 | 2 | 4 | e | 5 | 3 | 1 |
| 3 | 3 | 2 | 5 | e | 1 | 4 |
| 4 | 4 | 5 | 3 | 1 | e | 2 |
| 5 | 5 | 3 | 1 | 4 | 2 | e |

Let L = {$L_1 \cup L_2 \cup L_3$, $*_1$, $*_2$, $*_3$} is a 3-commutative loop groupoid and the order of L is 7 + 8 + 6 = 21. Clearly L is a finite 3-loop groupoid.

**DEFINITION 5.2.3:** *Let L = {$L_1 \cup L_2 \cup \ldots \cup L_N$, $*_1$, …, $*_N$} be a N-loop groupoid. A proper subset P of L is said to be a sub N-loop groupoid if P = {$P_1 \cup P_2 \cup \ldots \cup P_N$, $*_1$, …, $*_N$} be a N-loop groupoid. A proper subset P = {$P_1 \cup P_2 \cup \ldots \cup P_N$, $*_1$, …, $*_N$} is such that if P itself is a N-loop groupoid then we call P the sub N-loop groupoid of L.*

***Example 5.2.3:*** Let L = {$L_1 \cup L_2 \cup L_3$, $*_1$, $*_2$, $*_3$} where $L_1 = P \times L$; here P and L are loops given by the following tables:

| * | e | a | b | c | d |
|---|---|---|---|---|---|
| e | e | a | b | c | d |
| a | a | e | c | d | b |
| b | b | d | a | e | c |
| c | c | b | d | a | e |
| d | d | c | e | b | a |

and

| *' | e | $a_1$ | $a_2$ | $a_3$ | $a_4$ | $a_5$ |
|---|---|---|---|---|---|---|
| e | e | $a_1$ | $a_2$ | $a_3$ | $a_4$ | $a_5$ |
| $a_1$ | $a_1$ | e | $a_5$ | $a_4$ | $a_3$ | $a_2$ |
| $a_2$ | $a_2$ | $a_3$ | e | $a_1$ | $a_5$ | $a_4$ |
| $a_3$ | $a_3$ | $a_5$ | $a_4$ | e | $a_2$ | $a_1$ |
| $a_4$ | $a_4$ | $a_2$ | $a_1$ | $a_5$ | e | $a_3$ |
| $a_5$ | $a_5$ | $a_5$ | $a_4$ | $a_2$ | $a_1$ | e |

where P = {e, a, b, c, d} and L = {e $a_1$, $a_2$, $a_3$, $a_4$, $a_5$}. $L_2$ is the groupoid given by $Z_{12}$ = {a * b = 2a + 4 b (mod 12) for a, b $\in Z_{12}$} and $L_3$ = {$Z_5 \times Z_{10}$} where the operation in $Z_5$ is {a * b = a + 3b (mod 5)} and in $Z_{10}$ is {a $\otimes$ b = 2 a + b (mod 10)}. Clearly L is a 3-loop groupoid. Take $P = P_1 \cup P_2 \cup P_3$ where $P_1$ = {P ×

{e}} is a loop $P_2$ = {a, b ∈ {0, 2, 4, 6, 8, 10} ⊂ $Z_{12}$ and $P_3$ = {$Z_5$ × {ϕ}} P is a sub3-loop groupoid of L.

Next we define sub N group of a N-loop groupoid.

**DEFINITION 5.2.4:** *Let L = {$L_1 \cup L_2 \cup \ldots \cup L_N$, $*_1, \ldots, *_N$} be a N-loop groupoid. A proper subset G = {$G_1 \cup G_2 \cup \ldots \cup G_N$, $*_1, \ldots, *_N$} is called a sub N-group if each ($G_i$, $*_i$) is a group.*

Now we illustrate this by the following example:

***Example 5.2.4:*** L = $L_1 \cup L_2 \cup L_3$ be a 3-loop groupoid where $L_1$ is the loop given by the following table:

| $*_1$ | e' | a | b | c | d |
|---|---|---|---|---|---|
| e' | e' | a | b | c | d |
| a | a | e' | c | d | b |
| b | b | d | a | e' | c |
| c | c | b | d | a | e' |
| d | d | c | e' | b | a |

i.e.

$L_1$ = {e', a, b, c, d},

$L_2$ = {$Z_{24}$ | a $*_1$ b = 2a + 4b (mod 24) for all a, b, ∈ $Z_4$} is a groupoid and

$L_3$ is the loop given by the table:

| $*_3$ | e | $g_1$ | $g_2$ | $g_3$ | $g_4$ | $g_5$ |
|---|---|---|---|---|---|---|
| e | e | $g_1$ | $g_2$ | $g_3$ | $g_4$ | $g_5$ |
| $g_1$ | $g_1$ | e | $g_4$ | $g_2$ | $g_5$ | $g_3$ |
| $g_2$ | $g_2$ | $g_4$ | e | $g_5$ | $g_3$ | $g_1$ |
| $g_3$ | $g_3$ | $g_2$ | $g_5$ | e | $g_1$ | $g_4$ |
| $g_4$ | $g_4$ | $g_5$ | $g_3$ | $g_1$ | e | $g_2$ |
| $g_5$ | $g_5$ | $g_3$ | $g_1$ | $g_4$ | $g_2$ | e |

Now take P = $\{P_1 \cup P_2 \cup P_3, *_1, *_2, *_3\}$ where

$P_1$ = {e', a},
$P_2$ = {0, 6, 3, 12, 18, 9} under '+' and
$P_3$ = $\{e, g_1\}$.

Clearly P is sub N-group of the N-loop groupoid L.

Now we proceed on to define sub N-groupoid.

**DEFINITION 5.2.5:** *Let $L = \{L_1 \cup L_2 \cup \ldots \cup L_N, *_1, \ldots, *_N\}$ be a N-loop groupoid. A proper subset $T = \{T_1 \cup T_2 \cup \ldots \cup T_N, *_1, \ldots, *_N\}$ is said to be a sub N-groupoid of the N-loop groupoid if each $(T_i, *_i)$ is a groupoid.*

It is important to say given any N-loop groupoid we are not always guaranteed of a sub N-group or sub N-groupoid. It may or may not exist.

Next we proceed on to define sub N-loop of a N-loop groupoid.

**DEFINITION 5.2.6:** *Let $L = \{L_1 \cup L_2 \cup \ldots \cup L_N, *_1, \ldots, *_N\}$ be a N-loop groupoid. A non empty subset $S = \{S_1 \cup S_2 \cup \ldots \cup S_N, *_1, \ldots, *_N\}$ is said to be a sub N-loop if each $\{S_i, *_i\}$ is a loop.*

Next we define sub N-semigroup.

**DEFINITION 5.2.7:** *Let $L = \{L_1 \cup L_2 \cup \ldots \cup L_N, *_1, \ldots, *_N\}$ be a N-loop groupoid. A non empty subset $W = \{W_1 \cup W_2 \cup \ldots \cup W_N, *_1, \ldots, *_N\}$ of L said to be a sub N-semigroup if each $\{W_i, *_i\}$ is a semigroup.*

Next we define the notion of sub N-group semigroup.

**DEFINITION 5.2.8:** *Let $L = \{L_1 \cup L_2 \cup \ldots \cup L_N, *_1, \ldots, *_N\}$ be a N-loop groupoid. Let $R = \{R_1 \cup R_2 \cup \ldots \cup R_N, *_1, \ldots, *_N\}$be a proper subset of L. We call R a sub N-group groupoid of the N-loop groupoid L if each $\{R_i, *_i\}$ is either a group or a groupoid.*

Thus one of the main advantage of these mixed structures is that we can define many mixed substructures of varying varieties given one mixed structure.

We can give several illustrative examples of each. Now we give just an example of a sub-N-group groupoid of a N-loop groupoid.

***Example 5.2.5:*** Let L = $\{L_1 \cup L_2 \cup L_3, *_1, *_2, *_3\}$ be a N-loop groupoid. Take N = 3, where

$L_1$ = $\{L_7(4)$ is a commutative loop of order 8 given in page 19$\}$,
$L_2$ = $\{Z_{12}$, where for a, b $\in Z_{12}$; a $*_2$ b = 3a + 6b mod (12)$\}$ and
$L_3$ = $\{Z_{14}$, where for a, b $\in Z_{14}$; a $*_3$ b = 3a + 4b (mod 14)$\}$.

Consider the subset K = $K_1 \cup K_2 \cup K_3, *_1, *_2, *_3\}$ of L.

$\{K_1, *_1\}$ = $\{e, g_4, *_1\}$ is a group,
$\{K_2, *_2\}$ = $\{0, 3, 6, 9, *_2\}$ is a subgroupoid and
$\{K_3, *_3\}$ = $\{7, 0\}$ is a subgroupoid.

Thus K is a sub N group groupoid of L.

Having defined several substructures we will give more weightage only to the algebraic sub N-structure which is the same as the given N-structure. It is still important to see that these sub N-structures are such that they may be of any order. The classical Lagrange's theorem may not be in general true for all finite order N-loop groupoid.

First we illustrate this by an example.

***Example 5.2.6:*** Let L = $\{L_1 \cup L_2 \cup L_3, *_1, *_2, *_3\}$ be a 3-loop groupoid of finite order where

$L_1$ = $\{L_7(4)$ the commutative loop of order 8 given in page 19$\}$,
$L_2$ = $\{Z_{10}$, for a, b $\in Z_{10}$, with a * b = 2a + 3b (mod 10)$\}$ and

$L_3 = \{Z_{24} \mid \text{for } a, b \in Z_{24};\ a * b = 8a + 4b \pmod{24}\}$.

order of L = 8 + 10 + 24 = 42.

Let $K = \{K_1 \cup K_2 \cup K_3, *_1, *_2, *_3\}$ be a sub N-loop groupoid where

$K_1 = \{e, g_1\}$,
$K_2 = \{0, 2, 4, 6, 8\}$ is a subgroupoid of $(L_2, *_2)$ and
$K_3 = \{0, 2, 4, 6, 8, 10, 12, 14, 16, 18, 20, 22\} \subset Z_{24}$.

o(K) = 2 + 5 + 12 = 19 we see $19 \nmid 42$.

Thus a sub N loop groupoid in general does not divide the order of the N-loop groupoid which is of order 42.

Let $V = \{V_1 \cup V_2 \cup V_3, *_1, *_2, *_3\}$ where

$V_1 = L_1$,
$V_2 = \{0, 2, 4, 6, 8\} \subset Z_{10}$ and
$V_3 = \{0, 6, 12, 18\} \subseteq Z_4$ is a sub N-loop groupoid,

where o (V) = 8 + 5 + 4 = 17, $17 \nmid 42$. Thus even the order of the sub N-loop groupoid does not divide the order of the N-loop groupoid L.

Now we proceed on to define the notion of Lagrange N-loop groupoid.

**DEFINITION 5.2.9:** *Let $L = \{L_1 \cup L_2 \cup \ldots \cup L_N, *_1, \ldots, *_N\}$ be a N-loop groupoid of finite order. $K = \{K_1 \cup K_2 \cup K_3, *_1, \ldots, *_N\}$ be a sub N-loop groupoid of L. If every sub N-loop groupoid divides the order of the N-loop groupoid L then we call L a Lagrange N-loop groupoid. If no sub N-loop groupoid of L divides the order of L then we say L is a Lagrange free N-loop groupoid.*

Let us give an example of a Lagrange free N-loop groupoid.

***Example 5.2.7:*** Let $L = \{L_1 \cup L_2 \cup L_3, *_1, *_2, *_3\}$ be a N-loop groupoid where

$L_1$ = {$L_7(4)$ is the loop of order 8 given in page 19},
$L_2$ = {$Z_9$ | a * b = 3a + 4b (mod 9)} and
$L_3$ = {$Z_{14}$ | a * b = 2a + 5b (mod 14)}.

Clearly order of L, is 8 + 9 + 14 = 31. We see whatever be the sub N-loop groupoid none will divide 31 since 31 is a prime.

Next we define the notion of weakly Lagrange N-loop groupoid.

**DEFINITION 5.2.10:** *Let L = {$L_1 \cup L_2 \cup \ldots \cup L_N$, $*_1, \ldots, *_N$} be a finite N-loop groupoid. Suppose L contains atleast one proper sub N-loop groupoid say P such that o(P) / o(L) then we call L a weakly Lagrange N-loop groupoid.*

It can be easily verified that every Lagrange N-loop groupoid is a weakly Lagrange N-loop groupoid. For a N-loop groupoid L to be Lagrange or weakly Lagrange we need the order of L to be a composite number. All N-loop groupoids of finite order whose order is a prime is a Lagrange free N-loop groupoid. It is important to mention here that even N-loop groupoid of prime order can have sub N-loop groupoid.

***Example 5.2.8:*** Let L = {$L_1 \cup L_2 \cup L_3$, $*_1$, $*_2$, $*_3$} be 3-loop groupoid where

$L_1$ = {$L_7(3)$ is the commutative loop of order 8 in page 19},
$L_2$ = {$Z_{12}$ | a * b = 2a + 4b (mod 12)} and
$L_3$ = {$Z_{16}$ | a * b = 3a + 5b (mod 16)}.

$L_2$ and $L_3$ are groupoids. Now order of L = 8 + 12 + 16 = 36.

Let us find some of its sub N-loop groupoids. Let V = {$V_1 \cup V_2 \cup V_3$, $*_1$, $*_2$, $*_3$} where

$V_1$ = {e, $g_2$} is a subloop of $L_1$,
$V_2$ = {0, 2, 4, 6, 8, 10} $\subset Z_{12}$ is a subgroupoid and
$V_3$ = {0, 2, 4, 6, 8, 10, 12, 14} $\subset Z_{16}$ is a subgroupoid.

Now $o(V) = 2 + 6 + 8 = 16$ i.e. $o(V) \nmid o(L)$ i.e. $16 \nmid 36$. Now consider the sub N-loop groupoid $P = \{P_1 \cup P_2 \cup P_3, *_1, *_2, *_3\}$ where

$P_1 = \{e, g_1\}$,
$P_2 = \{0, 8, 4\}$ and
$P_3 = V_3$.

Now $o(P) = 2 + 3 + 8 = 13$. i.e. $o(P) \nmid o(L)$ so $13 \nmid 36$. Thus L is not a Lagrange N-loop groupoid. Take $S = \{S_1 \cup S_2 \cup S_3, *_1, *_2, *_3\}$ where

$S_1 = \{e, g_4\}$,
$S_2 = \{0, 8, 4\}$ and
$S_3 = \{0, 4, 8, 12\}$.

Now $o(S) = 2 + 3 + 4 = 9$, $o(S) / o(L)$ i.e. $9 / 36$. Thus L is a weakly Lagrange N-loop groupoid, so L is not a Lagrange free or Lagrange N-loop groupoid.

***Example 5.2.9:*** Let $L = \{L_1 \cup L_2 \cup L_3, *_1, *_2, *_3\}$ be a N-loop groupoid where

$L_1 = \{$the loop $L_5(3)\}$ a loop of order 6},
$L_2 = \{Z_6 \mid a * b = 2a + 4b \pmod 6\}$ for $a, b \in Z_{16}$ and
$L_3 = \{Z_4 \mid a * b = 3a + b \pmod 4$ for $a,b \in Z_4\}$ are groupoids.

$o(L) = 6 + 6 + 4 = 16 = 2^4$ . Take $V = \{V_1 \cup V_2 \cup V_3, *_1, *_2, *_3\}$ where

$V_1 = \{e, g_3\}$,
$V_2 = \{0, 2, 4\}$ and
$V_3 = \{0, 2\}$.

$o(V) = 7$. $7 \nmid 16$. Consider $P = \{P_1 \cup P_2 \cup P_3, *_1, *_2, *_3\}$ where

$P_1 = \{e, g_3\}$,

$P_2$ = {0, 2, 4, 3} and
$P_3$ = {0, 2}.

o (P) = 8 that is o (P) / o(L). So L is only a weakly Lagrange N loop groupoid.

***Example 5.2.10:*** Let L = {$L_1 \cup L_2 \cup L_3 \cup L_4, *_1, *_2, *_3, *_4$} be a 4-loop groupoid where

$L_1$ = $L_5(2)$,
$L_2$ = $L_7(4)$ be loops or order 6 and 8 respectively,
$L_3$ = {$Z_5$ | for a, b $\in Z_5$, a * b = 2 a + 3b (mod 5)} and
$L_4$ = {$Z_4$ | for a, b $\in Z_4$, a * b = 2a + b (mod 9)}.

be groupoids of order 5 and 4 respectively.

Now o (L) = 6 + 8 + 5 + 4 = 23. Clearly L cannot have any proper sub N-loop groupoid P such that o (P) / o(L). Thus L is a Lagrange free N-loop groupoid.

***Remark****:* It is interesting to note all prime order N-loop groupoids are Lagrange free N-loop groupoids. Still it is pertinent to note that even if order of a N-loop groupoid is a prime still it can have non trivial sub N loop groupoids. But we call them Lagrange free as no number other than itself can divide a prime number. Now we proceed on to the next concept of Sylow sub N- loop groupoids.

Now we see in case of N-loop groupoids L of finite order, even if $p^{\alpha}$/ o(L) but $p^{\alpha+1} \not|$ o(L) where p is a prime, L may not have sub N-loop groupoids of order $p^{\alpha}$.

**DEFINITION 5.2.11:** *Let L = {$L_1 \cup L_2 \cup \ldots \cup L_N, *_1, \ldots, *_N$} be a N-loop groupoid of finite order. If p is a prime such that $p^{\alpha}$/ o(L) but $p^{\alpha+1} \not|$ o(L) and if $L_t$ has a sub N-loop groupoid of order $p^{\alpha}$ then we call P the p Sylow sub N-loop groupoid of L.*

*If L is a N-loop groupoid of finite order; L has for every prime p, p / o(L) and $p^{\alpha}$ / o(L) but $p^{\alpha+1} \not|$ o(L) and if L has a p*

*Sylow sub N-loop groupoids then we call L a Sylow N-loop groupoid.*

We illustrate this by the following example:

***Example 5.2.11:*** Let L = $\{L_1 \cup L_2 \cup L_3, *_1, *_2, *_3\}$ be a 3-loop groupoid of finite order where

$L_1$ = $L_5(3)$, a loop of order 6,
$L_2$ = {a, b $\in Z_6$ such that a * b = 4a + 2b (mod 6)} and
$L_3$ = {For a, b $\in Z_8$ define a * b = 6, a + 2b (mod 8)} be groupoids.

Clearly o (L) = 6 + 6 + 8 = 20. 5 / 20 and $5^2 \nmid 20$, 2/20 and $2^2$/ 20 and $2^3 \nmid 20$. To find whether L has sub N-loop groupoids of order 5 and 4.

Let V = $\{V_1 \cup V_2 \cup V_3, *_1, *_2, *_3\}$ where

$V_1$ = {e},
$V_2$ = {0, 4} and
$V_3$ = {0, 3}.
o (V) = 1 + 2 + 2 = 5.

Thus V is a 5 Sylow sub N-loop groupoid. Yet we cannot have a sub N-loop groupoid of order 4. For both groupoids or a loop and groupoid must be an identity.

Consider T = $\{T_1 \cup T_2 \cup T_3, *_1, *_2, *_3\}$, be a sub N-loop groupoid where

$T_1$ = $\{e, g_1\}$,
$T_2$ = {0, 3} and
$T_3$ = $\{0, 2, 4, 6\} \subset Z_8$.

o(T) = 2 +2 + 4 = 8 = $2^3$; $2^2$ / o(L) but L has a sub N-loop groupoid of order $2^3$.

Thus we make some more modified definitions in this case.

**DEFINITION 5.2.12:** *Let $L = \{L_1 \cup L_2 \cup \ldots \cup L_N, *_1, \ldots, *_N\}$ be a N-loop groupoid of finite order. If for at least one prime p such that $p^{\alpha}/ o(L)$ but $p^{\alpha+1} \not| \; o(L)$ we have a sub N-loop groupoid of order $p^{\alpha}$ then we call L a weakly Sylow N-loop groupoid.*

**DEFINITION 5.2.13:** *Let $L = \{L_1 \cup L_2 \cup \ldots \cup L_N, *_1,\ldots, *_N\}$ be a N-loop groupoid of finite order. If for all primes p such that p / o(L) we have non trivial sub N-loop groupoid of order $p^t$ ($t \geq 1$) then we call L a semi Sylow N-loop groupoid.*

*If this happens only for one prime p then we call L a weakly semi Sylow N-loop groupoid.*

It is interesting to note all semi Sylow N-loop groupoids are trivially weakly semi-Sylow N-loop groupoid. Also all Sylow N-loop groupoids are weakly Sylow N-loop groupoids. We give some examples before we proceed on to define some more properties about N-loop groupoids.

***Example 5.2.12:*** Let $L = \{L_1 \cup L_2 \cup L_3, *_1, *_2, *_3\}$ be a 3-loop groupoid where

$L_1$ = $\{Z_{12}$; for a, b $\in Z_{12}$, a * b = 2a + 4b mod (12)},
$L_2$ = $\{Z_4$ | for all a, b $\in Z_4$ a * b = a + 3b (mod 4)} are groupoids and
$L_3$ = $L_7(4)$ a loop of order 8.

Clearly o(L) = 12 + 4 + 8 = 24. Now we find what type of N-loop groupoid L has.

Let $V = \{V_1 \cup V_2 \cup V_3, *_1, *_2, *_3\}$ where

$V_1$ = {0, 3, 6, 9},
$V_2$ = {0, 2} and
$V_3$ = $\{e, g_3\}$.

Clearly o(V) = 4 + 2 + 2 = 8, 8/24 so L has a 2-Sylow sub N-loop groupoid. Clearly L cannot have a nontrivial 3-Sylow sub

N-loop groupoid. So at the out set L is only a weakly Sylow N-loop groupoid.

Now we proceed on to define the notion of super Sylow N-loop groupoid.

**DEFINITION 5.2.14:** *Let $L = \{L_1 \cup L_2 \cup \ldots \cup L_N, *_1, \ldots, *_N\}$be a N-loop groupoid of finite order. Suppose L is a Sylow N-loop groupoid i.e. for every prime p such that $p^{\alpha}$ / o(L) but $p^{\alpha+1} \not| \; o(L)$; L has a sub- N-loop groupoid of order $p^{\alpha}$. If in addition for every prime p we have a sub N-loop groupoid of order $p^{\alpha+t}$ $(t \geq 1)$ then we call the Sylow N-loop groupoid to be a super Sylow N-loop groupoid.*

It can be easily verified that all super Sylow N-loop groupoids are Sylow N-loop-groupoids, but in general a Sylow N-loop groupoid need not be a super Sylow N-loop groupoid.

We proceed on to define the notion of Cauchy element in a N-loop groupoid of finite order.

**DEFINITION 5.2.15:** *Let $L = \{L_1 \cup L_2 \cup \ldots \cup L_N, *_1, \ldots, *_N\}$ be a N-loop groupoid of finite order. Let $a \in L$ if there is a positive integer say m such that $a^m = e$ and if m / o(L) then we call a the Cauchy element of L.*

*Note:* It is important to observe that in general all elements in L such that $a^m = e$ need not be a Cauchy element. Secondly it is also important to note all elements in L need not be such that $a^t = e$ for some integer t for L can have elements x that $x^n = x$ or $x^m = 0$.

***Example 5.2.13*:** Let $L = \{L_1 \cup L_2 \cup L_{3\,/}\, *_1, *_2, *_3\}$ be a 3-loop groupoid of finite order where

$L_1$ = $L_5$ (3) be a loop of order 6,
$L_2$ = $\{Z_6$, such that for a, b $\in Z_6$, define a * b = 4a +2b (mod 6) and
$L_3$ = {For all a, b $\in Z_7$ such that a * b = 3a + 4b (mod 7)},

where $L_2$ and $L_3$ are groupoids of finite order.

Now o(L) = 6 + 6 + 7 = 19 i.e. L is a N-loop groupoid of order 19, a prime. Now we see every element in $L_1 \setminus \{e\}$ is such that $g_i^2 = e$ but $2 \nmid 19$. L cannot have any Cauchy element though L is of finite order.

Thus we can make the following. Let $L = \{L_1 \cup L_2 \cup \ldots \cup L_N, *_1, \ldots, *_N\}$ be a N-loop groupoid of a prime order then L has no Cauchy element. L is said to be "Cauchy free". It may so happen that the N-loop groupoid is of order n, n a composite number still we may not have order of every element of finite order of L to divide the order of L.

***Example 5.2.14:*** $L = \{L_1 \cup L_2 \cup L_3 *_1, *_2, *_3\}$ be a 3-loop groupoid of finite order where

$L_1$ = $\{L_7 (3)\}$ is a loop of order 8,
$L_2$ = $\{Z_{12} \mid$ for a, b $\in Z_{12}$, a * b = 6a + 8 b (mod 12)$\}$ and
$L_3$ = $\{Z_5$ / for a, b $\in Z_5$, a * b = 4a + b (mod 5)$\}$.

Clearly o (L) = 8 + 12 + 5 = 25.

Now we see the order of L is not a prime still we cannot say every element x in L such that $x^m = e$ is such that m / o(L), for take $g_i \in L_2 \setminus \{e\}$, each $g_i$ is such that $g_i^2 = e$ but $2 \nmid 25$. Hence the claim.

Now we proceed on to define the notion of Moufang N-loop groupoid.

**DEFINITION 5.2.16**: *Let $L = \{L_1 \cup L_2 \cup \ldots \cup L_N, *_1, \ldots, *_N\}$ be a N-loop groupoid. We call L a Moufang N-loop groupoid if each $(L_i, *_i)$ satisfies the following identities:*

i. *(xy) (zx) = (x (yz))x.*
ii. *((x y)z)y = x (y (2y)).*
iii. *x (y (xz)) = (xy)x)z for x, y, z $\in L_i$, $1 \leq i \leq N$.*

*Thus for a N-loop groupoid to be Moufang both the loops and the groupoids must satisfy the Moufang identity.*

Now we proceed on to define Bruck N-loop groupoid.

**DEFINITION 5.2.17:** *A N-loop groupoid; $L = \{L_1 \cup L_2 \cup \ldots \cup L_N, *_1, \ldots, *_N\}$ is called a Bruck N-loop groupoid if $\{L_i, *_i\}$ are either Bruck loop $(xyz)\, z = x\, (y\, (xz))$ and $(xy)^{-1} = x^{-1}\, y^{-1}$ for all $x, y, z \in L_i$; $i = 1, 2, \ldots, N$. The groupoid need only satisfy $(x(yz)z) = x(y(xz))$ for all $x, y, z \in L_i$.*

On similar lines we can define Bol N-loop groupoids, WIP – N-loop groupoid or alternative N-loop groupoid. We call a non empty proper subset to be normal sub N-loop groupoid if the following conditions are satisfied.

**DEFINITION 5.2.18:** *Let $L = \{L_1 \cup L_2 \cup \ldots \cup L_N, *_1, \ldots, *_N\}$ be a N-loop groupoid. A proper subset P ($P = P_1 \cup P_2 \cup \ldots \cup P_N, *_1, \ldots, *_N$) of L is a normal sub N-loop groupoid of L if*

1. *If P is a sub N-loop groupoid of L.*
2. *$x_i\, P_i = P_i\, x_i$ (where $P_i = P \cap L_i$)*
3. *$y_i\, (x_i\, P_i) = (y_i\, x_i)\, P_i$ for all $x_i, y_i \in L_i$.*

*This is true for each $P_i$, i.e., for $i = 1, 2, \ldots, N$.*

Now we can define homomorphism of N-loop groupoids.

**DEFINITION 5.2.19:** *Let $L = \{L_1 \cup L_2 \cup \ldots \cup L_N, *_1, \ldots, *_N\}$ and $K = \{K_1 \cup K_2 \cup \ldots \cup K_N, *_1, \ldots, *_N\}$ be two N-loop groupoids such that if $(L_i, *_i)$ is a groupoid then $\{K_i, *_i\}$ is also a groupoid. Likewise if $(L_j, *_j)$ is a loop then $(K_j, *_j)$ is also a loop true for $1 \leq i, j \leq N$. A map $\theta = \theta_1 \cup \theta_2 \cup \ldots \cup \theta_N$ from L to K is a N-loop groupoid homomorphism if each $\theta_i : L_i \to K_i$ is a groupoid homomorphism and $\theta_j : L_j \to K_j$ is a loop homomorphism $1 \leq i, j \leq N$.*

Thus N-loop groupoid homomorphisms are defined on very special conditions only. We cannot think of picking and mapping for it may create problems of undefinedness, however if any such maps be sought we do not call them as a N-loop groupoid homomorphisms but we call them as pseudo homomorphisms. This pseudo homomorphisms can take place even when N are different.

Now we proceed on to define pseudo homomorphisms of N-loop groupoids and M-loop groupoids. ($M \neq N$)

**DEFINITION 5.2.20:** *Let $L = \{L_1 \cup L_2 \cup \ldots \cup L_N, *_1, \ldots, *_N\}$ be a N-loop groupoid and $K = \{K_1 \cup K_2 \cup \ldots \cup K_M, *_1, \ldots, *_M\}$ be a M-loop groupoid. A map $\phi = \phi_1 \cup \phi_2 \cup \ldots \cup \phi_{N'}$ from L to K is called a pseudo N-M-loop groupoid homomorphism if each $\phi_i$: $L_t \rightarrow K_s$ is either a loop homomorphism or a groupoid homomorphism, $1 \leq t \leq N$ and $1 \leq s \leq M$, according as $L_t$ and $K_s$ are loops or groupoids respectively (we demand $N \leq M$ for if $M > N$ we have to map two or more $L_i$ onto a single $K_j$ which can not be achieved easily).*

Next we proceed on to define the notion of Smarandache N-loop groupoids as this structure has several sub N-structures we have several types of Smarandache N-loop groupoids, which we would be defining.

Now we proceed on to define several types of Smarandache N-loop groupoids and give some of its properties. It is in fact pertinent to state one can easily derive and prove several properties but we leave for the reader to develop such concepts as our motivation is mainly to attract young researchers in this subject.

**DEFINITION 5.2.21:** *Let $L = \{L_1 \cup L_2 \cup \ldots \cup L_N, *_1, \ldots, *_N\}$be a N-loop groupoid. We call L a Smarandache loop N-loop groupoid (S-loop N-loop groupoid) if L has a proper subset $P = \{P_1 \cup P_2 \cup \ldots \cup P_N, *_1, \ldots, *_N\}$ such that each $P_i$ is a loop i.e. P is a N-loop.*

We illustrate this by the following example:

***Example 5.2.15:*** Let L = $\{L_1 \cup L_2 \cup L_3, *_1, *_2, *_3\}$ be a N-loop groupoid where

$L_1$ = $L_{15}(2)$ is a loop of order 16 and (e, 2, 5, 8, 14) is a subloop in $L_{15}(2)$.

$L_2$ = $L_{45}(8)$, the loop of order 46 where {e, 1, 6, 11, 16, 21, 26, 31, 36, 41} is a subloop in $L_{45}(8)$.

$L_3$ = L' × G, the groupoid where L' = $L_5(3)$ the loop of order 6 and G = $\{Z_6 \mid a * b = 2a + 4b \pmod 6$ for all a, b $\in Z_6\}$ is a groupoid and L × {0} is a subloop of the groupoid $L_3$.

Now take P = $\{P_1 \cup P_2 \cup P_3, *_1, *_2, *_3\}$ where

$P_1$ = {e, 2, 5, 8, 14},
$P_2$ = {e, 1, 6, 11, 16, 21, 26, 31, 36, 41} and
$P_3$ = L × {0}, P is a loop.

Clearly P is a 3-loop. Hence L is a Smarandache loop N-loop groupoid.

We call the Smarandache loop N-loop groupoid to be a Smarandache Bruck loop N-loop groupoid if P ⊂ L is a Smarandache Bruck N-loop.

On similar lines we call the N-loop groupoid L to be a Smarandache Bol loop N-loop groupoid if the Smarandache loop N-loop groupoid is a Smarandache Bol N-loop. We call a N-loop groupoid L to be a Smarandache Moufang N-loop groupoid if P ⊂ L where P is a N-loop is a Smarandache Moufang N-loop.

On similar lines we can define Smarandache WIP – loop N-loop groupoid, Smarandache diassociative loop N-loop groupoid and so on.

Next we proceed on to define the notion of Smarandache group N-loop groupoid.

**DEFINITION 5.2.22:** *Let L = $\{L_1 \cup L_2 \cup \ldots \cup L_N, *_1, \ldots, *_N\}$ be a N-loop groupoid. We call L a Smarandache group N-loop groupoid (S-group N-loop groupoid) if L has a proper subset P* *where $P = \{P_1 \cup P_2 \cup \ldots \cup P_N, *_1, \ldots, *_N\}$ is a N-group i.e., each $P_i$ is a group and $P_i \subset L_i$; $P_i = P \cap L_i$.*

***Example 5.2.16:*** Let $L = \{L_1 \cup L_2 \cup L_3 \cup L_4, *_1, *_2, *_3, *_4\}$ where

| | | |
|---|---|---|
| $L_1$ | = | $L_5(3)$ a loop of order 6, |
| $L_2$ | = | $L_7$ (4) a loop of order 8, |
| $L_3$ | = | $\{a, b \in Z_4$ such that $a * b = a + 3b \pmod 4\}$ the groupoid of order 4 and |
| $L_4$ | = | $G = \langle g \mid g^3 = 1\rangle, \times \{Z_8 \mid$ for $a, b \in Z_8$, $a * b = 2a + 6b \pmod 8\}$ $L_4$ is a groupoid of order 24. |

Take $P = \{P_1 \cup P_2 \cup P_3 \cup P_4; *_1, *_2, *_3, *_4\}$ where $P_1 = \{e, g_1\}$, $P_2 = \{e, g_4\}$, $P_3 = \{0, 2\}$ and $P_4 = G\ \{0\}$.

Clearly P is a proper subset of L which has a 4-group structure. So L is a Smarandache group N-loop groupoid.

Now we give some more properties of Smarandache group N-loop groupoids.

**DEFINITION 5.2.23:** *Let $L = \{L_1 \cup L_2 \cup \ldots \cup L_N, *_1, \ldots, *_N\}$ be a N-loop groupoid of finite order. If L is S-loop N-loop groupoid and if every proper subset $P = \{P_1 \cup P_2 \cup \ldots \cup P_N, *_1, \ldots, *_N\}$ which is a N-loop is such that o(P) / o(L) then we call L a Smarandache Lagrange loop N-loop groupoid (S-Lagrange loop N-loop groupoid). If the finite N-loop groupoid has atleast one proper subset $P = \{P_1 \cup P_2 \cup \ldots \cup P_N, *_1, \ldots, *_N\}$ which is a N-loop and is such that o(P) / o(L) then we call L a Smarandache weakly Lagrange loop N-loop groupoid (S-weakly Lagrange loop N-loop groupoid). If the S-loop N-loop groupoid L is such that the order of no proper N-loop P in L divides the order of L then we call L a Smarandache free Lagrange loop N-loop groupoid (S-Smarandache free Lagrange loop N-loop groupoid).*

*The same definition can be given in case of Smarandache group N-loop groupoid. L is such that the order of no proper N-group P in L divides the order of L then we call L a Smarandache free Lagrange group N-loop groupoid. The same*

*definition can be given in case of Smarandache group N-loop groupoid.*

Now we proceed on to define the notion of Smarandache N-loop groupoid.

**DEFINITION 5.2.24:** *Let $L = \{L_1 \cup L_2 \cup \ldots \cup L_N, *_1, \ldots, *_N\}$ be a N-loop groupoid. We call L a Smarandache N-loop groupoid (S-N-loop groupoid) if each $L_i$ is either a S-loop or a S-groupoid. i.e. L has a proper subset set $P = \{P_1 \cup P_2 \cup \ldots \cup P_N, *_1, \ldots, *_N\}$ such that each $P_i = P \cap L_i$ $(i = 1, 2, \ldots, N)$ is either a group or a semigroup.*

***Example 5.2.17:*** Let $L = \{L_1 \cup L_2 \cup L_3 \cup L_4, *_1, *_2, *_3, *_4\}$ be a 4-loop groupoid. Here

$L_1$ = $L_5(3)$ a loop of order 6,
$L_2$ = $L_7(4)$ a loop of order 8 and
$L_3$ = $S(3) \times \{Z_{10}$ | for a, b in $Z_4$, $a * b = a + 3b \pmod 4\}$.

L is a Smarandache N-loop groupoid for, take $P = \{P_1 \cup P_2 \cup P_3 \cup P_4\ *_1, *_2, *_3, *_4\}$ where

$P_1$ = $\{e, g_1\}$,
$P_2$ = $\{e, g_6\}$,
$P_3$ = $S(3)$ and
$P_4$ = $\{0, 2\}$.

Clearly P is a non empty set such that $L_1$ and $L_2$ are S-loops and $L_3$ and $L_4$ are S-groupoids. Hence L is a Smarandache N-loop groupoid.

We now define sub Smarandache N-loop groupoid of a N-loop groupoid L.

**DEFINITION 5.2.25:** *Let $L = \{L_1 \cup L_2 \cup \ldots \cup L_N, *_1, \ldots, *_N\}$ be a N-loop groupoid. Let $P = \{P_1 \cup P_2 \cup \ldots \cup P_N, *_1, \ldots, *_N\}$ be a proper subset of L and P be a sub N-loop groupoid, if P has a proper subset $T = \{T_1 \cup T_2 \cup \ldots \cup T_N, *_1, \ldots, *_N\}$ such that* *each $T_i$ is a group or a semigroup then P is called a Smarandache sub N-loop groupoid (S-sub N-loop groupoid) i.e. if P is a sub N-loop groupoid then P must be S-N loop groupoid for P to be a S-sub N-loop groupoid.*

We have the following result.

**THEOREM 5.2.1:** *Let $L = \{L_1 \cup L_2 \cup \ldots \cup L_N, *_1, \ldots, *_N\}$ be a N-loop groupoid. If L has a S-sub N-loop groupoid then L is a S-N-loop groupoid.*

*Proof:* Suppose L is a N-loop groupoid such that L has a S-sub N-loop groupoid we see L contains a proper subset $P = \{P_1 \cup P_2 \cup \ldots \cup P_N, *_1, \ldots, *_N\}$ which is a Smarandache sub N-loop groupoid so P contains a proper subset $X = \{X_1 \cup X_2 \cup \ldots \cup X_N, *_1, \ldots, *_N\}$ such that $X_i = X \cap P_i$ and $X_i$ is a group or a semigroup for each i = 1, 2, …, N. Now $X \subset P \subset L$ so $X \subset L$ i.e. X is a proper subset of L such that the presence of X makes L a S-N-loop groupoid. Hence the result.

It is still important to note that all S-sub N-loop groupoid is a sub N-loop groupoid but in general every sub N-loop groupoid need not be a S-sub N-loop groupoid.

Now we proceed on to define Smarandache groupoid N-loop groupoid.

**DEFINITION 5.2.26:** *Let $L = \{L_1 \cup L_2 \cup \ldots \cup L_N, *_1, \ldots, *_N\}$ be a N-loop groupoid. We call L a Smarandache groupoid N-loop groupoid (S-groupoid N-loop groupoid) if L contains a proper subset $Y = \{Y_1 \cup Y_2 \cup \ldots \cup Y_N, / *_1, \ldots, *_N\}$ such that each $Y_i$ is a S-N groupoid.*

We cannot always say given any N-loop groupoid it can be S-loop N-loop groupoid or S-group N-loop groupoid or S-groupoid N-loop groupoid. Thus given a N-loop groupoid it can be anything or nothing.

We can easily construct examples of S-loop N-loop groupoid or so on. Thus we are always guaranteed of the existence of such structures.

Now we define Smarandache commutative N-loop groupoid.

**DEFINITION 5.2.27:** *Let* $L = \{L_1 \cup L_2 \cup \ldots \cup L_N, *_1, \ldots, *_N\}$ *be a N-loop groupoid. We call L a Smarandache commutative N-loop groupoid (S-commutative N-loop groupoid) if each* $(L_i, *_i)$ *is either a S-commutative loop if* $(L_i, *_i)$ *is a loop or* $(L_j, *_j)$ *is a S-commutative groupoid if* $(L_i, *_i)$ *is a groupoid.*

*It is interesting to note that a S-commutative N-loop groupoid need not in general be commutative. But every commutative N-loop groupoid is always a S-commutative N-loop groupoid.*

Now we proceed on to define Smarandache homomorphism of N-loop groupoids.

**DEFINITION 5.2.28:** *Let* $L = \{L_1 \cup L_2 \cup \ldots \cup L_N, *_1, \ldots, *_N\}$ *and* $K = \{K_1 \cup K_2 \cup \ldots \cup K_N, *_1, \ldots, *N\}$ *be S-N-loop groupoids. A map* $\phi = \{\phi_1 \cup \ldots \cup \phi_N\}$ *from L to K is called the Smarandache homomorphism of N-loop groupoids (S-homomorphism of N-loop groupoids) if each* $\phi_i$ *is a group homomorphism of a semigroup homomorphism according as* $L_i$ *is a S-loop or a S-groupoid respectively.*

It is important to note that the map may not be defined on the whole of domain space only on the group or semigroup contained in each loop or groupoid $L_i$ respectively. Now we can call the N-loop groupoid to be Smarandache Lagrange semigroup in the following definitions.

**DEFINITION 5.2.29:** *Let* $L = \{L_1 \cup L_2 \cup \ldots \cup L_N, *_1, \ldots, *_N\}$ *be a N-loop groupoid suppose* $L_{i_1}, \ldots, L_{i_K}$ *be K groupoid in the set* $\{L_1, L_2, \ldots, L_N\}$. *We call L a Smarandache groupoid N-loop groupoid (S-groupoid N-loop groupoid) only if each groupoid* $\{L_{i_t} \mid 1 \leq t \leq K\}$ *is a S-groupoid.*

We do not demand the loops in the N-loop groupoid to be S-loops. Now we proceed on to define Smarandache Lagrange groupoid N-loop groupoid, in case the order of the groupoids in the N-loop groupoid is finite. We don't demand the loops to be of finite order.

**DEFINITION 5.2.30:** *Let L = $\{L_1 \cup L_2 \cup \ldots \cup L_N, *_1, \ldots, *_N\}$ be a N-loop groupoid; which has G = $\{L_{i_1}, ..., L_{i_K}\}$ to be K-groupoids. Let the order of each groupoid be of finite order. Suppose L is a S-groupoid N-loop groupoid.*

*Let P = $\left\{P_{i_1}, ..., P_{i_K}\right\}$ be the K-semigroup, if o(P) / o(G) for every K-semigroup in G then we call L a Smarandache Lagrange groupoid N-loop groupoid (S-Lagrange groupoid N-loop groupoid) .*

*If we have atleast one K-semigroup P such that o(P) / o(G) then we call L a Smarandache weakly Lagrange groupoid N-loop groupoid (S-weakly Lagrange groupoid N-loop groupoid).*

*We call L a Smarandache Lagrange free groupoid N-loop groupoid (S-Lagrange free groupoid N-loop groupoid) if G has no K-semigroup P such that o(P) / o(G).*

Now we illustrate this by the following example:

***Example 5.2.18:*** Let L = $\{L_1 \cup L_2 \cup L_3 \cup L_4, *_1, *_2, *_3, *_4\}$ be a S-4-groupoid N-loop groupoid where

$L_1$ and $L_2$ are loops and

$L_3$ = $\{Z_{12}$ | for a, b $\in Z_{12}$ a * b = 2a + 4b (mod 12)$\}$ and

$L_4$ = $\{Z_{18}$ | for a, b $\in Z_{18}$ a * b = 3a + 6 b (mod 18)$\}$ be S-groupoids of finite order.

Clearly G = $\{L_3 \cup L_4\}$ is a 2-groupoid. o(G) = 12 + 18 = 30.

Let P = $P_1 \cup P_2$ = {0, 6} ∪ {0, 6, 12}, P is a 2-semigroup, and o(P) / o(G) i.e. 5/30.

Take T = $T_1 \cup T_2$ = {0, 6} ∪ {0, 6, 9, 12} is a 2-semigroup o(T) / o(G) i.e. 6 / 30.

Take $S = S_1 \cup S_2 = \{0, 6\} \cup \{0, 9\}$ a 2-semigroup, $o(S) \not| \; o(G)$ i.e. $4 \not| \; 30$.

Thus we see L is not a S-Lagrange groupoid N-loop groupoid. In fact L is only a S-weakly Lagrange groupoid N-loop groupoid. Also L is not S-Lagrange free groupoid N-loop groupoid.

Now we proceed on to define Smarandache loop N-loop groupoid.

**DEFINITION 5.2.31:** *Let $L = \{L_1 \cup L_2 \cup \ldots \cup L_N, *_1, \ldots, *_N\}$ be a N-loop groupoid. Let $R = \left\{L_{j_i} \cup \ldots \cup L_{j_s}\right\}$ where each $L_{j_i}$ is a loop from the set $\{L_1, \ldots, L_N\}$. ($i = 1, 2, \ldots, s$) (i.e. all the loops from the set $\{L_1, \ldots, L_N\}$ is in R). Suppose each $L_{j_i}$ is a S-loop then we call L a Smarandache loop N-loop groupoid (S-loop N-loop groupoid).*

We do not demand the groupoids in the set to be S-groupoids. We have the following interesting theorem.

**THEOREM 5.2.2:** *Let $L = \{L_1 \cup L_2 \cup \ldots \cup L_N, *_1, \ldots, *_N\}$ be S-N-loop groupoid, then L is S-loop N-loop groupoid and S-groupoid N-loop groupoid.*

*Proof:* If L is a S-N-loop groupoid then each $L_i$ is either a S-loop or a S-groupoid, so the collection of all loops in $\{L_1, \ldots, L_N\}$ are S-loops and the collection of all groupoids in $\{L_1, \ldots, L_N\}$ are S-groupoids. Hence L a S-loop N-loop groupoid and S-groupoid N-loop groupoid.

Now we proceed on to define Smarandache Lagrange loop N-loop groupoid.

**DEFINITION 5.2.32:** *Let $L = \{L_1 \cup L_2 \cup \ldots \cup L_N, *_1, \ldots, *_N\}$ be a S-loop N-loop groupoid, suppose $R = \left\{L_{j_1}, \ldots, L_{j_t}\right\}$ be the collection of all loops from the set $\{L_1, \ldots, L_N\}$. If each $L_{j_p}$, $p = 1, 2, \ldots, t$ is of finite order, then $o(R) < \infty$. If V =* $\left\{V_{j_1} \cup V_{j_2} \cup ... \cup V_{j_t}, *_{j_i}, ..., *_{j_t}\right\}$ *be a t-group of R and if o(V) / o(R); then we call L a Smarandache weakly Lagrange loop N-loop groupoid (S-weakly Lagrange loop N-loop groupoid).*

*Suppose for each t-group V in R o(V) / o(R) then we call L a Smarandache Lagrange loop N-loop groupoid (S-Lagrange loop N-loop groupoid).*

*We call L a Smarandache Lagrange free loop N-loop groupoid (S-Lagrange free loop N-loop groupoid) if R has no t-group V such that o(V) / o(R).*

It is important to note that L need not be a finite N-loop groupoid in these cases.

**DEFINITION 5.2.33:** *Let L = {$L_1 \cup L_2 \cup \ldots \cup L_N$, $*_1, \ldots, *_N$} be a S-groupoid N-loop groupoid. A proper subset of the collection of all groupoids G =* $\left\{L_{i_1} \cup ... \cup L_{i_K}\right\}$ *of {$L_1 \cup \ldots \cup L_N$} is said to be Smarandache subgroupoid N-loop groupoid (S-subgroupoid N-loop groupoid) if G has a proper subset T =* $\left\{T_{i_1} \cup T_{i_2} \cup ... \cup T_{i_K}\right\}$ *such that T is itself a Smarandache sub K-groupoid.*

Now we proceed on to define ideal of S-groupoid N-loop groupoid.

**DEFINITION 5.2.34:** *Let L = {$L_1 \cup L_2 \cup \ldots \cup L_N$, $*_1, \ldots, *_N$} be a S-groupoid N-loop groupoid. Let G =* $\left\{L_{i_1} \cup ... \cup L_{i_K}\right\}$ *be the collection of groupoids in L. A non empty proper subset P =* $\left\{P_{i_1} \cup ... \cup P_{i_K}\right\}$ *of the K-groupoid G is said to be a K-left ideal of the K-groupoid G if*

i. *P is a sub K-groupoid.*
ii. *For all $x_i \in L_{i_t}$ and $a_i \in P_{i_t}$, $x_i a_i \in P_{i_t}$, t = 1, 2, …, K.*

*One can similarly define K-right ideal of the K-groupoid. We say P is an K-ideal of the K-groupoid G if P is simultaneously a K-left and a K-right ideal of G.*

**DEFINITION 5.2.35:** *Let $L = \{L_1 \cup L_2 \cup \ldots \cup L_N, *_1, \ldots, *_N\}$ be a N-loop groupoid. Let $G = G_{i_1} \cup ... \cup G_{i_K}$ be the K-groupoid; i.e. $G_{i_1}, ..., G_{i_K}$ be the collection of all groupoids in the set $\{L_1, \ldots, L_N\}$. A sub K-groupoid $V = \{V_{i_1} \cup V_{i_2} \cup ... \cup V_{i_K}\}$ of G is said to be normal K subgroupoid of G if*

$$a_i V_{i_t} = V_{i_t} a_i$$

$$(V_{i_t} x_i) y_i = V_{i_t} (x_i y_i)$$

$$y_i (x_i V_{i_t}) = (y_i x_i) V_{i_1}$$

*for all $a_i, x_i, y_i \in G_{i_t}$. The S-groupoid N-loop groupoid is K-simple if it has no nontrivial normal K-subgroupoids.*

We now state when is the S-groupoid N-loop groupoid normal.

**DEFINITION 5.2.36:** *Let $L = \{L_1 \cup L_2 \cup \ldots \cup L_N, *_1, \ldots, *_N\}$ be a S-groupoid N-loop groupoid. Let $G = G_{i_1} \cup ... \cup G_{i_K}$ be the K-groupoid (i.e. $G_{i_1}, ....., G_{i_K}$ be the collection of all groupoids from the set $L_1, \ldots, L_N$). We call L a Smarandache normal groupoid N-loop groupoid (S-normal groupoid N-loop groupoid) if*

i. *xG = Gx where*

   $X = x_{i_1} \cup ... \cup x_{i_K}$ *and*

   $G = G_{i_1} \cup G_{i_2} \cup ... \cup G_{i_K}$ *i.e.*

   $$xG = \left(x_{i_1} G_{i_1} \cup ... \cup x_{i_K} G_{i_K}\right) = \left(G_{i_1} x_{i_1} \cup ... \cup G_{i_K} x_{i_K}\right) = Gx.$$

ii. *G(xy) = (Gx) y; x, y ∈ G where*

   $x = x_{i_1} \cup ... \cup x_{i_K}$ *and*

   $y = y_{i_1} \cup ... \cup y_{i_K}$.

$$G(xy) = \left(G_{i_1} \cup \ldots \cup G_{i_K}\right)\left(x_{i_1} \cup \ldots \cup x_{i_K}\right)\left(y_{i_1} \cup \ldots \cup y_{i_K}\right)$$
$$= \left(G_{i_1}(x_{i_1} y_{i_1}) \cup \ldots \cup G_{i_K}(x_{i_K} y_{i_K})\right)$$
$$= \left\{\left(G_{i_1} x_{i_1}\right) y_{i_1} \cup \ldots \cup \left(G_{i_K} x_{i_K}\right) y_{i_K}\right\}$$
$$= (Gx)\, y.$$

*iii.* $y\,(xG) = (yx)\,G$ where

$$y = y_{i_l} \cup \ldots \cup y_{i_K} \text{ and}$$
$$x = x_{i_l} \cup \ldots \cup x_{i_K}$$
$$y(xG) = \left(y_{i_1} \cup \ldots \cup y_{i_K}\right)\left[\left(x_{i_l} \cup \ldots \cup x_{i_K}\right)\left(G_{i_1} \cup \ldots \cup G_{i_K}\right)\right]$$
$$= y_{i_1}\left(x_{i_1} G_{i_1}\right) \cup y_{i_2}\left(x_{i_2} G_{i_2}\right) \cup \ldots \cup y_{i_K}(x_{i_K} G_{i_K})$$
$$= \left(y_{i_1} x_{i_1}\right) G_{i_1} \cup \left(y_{i_2} x_{i_2}\right) G_{i_2} \cup \ldots \cup (y_{i_K} x_{i_K}) G_{i_K}$$
$$= (yx)G$$

for all $x, y \in G$.

**DEFINITION 5.2.37:** *Let $L = \{L_1 \cup L_2 \cup \ldots \cup L_N, *_1, \ldots, *_N\}$ be a S-groupoid N-loop groupoid. Suppose H and P be two S-K-subgroupoids of $G = \{ L_{i_1} \cup L_{i_2} \cup \ldots \cup L_{i_K} \}$ where $\left\{L_{i_1}, \ldots, L_{i_K}\right\}$ is the collection of all groupoids from the set $\{L_1, L_2, \ldots, L_N\}$ we say H and P are Smarandache conjugate if there exists*

$$x = \left(x_{i_1} \cup \ldots \cup x_{i_K}\right) \in H$$

*such that H = xP or Px where*

$$H = \left\{H_{i_1} \cup \ldots \cup H_{i_K}\right\} \text{ and}$$
$$P = \left\{P_{i_1} \cup \ldots \cup P_{i_K}\right\} \text{ and}$$
$$H = xP$$
$$\text{so } H = \left\{H_{i_1} \cup \ldots \cup H_{i_K}\right\} = \left(x_{i_1} \cup \ldots \cup x_{i_K}\right)\left(P_{i_1} \cup \ldots \cup P_{i_K}\right).$$
$$= \left\{x_{i_1} P_{i_1} \cup \ldots \cup P_{i_K} x_{i_K}\right\} = \left\{P_{i_1} x_{i_1} \cup \ldots \cup P_{i_K} x_{i_K}\right\}$$

*Clearly $H \cap P = \phi$.*

Now we proceed on to define the notion of center of the Smarandache groupoid N-loop groupoid.

**DEFINITION 5.2.38:** *Let L = {$L_1 \cup L_2 \cup \ldots \cup L_N$, $*_1, \ldots, *_N$} be a S-groupoid N-loop groupoid. Let G = $\left\{L_{i_1}, ..., L_{i_K}\right\}$ be the K-groupoid (i.e. $\left\{L_{i_1}, ..., L_{i_K}\right\}$ is the collection of all groupoids from the set {$L_1, \ldots, L_N$} we define Smarandache centre of the S-groupoid N-loop groupoid (S-centre of the S-groupoid N-loop groupoid) L to be $SC_S$ (G) = {x $\in$ G | x a = ax for all a $\in$ G} i.e.*

$$\left\{\left(x_{i_1} \cup ... \cup x_{i_K}\right) \in L_{i_1} \cup ... \cup L_{i_K} \mid ax = \left(a_{i_1} \cup ... \cup a_{i_K}\right)\left(x_{i_1} \cup ... \cup x_{i_K}\right)\right.$$
$$= \left(a_{i_1} x_{i_1} \cup ... \cup a_{i_K} x_{i_K}\right) = xa$$
$$= \left(x_{i_1} \cup ... \cup x_{i_K}\right)\left(a_{i_1} \cup ... \cup a_{i_K}\right)$$
$$= \left.\left(x_{i_1} a_{i_1} \cup ... \cup x_{i_K} a_{i_K}\right)\right\}.$$

*We also call a pair of elements a, b in the K groupoid G = $\left\{L_{i_1} \cup ... \cup L_{i_K}\right\}$, Smarandache conjugate pair (S-conjugate pair) if a = b x (x. b for some x $\in$ G) and b = a. y (y a for some y $\in$ G).*

*A = ($a_{i_1} \cup ... \cup a_{i_K}$) is right Smarandache conjugate (right S-conjugate) with b in G if we can find x, y $\in$ G such that ax =b and by = a (xa = b and yb = a),*

*i.e.* $(b_{i_1} \cup ... \cup b_{i_K})(y_{i_1} \cup ... \cup y_{i_K})$

$$= (b_{i_1} y_{i_1} \cup ... \cup b_{i_K} y_{i_K})$$
$$= (a_{i_1} \cup ... \cup a_{i_K}).$$

$$(x\,a = (x_{i_1} \cup ... \cup x_{i_K})(a_{i_1} \cup ... \cup a_{i_K})$$
$$= \left(x_{i_1} a_{i_1} \cup ... \cup x_{i_K} a_{i_K}\right) = \left(b_{i_1} \cup ... \cup b_{i_K}\right) = b\,)$$

*and*

$$(y\,b = (y_{i_1} \cup ... \cup y_{i_K})(b_{i_1} \cup ... \cup b_{i_K})$$
$$= (y_{i_1} b_{i_1} \cup ... \cup y_{i_K} b_{i_K}) = (a_{i_1} \cup ... \cup a_{i_K}) = a\,)$$

**DEFINITION 5.2.39:** *Let $L = \{L_1 \cup L_2 \cup \ldots \cup L_N, *_1, \ldots, *_N\}$ be a S-groupoid N-loop groupoid. $G = \left\{L_{i_1} \cup \ldots \cup L_{i_K}\right\}$ be a S-K-groupoid (i.e. $L_{i_1} \cup \ldots \cup L_{i_K}$ be the collection of all of groupoids of the set $(L_1, \ldots, L_N)$) .*

*We call a S-sub K groupoid $H = \{H_{i_1} \cup \ldots \cup H_{i_K}\}$ where $H_{i_t} = H \cap L_{i_t}$ $(1 \le t \le K)$ to be Smarandache semi normal K-groupoid (S-semi normal K-groupoid) of N-loop groupoid L if*

*i. $aH = X$ for all $a \in G$ i.e.*

$$\left(a_{i_1} \cup a_{i_2} \cup \ldots \cup a_{i_K}\right)\left(H_{i_1} \cup \ldots \cup H_{i_K}\right)$$
$$= \left(a_{i_1} H_{i_1} \cup \ldots \cup a_{i_K} H_{i_K}\right) = X$$
$$= \left(X_{i_1} \cup \ldots \cup X_{i_K}\right).$$

*ii. $Ha = Y$ for all $a \in G$ i.e.*

$$\left(H_{i_1} \cup \ldots \cup H_{i_K}\right)\left(a_{i_1} \cup \ldots \cup a_{i_K}\right)$$
$$= \left(H_{i_1} a_{i_1} \cup \ldots \cup H_{i_K} a_{i_K}\right) = Y$$
$$= \left(y_{i_1} \cup \ldots \cup y_{i_K}\right),$$

*where both X and Y are sub K- groupoids.*

Interested reader may find some interesting relations between the S-semi normal K- groupoid and S- normal K-groupoid of a N-loop groupoid.

One can also define Smarandache semi conjugate sub K-groupoids of a N-loop groupoid.

**DEFINITION 5.2.40:** *Let $L = \{L_1 \cup L_2 \cup \ldots \cup L_N, *_1, \ldots, *_N\}$ be a S-groupoid N-loop groupoids where $G = \left\{L_{i_1} \cup \ldots \cup L_{i_K}\right\}$ is the K-groupoid (Here $L_{i_1} \cup \ldots \cup L_{i_K}$ is the collection of all groupoids from $L_1, \ldots, L_N$). Let H and P be two sub K-groupoids of K. We say H and K are Smarandache K-semi conjugate sub K-groupoids (S-K-semi conjugate sub K-groupoids) of G if*

*i. H and P are S- sub K-groupoids of G.*
*ii. H = xP or Px or*

$$H = \left(H_{i_1} \cup \ldots \cup H_{i_K}\right)$$
$$= \left(x_{i_1} \cup \ldots \cup x_{i_K}\right) \left(P_{i_1} \cup P_{i_2} \cup \ldots \cup P_{i_K}\right)$$
$$= \left(x_{i_1} P_{i_1} \cup \ldots \cup x_{i_K} P_{i_K}\right)$$
$$[\,or \left(P_{i_1} \cup \ldots \cup P_{i_K}\right) \left(x_{i_1} \cup \ldots \cup x_{i_K}\right)$$
$$= \left(P_{i_1} x_{i_1} \cup \ldots \cup P_{i_K} x_{i_K}\right)$$
$$= Px]\ or$$

*P = x H i.e.*

$$\left(P_{i_1} \cup \ldots \cup P_{i_K}\right)$$
$$= \left(x_{i_1} \cup \ldots \cup x_{i_K}\right) \left(H_{i_1} \cup \ldots \cup H_{i_K}\right)$$
$$= \left(x_{i_1} H_{i_1} \cup \ldots \cup x_{i_K} H_{i_K}\right)$$
$$[or\ Hx = \left(H_{i_1} \cup \ldots \cup H_{i_K}\right) \bullet \left(x_{i_1} \cup \ldots \cup x_{i_K}\right)$$
$$= \left(H_{i_1} x_{i_1} \cup \ldots \cup H_{i_k} x_{i_k}\right)]$$

*for some* $\left(x_{i_1} \cup \ldots \cup x_{i_k}\right) \in G$ *i.e. each* $x_{i_t} \in L_{i_t}$ *(*$1 \le t \le K$*).*

On similar lines one can define for S loop N-loop groupoids i.e. by using t-loops instead of K-groupoids. Several interesting results can be obtained, which is left as an exercise for the reader.

## 5.3 N-group loop semigroup groupoid (glsg) algebraic structures

In this section we proceed on to define the mixed N-algebraic structures, which include both associative and non associative structures. Here we define them and give their substructures and a few of their properties.

**DEFINITION 5.3.1:** *Let A be a non empty set on which is defined N-binary closed operations $*_1, \ldots, *_N$. A is called as the N-group-loop-semigroup-groupoid (N-glsg) if the following conditions, hold good.*

i. *$A = A_1 \cup A_2 \cup \ldots \cup A_N$ where each $A_i$ is a proper subset of A (i.e. $A_i \not\subseteq A_j \not\subseteq$ or $A_j \not\subseteq A_i$ if $(i \neq j)$.*
ii. *$(A_i, *_i)$ is a group or a loop or a groupoid or a semigroup (or used not in the mutually exclusive sense) $1 \leq i \leq N$. A is a N –glsg only if the collection $\{A_1, \ldots, A_N\}$ contains groups, loops, semigroups and groupoids.*

*Note:* If each $(A_i, *_i)$ is just only group or (and) semigroup we get the already defined N-group semigroup. If each $(A_j, *_j)$ is just a loop and (or) groupoid we obtain the already defined N-loop groupoid structure. If $\{A_i, *_i\}$ are only group, loop and semigroup then we call A the N-group loop semigroup i.e. N-gls structure.

If each $\{A_i, *_i\}$ are only just a group, semigroup and groupoid then we call A just the N-group semigroup groupoid i.e. N-gsg structure.

Likewise we can have a N-lsg structure or just N-gls structure also. Whenever situation arises we will make a mention of them.

We just give some examples in order to make the definition more concrete.

***Example 5.3.1:*** Let $A = \{A_1 \cup A_2 \cup A_3 \cup A_4; *_1, *_2, *_3, *_4\}$ where

$A_1$ = $G = \langle g \mid g^8 = 1\rangle$,
$A_2$ = Loop of order six given in page 135,
$A_3$ = $\{Z_{12}$, the semigroup under multiplication modulo 12$\}$ and
$A_4$ = $\{Z_8$, the groupoid defined by $a * b = 2a + 6b \pmod 8\}$.

Clearly A is a 4-glsg.

Now we can define the order of any N-glsg. The order of N-glsg is the number of distinct elements in A. If the number of elements in A is finite we say A is a finite N-glsg. If the order of A is infinite we say A is an infinite N-glsg.

Now we can also define N-gls or N-gls or N-lsg or Ngsg and so on. Now we define the notion of sub N-group-loop-semigroup groupoid of a N-glsg.

**DEFINITION 5.3.2:** *Let $A = \{A_1 \cup \ldots \cup A_N, *_1, \ldots, *_N\}$ where $A_i$ are groups, loops, semigroups and groupoids. We call a non empty subset $P = \{P_1 \cup P_2 \cup \ldots \cup P_N, *_1, \ldots, *_N\}$ of A, where $P_i = P \cap A_i$ is a group or loop or semigroup or groupoid according as $A_i$ is a group or loop or semigroup or groupoid. Then we call P to be a sub N-glsg.*

It is important to note that even if A is a finite N-glsg and P a sub N-glsg; o(P) need not in general divide order of A (o(A)).

We just illustrate this by an example.

***Example 5.3.2:*** Let $A = \{A_1 \cup A_2 \cup A_3 \cup A_4; *_1, *_2, *_3, *_4\}$ be a finite N-glsg, where

$A_1$ = $\{D_{2,6} = \{a, b \mid a^2 = b^6 = 1, b a b = a\}$ the dihedral group of order 12,
$A_2$ = $\{L_5(3)\}$ be the loop of order 6,
$A_3$ = $\{Z_8 = (0, 1, 2, \ldots, 7)$ the semigroup under multiplication modulo 8 and
$A_4$ = $\{Z_{10} \mid$ for $a, b \in Z_4$, $a *_4 b = a + 4b \pmod{10}\}$.

Now o(A) = 12 + 6 + 8 + 10 = 36, a N-glsg structure of finite order.

Let $P = \{P_1 \cup P_2 \cup P_3 \cup P_4 *_1, *_2, *_3, *_4\}$ be a proper subset of A, where

$P_1$ = $\{1, b, \ldots, b^5\}$,
$P_2$ = $\{e, g_1\}$,
$P_3$ = $\{0, 4\}$ and
$P_4$ = $\{0, 5\}$.

Clearly P is a sub N-glsg. o(P) = 6 + 2 + 2 + 2 = 12, 12 / 36.

Take R = $\{R_1 \cup R_2 \cup R_3 \cup R_4 *_1, *_2, *_3, *_4\}$ be a proper subset of A. Take

$R_1 = \{1, b, ..., b^5\}$,
$R_2 = \{e, g_3\}$,
$R_3 = \{0, 2, 4, 6\}$ and
$R_4 = \{0, 2, 4, 6, 8\}$.

o(R) = 6 + 2 + 4 + 5 = 17. o(R) $\nmid$ o(A). Thus we see the order all sub-N-glsg P of a finite N-glsg A does not divide order of A.

Thus we make the following definition.

**DEFINITION 5.3.3:** *Let A be a N-glsg structure where A = $\{A_1 \cup A_2 \cup \ldots \cup A_N; *_1, \ldots, *_N\}$ be of finite order. If every sub N-glsg P of A is such that o(P) / o(A) then we call A a Lagrange N-glsg.*

*It at least one sub N-glsg R exists in A such that o(R) / o(A) then we call A the weakly Lagrange N-glsg. If A has no sub N-glsg T such that o(T) / o(A) then we call A a Lagrange free N-glsg.*

We have clearly several Lagrange free N-glsg. For take the order of the N-glsg. For take the order of the N-glsg, A = $\{A_1 \cup A_2 \cup \ldots \cup A_N; *_1, \ldots, *_N\}$ such that o(A) = p, p a prime. Then 1 or p alone divides p so no sub N-glsg of A divides o(A). It is interesting to note that every Lagrange N-glsg is a weakly Lagrange N-glsg.

***Example 5.3.3:*** Let A = $\{A_1 \cup A_2 \cup A_3 \cup A_4; *_1, *_2, *_3, *_4\}$ where

$A_1$ = $S_3$, the symmetric group of degree 3,
$A_2$ = $L_7(3)$, the loop of order 8,
$A_3$ = $\{Z_6$, the semigroup under multiplication modulo 6$\}$ and
$A_4$ = $\{Z_9 \mid \text{for } a, b \in Z_9, a * b = 3a + 6b \pmod{19}\}$.

A is a finite N-glsg. Now o(A) = 6 + 8 + 6 + 9 = 29 a prime.

Let P = {$P_1 \cup P_2 \cup P_3 \cup P_4$, $*_1$, $*_2$, $*_3$, $*_4$} where

$$P_1 = \left\{\begin{pmatrix} 1 & 2 & 3 \\ 1 & 2 & 3 \end{pmatrix}, \begin{pmatrix} 1 & 2 & 3 \\ 2 & 3 & 1 \end{pmatrix}, \begin{pmatrix} 1 & 2 & 3 \\ 3 & 1 & 2 \end{pmatrix}\right\} = A_3,$$

$P_2$ = {e, $g_2$} $\subset L_7(3)$,
$P_3$ = {0, 2, 4} $\subset Z_6$ and
$P_4$ = {0, 6, 3} $\subset Z_9$.

o(P) = 3 + 2 + 3 + 3 = 11. Thus A is a Lagrange free N-glsg. It is to be noted when we say A is Lagrange free a N-glsg we do not mean A cannot have proper subsets P such that P is a sub N-glsg.

***Example 5.3.4:*** Let A = {$A_1 \cup A_2 \cup A_3 \cup A_4$; $*_1$, $*_2$, $*_3$, $*_4$} where

$A_1$ = {g | $g^4 = 1$},
$A_2$ = $L_5(3)$, the loop of order 6,
$A_3$ = {$Z_6$, multiplication modulo 6} is a semigroup and
$A_4$ = {$Z_3$ / a, b $\in Z_5$, a $*_3$ b = a + 4b (mod 5)} is a groupoid.

A is a 4-glsg, o(A) = 4 + 6 + 6 + 5 = 21. Clearly A cannot have sub4-glsg of order 3.

Take P = {$P_1 \cup P_2 \cup P_3 \cup P_4$, $*_1$, $*_2$, $*_3$, $*_4$} where

$P_1$ = {1, $g^2$},
$P_2$ = {$p_1$, $g_1$},
$P_3$ = {0, 2, 4} $\subset Z_6$ and
$P_4$ = {$Z_5$}, P is a 4 sub-glsg.

o(P) = 2 + 2 + 3 + 5 = 12. 12 $\nmid$ 21. Since $P_4$ cannot have any proper subgroupoid other then {0} and itself we see any sub4-glsg will have order strictly greater than 8 so A is a Lagrange free 4-glsg where the order of the 4-glsg is a composite number.

Now we proceed on to define the special notion of sub K-group of N-glsg, sub m-loop of N-glsg, sub t-semigroup of N-glsg and sub q-groupoids of N-glsg.

**DEFINITION 5.3.4**: *Let A = $\{A_1 \cup A_2 \cup \ldots \cup A_N; *_1, \ldots, *_N\}$ be a N-glsg. A proper subset T = $\left\{T_{i_1} \cup \ldots \cup T_{i_K}, *_{i_1}, \ldots, *_{i_K}\right\}$ of A is called the sub K-group of N-glsg if each $T_{i_t}$ is a group from $A_r$ where $A_r$ can be a group or a loop or a semigroup of a groupoid but has a proper subset which is a group.*

***Example 5.3.5:*** A = $\{A_1 \cup A_2 \cup A_3 \cup A_5, *_1, \ldots, *_5\}$ be a 5-glsg, where

$A_1$ = $\{S_4$, the symmetric group of degree 4$\}$,
$A_2$ = $L_5(3)$ the loop of order 6,
$A_3$ = $L_7(4)$ the loop of order 8,
$A_4$ = S(3), the symmetric semigroup with $3^3$ elements and
$A_5$ = $\{Z_{12}$, for a, b $\in Z_6$, a * b = 4a + 2b mod 6$\}$ is a groupoid.

Take P = $\{P_1 \cup P_2 \cup P_3 \cup P_4, *_1, \ldots, *_5\}$ where

$P_1$ = $A_4$, the alternating subgroup of $S_4$,
$P_2$ = $\{e, g_2\}$,
$P_3$ = $\{e, g_6\}$,
$P_4$ = $\left\{\begin{pmatrix} 1 & 2 & 3 \\ 1 & 2 & 3 \end{pmatrix}, \begin{pmatrix} 1 & 2 & 3 \\ 2 & 1 & 3 \end{pmatrix}\right\} \subset$ S(3) is a sub 4-group of the 5-glsg.

Thus the K cannot be said as the maximum number of groups in $A_1, \ldots, A_N$ we can only say minimum of K will be maximum number of groups in the collection $\{A_1, \ldots, A_N\}$ provided no group $A_i$ is a cyclic group of prime order.

**DEFINITION 5.3.5:** *Let A = $\{A_1 \cup A_2 \cup \ldots \cup A_N; *_1, \ldots, *_N\}$ be a N-glsg. A proper subset T = $\left\{T_{i_1} \cup \ldots \cup T_{i_r}\right\}$ is said to be sub r-* *loop of A if each $T_{i_j}$ is a loop and $T_{i_j}$ is a proper subset of some $A_p$. As in case of sub K-group r need not be the maximum number of loops in the collection $A_1, \ldots, A_N$.*

***Example 5.3.6:*** Let A = $\{A_1 \cup A_2 \cup A_3 \cup A_5; *_1, *_2, *_3, *_4, *_5\}$ be a 5-glsg where

$A_1$ = $\{g \mid g^7 = 1\}$,
$A_2$ = $L_5(3) \times \{D_{2.3}\}$,
$A_3$ = $\{L_7(4) \times S_3\}$,
$A_4$ = $\{Z_{16}$, a semigroup under multiplication modulo 16$\}$and
$A_5$ = $\{Z_{10} \mid$ for a, b $\in Z_{10}$, a * b = 2a + 3b (mod 10)$\} \times L_{11}(3)$.

Clearly $A_1$ is a group. $A_2$ and $A_3$ are loops, $A_4$ is a semigroup and $A_5$ is a groupoid.

Take P = $\{P_1 \cup P_2 \cup P_3\}$ where

$P_1$ = $L_5(3) \subset A_2$ ,
$P_2$ = $L_7(4) \subset A_3$ and
$P_3$ = $L_{11}(3) \subset A_5$

P is a sub3-loop of the 5-glsg. However we see that there are only 2 loops in the set $\{A_1, A_2, A_3, A_4, A_5\}$.

Now we proceed on to define the sub u semigroup of the N-glsg (u < N).

**DEFINITION 5.3.6:** *Let A = $\{A_1 \cup A_2 \cup \ldots \cup A_N; *_1, \ldots, *_N\}$ be a N-glsg. Let P = $\{P_1 \cup P_2 \cup \ldots \cup P_N\}$ be a proper subset of A where each $P_i$ is a semigroup then we call P the sub u-semigroup of the N-glsg.*

Now we illustrate this by the following example:

***Example 5.3.7:*** Let A = $\{A_1 \cup A_2 \cup A_3 \cup A_4 \cup A_5; *_1, \ldots, *_5\}$ be a 5-glsg. Take

$A_1$ = $\{S_3\}$ be the group,
$A_2$ = $\{L_7(3)$, loop of order 8$\}$,
$A_3$ = $\{Z_{18}$, semigroup under multiplication modulo 18$\}$,
$A_4$ = $\{Z_{10}$, semigroup under multiplication modulo 10$\}$ $\times$ $\{Z_6 \mid$ for a, b $\in Z_6$, a * b = 2a + b (mod 6)$\}$ = $B_1 \times B_2$ is a groupoid, and
$A_5$ = $\{Z_7 \mid$ for a, b $\in Z_7$, a * b = 3a + 4b (mod 7)$\}$ $\times$ $\left\{\begin{pmatrix} a & b \\ c & d \end{pmatrix}\right.$ / a b c d $\in$ Q matrix under multiplication$\}$ = $C_1 \times C_2$.

Take S = $\{S_1 \cup S_2 \cup S_3\}$ where

$S_1$ = $\{0, 2, 4, 6, 8\}$ semigroup under multiplication modulo 10,
$S_2$ = $B_1 \times \{0\} \subseteq B_1 \times B_2 = A_4$ a semigroup and
$S_3$ = $\{0\} \times C_2 \subseteq C_1 \times C_2 \subset A_5$ semigroup.

Clearly S is a sub3-semigroup of 5 glsg.

Now we proceed on to define the concept of sub t-groupoid of N-glsg (t < N).

**DEFINITION 5.3.7:** *Let A = $\{A_1 \cup A_2 \cup \ldots \cup A_N; *_1, \ldots, *_N\}$ be a N-glsg. A proper subset C = $\{C_1 \cup C_2 \cup \ldots \cup C_t\}$ of A is said to be a sub-t-groupoid of A if each $C_i$ is a groupoid.*

Now we illustrate this by the following example:

***Example 5.3.8:*** Let A = $\{A_1 \cup A_2 \cup A_3 \cup A_4 \cup A_5; *_1, \ldots, *_5\}$ be a 5-glsg, where

$A_1$ = {Group of quaternions},
$A_2$ = $L_5(3)$, the loop of order 6,
$A_3$ = $\{S(5)$, the symmetric semigroup$\}$;
$A_4$ = $\{Z_{12} \mid$ for a, b $\in Z_{12}$, a $*_4$ b = 2a + 4b (mod 12)$\}$ and
$A_5$ = $\{Z_{15}$ / for a, b $\in Z_{15}$ define a $*_5$ b = a + 4b (mod 15)$\}$ are groupoids.

Now define $D = D_1 \cup D_2$ to be the sub 2- groupoid where

$D_1 = \{0, 2, 4, 6, 8, 10\} \subset Z_{12}$ and
$D_2 = \{0, 3, 6, 9, 12\} \subseteq Z_{15}$.

Now we proceed on to analyze when these special substructures of this finite N-glsg are such that they divide the order of N-glsg, never divide the order of N-glsg and few of them divide the order of N-glsg.

**DEFINITION 5.3.8:** *Let $A = \{A_1 \cup A_2 \cup \ldots \cup A_N; *_1, \ldots, *_N\}$ be a N-glsg. Suppose A contains a subset $P = P_{L_1} \cup \ldots \cup P_{L_k}$ of A such that P is a sub K-group of A. If every P-sub K-group of A is commutative we call A to be a sub-K-group commutative N-glsg.*

*If atleast one of the sub-K-group P is commutative we call A to be a weakly sub K-group-commutative N-glsg. If no sub K-group of A is commutative we call A to be a non commutative sub-K-group of N-glsg.*

*Note:* These definitions can be made by replacing the sub-K-group by sub-t-groupoid, sub-p-loop and semigroup-semigroup of a N-glsg on similar lines. Now we will represent all these concepts by examples.

***Example 5.3.9:*** Let $A = \{A_1 \cup A_2 \cup A_3 \cup A_4 \cup A_5 \cup A_6; *_1, \ldots, *_6\}$ be a 6-glsg, where

$A_1 = \{g \mid g^{12} = 1\}$, cyclic group of order 12,
$A_2 = \{A_4$, the alternating group$\}$,
$A_3 = \{L_3(3)$, the loop of order 6$\}$,
$A_4 = \{Z_{20}$, the semigroup under multiplication modulo 20$\}$,
$A_5 = \{Z_{12}$, semigroup under multiplication modulo 12$\} \times \{Z_9 \mid$ for $a, b \in Z_9$, $a * b = 3a + 6b \pmod 9\} = B_1 \times B_2$, and
$A_6 = \{Z_{15} \mid$ for $a, b \in Z_{15}$, $a * b = 2a + 3b \pmod{15}\} \times \{Z_6$ semigroup under multiplication modulo 6$\} = C_1 \times C_2$.

The subset P = $\{P_1 \cup P_2 \cup P_3 \cup P_5 \cup P_6\}$ of A where

$P_1$ = $\{1, g^2, g^4, g^6, g^8, g^{10}\} \subseteq A_1$,

$$P_2 = \left\{\begin{pmatrix} 1 & 2 & 3 & 4 \\ 4 & 3 & 2 & 1 \end{pmatrix}, \begin{pmatrix} 1 & 2 & 3 & 4 \\ 2 & 1 & 4 & 3 \end{pmatrix}, \begin{pmatrix} 1 & 2 & 3 & 4 \\ 3 & 4 & 1 & 2 \end{pmatrix}, \begin{pmatrix} 1 & 2 & 3 & 4 \\ 4 & 3 & 2 & 1 \end{pmatrix}\right\}$$

$\subseteq A_4$,

$P_3$ = $\{e, g_1\} \subseteq L_5(3)$,

$P_4$ = $\{1,11\} \subset Z_{20}$,

$P_5$ = $\{3, 9\} \subseteq Z_{12}$ {9 acts as the identity of the group for 9.9 = 9 (mod 6) 3.9 = 9.3 = 3 (mod 6) and 3.3 = 9 (mod 6)}

and

$P_6$ = $\{1, 5\} \subseteq Z_6$ {where 5.5 = 1 (mod 6), 5.1 = 1.5 = 5 (mod 6)}.

Clearly P is sub 6-group of A. Clearly P is a 6-commutative 6-glsg. We see A is not a commutative structure for the group $A_4$ is not a commutative group.

It is interesting to note that in this example we could get a sub K-group of the 6 glsg where K = N = 6.

We see this example 5.3.9 has no sub t-loop for both the loops do not have any proper subloops. Now we find the sub r-semigroup of A = $\{A_1 \cup A_2 \cup \ldots \cup A_6\}$.

Let S = $\{S_1 \cup S_2 \cup S_3\}$ be a proper subset of A where

$S_1$ = $\{0, 2, 6, 8, 10, 12, \ldots, 18\} \subset Z_{20}$,

$S_2 \times \{0\}$ = $\{0, 3, 6, 9\} \times \{0\} \subset Z_{12} \times Z_9$ and

$S_3 \times \{0\}$ = $\{0\} \times \{0, 2, 4\} \subseteq C_1 \times C_2$.

Clearly S is a sub 3-semigroup of A which is commutative. Now we see whether this N-glsg A has any sub t-groupoids.

Take T = $\{T_1 \cup T_2\}$ where

$T_1$ = $\{0\} \times \{0, 3, 6\} \subseteq B_1 \times B_3$ and

$T_2$ = $\{0, 3, 6, 9, 12\} \times \{0\} \subseteq C_1 \times C_2$

are subgroupoids of $A_5$ and $A_6$ respectively. So the N-glsg has sub2-groupoids which is not commutative sub2-groupoids of A.

**DEFINITION 5.3.9**: *Let $A = \{A_1 \cup A_2 \cup \ldots \cup A_N; *_1, \ldots, *_N\}$ be a N-glsg of finite order. Suppose A has sub K-groups H and if every sub K-group H is such that o(H) / o(A) then we call A to be a Lagrange sub K-group. If A has atleast one sub K-group T such that o(T) / o(A) then we call A to be a weakly sub K-group. If A has sub K-group such that the order of none of them divide the order of A then we call A to be a Lagrange free sub K-group.*

It is important to note A may not at times have any sub K-group.

Now we illustrate them by examples.

***Example 5.3.10:*** Let $A = \{A_1 \cup A_2 \cup A_3 \cup A_5 \cup A_6; *_1, \ldots, *_6\}$ be a 6-glsg of finite order where

$A_1$ = $\{g \mid g^6 = 1\}$ cyclic group of order 6,
$A_2$ = $S_4$, the symmetric group of degree 4,
$A_3$ = $L_5(2)$, the loop of order 6,
$A_4$ = $L_7(3)$, the loop of order 8,
$A_5$ = $\{Z_{12}$ , semigroup under multiplication modulo 12$\}$ and
$A_6$ = $\{Z_{18} \mid$ for a, b $\in Z_8$, a $*_6$ b = 3a + 5b (mod 8)$\}$ groupoid of order 8.

Now order of A = o (A) = 6 + 24 + 6 + 8 + 12 + 8 = 64.

Take $P = \{P_1 \cup P_2 \cup P_3 \cup P_5\}$ a proper subset of A where

$P_1$ = $\{1, g^3\}$,

$$P_2 = \left\{1, \begin{pmatrix} 1 & 2 & 3 & 4 \\ 2 & 1 & 4 & 3 \end{pmatrix}\right\},$$

$P_3$ = $\{e, g_1\}$,
$P_4$ = $\{e, g_6\}$ and
$P_5$ = $\{1, 11\} \subset Z_{12}$.

Clearly P is a sub- 5- group of A. o(P) = 2 + 2 + 2 + 2 + 2 = 10. We see 10 $\nmid$ 64. Now consider T = $\{T_1 \cup T_2 \cup T_3 \cup T_4 \cup T_5\}$ a proper set subset of A where

$T_1$ = $\{g^2, g^4, 1\}$

$$T_2 = \left\{\begin{pmatrix}1 & 2 & 3 & 4\\1 & 2 & 3 & 4\end{pmatrix}, \begin{pmatrix}1 & 2 & 3 & 4\\2 & 3 & 1 & 4\end{pmatrix}, \begin{pmatrix}1 & 2 & 3 & 4\\3 & 1 & 2 & 4\end{pmatrix}\right\}$$

$T_3$ = $\{1, g_2\}$,
$T_4$ = $\{1, g_3\}$ and
$T_5$ = $\{4, 8\}$ is a sub5-group of A.

o(T) = 3 + 3 + 2 + 2 + 2 = 12, 12 $\nmid$ 64. Thus we see we have sub t-groups but their order does not divide the order of A.
Now we consider yet another example.

***Example 5.3.11:*** Let A = $\{A_1 \cup A_2 \cup A_3 \cup A_4 \cup A_5; *_1, *_2, *_3, *_4, *_5\}$ be a 5-glsg where

$A_1$ = $\{g \mid g^{12} = 1\}$,
$A_2$ = $\{Z_8$, group under addition modulo 8$\}$,
$A_3$ = $L_5(3)$ loop of order 6,
$A_4$ = $\{Z_9$, semigroup under multiplication modulo 9$\}$ and
$A_5$ = $\{Z_{10}$ / for a, b $\in Z_{10}$, a * b = 3a + 2b (mod 10)$\}$.

o(A) = 12 + 8 + 6 + 9 + 10 = 45.

Consider P = $\{P_1 \cup P_2 \cup P_3 \cup P_4\}$ a proper subset of A where

$P_1$ = $\{g^4, g^8, 1\}$, $P_2$ = $\{0, 4\}$,
$P_3$ = $\{e, g_1\}$ and $P_4$ = $\{1, 8\}$.

Now o (P) = 3 + 2 + 2 + 2 = 9. o (P) / o (A) i.e. 9 / 45.

Now consider the subset T = $\{T_1 \cup T_2 \cup T_3 \cup T_4\}$ where

$T_1$ = $\{g^4, g^8\ 1\}$,

$T_2$ = {0, 2, 4, 6},
$T_3$ = {e, $g_1$} and
$T_4$ = {1, 2, 4, 5, 7, 8} is a 4-group.

o (T) = 3 + 4 + 2 + 6 = 15, i.e. o (T) /o (A) i.e. 15 / 45.

Take a proper subset R = {$R_1 \cup R_2 \cup R_3 \cup R_4$} of A where

$R_1$ = {1, $g^2$, $g^4$, $g^6$, $g^8$, $g^{10}$},
$R_2$ = {0, 4},
$R_3$ = {e, $g_3$} and
$R_4$ = {1, 2, 4, 5, 7, 8}.

o(R) = 6 + 2 + 2 + 6 = 16. Clearly o(R) $\nmid$ o (A) i.e. 16 $\nmid$ 45. Just we see A in the above example is a weakly Lagrange sub4-group of A. Further A is not a Lagrange 4-group of A. We now construct one more example.

***Example 5.3.12:*** Let A = {$A_1 \cup A_2 \cup A_3 \cup A_4 \cup A_5$; $*_1$,$*_2$, $*_3$, $*_4$,$*_5$} where

$A_1$ = {g | $g^6$ = 1},
$A_2$ = $L_5$ (3), a loop of order 6,
$A_3$ = {$Z_{10}$, semigroup under multiplication modulo 10},
$A_4$ = {$Z_8$, semigroup under multiplication modulo 8} and
$A_5$ = {$Z_7$ / for a, b $\in$ we have, a * b = 3a + 4b (mod 7)}.

o(A) = 6 + 6 + 10 + 8 + 7 = 37. Clearly A is a prime order N-glsg. We see A has sub-4 groups but A is a Lagrange free N-glsg.

For take W = {$W_1 \cup W_2 \cup W_3 \cup W_4$} where

$W_1$ = {1, $g^2$ $g^4$} $\subset A_1$,
$W_2$ = {e, $g_3$} $\subseteq L_5$ (3),
$W_3$ = {1, 9} $\subset Z_{10}$ and
$W_4$ = {1, 7} $\subseteq Z_8$ .

o (W) = 3 + 2 + 2 + 2 = 9, 9 $\nmid$ 37.

It is easily verified that A has sub K-group but since o(A) is prime none of them will divide o (A). So A is a Lagrange free sub K-group.

Now we define when is the N-glsg a Lagrange sub K-loop and so on.

**DEFINITION 5.3.10:** *Let A = {$A_1 \cup A_2 \cup \ldots \cup A_N$; $*_1, \ldots, *_N$} be a N-glsg of finite order. Suppose for every P = {$P_1 \cup P_2 \cup \ldots \cup P_N$} a sub-t-loop of A we have o(P) / o(A) then we call A to be a Lagrange sub t- loop.*

*On the other hand if A has atleast one sub t loop, P such that o(P) / o(A) then we call A to be a weakly Lagrange loop-loop. If no loop-loop P of A is such that o(P) / o(A) then we call A to be a Lagrange free loop-loop.*

Now we give illustration of these by some suitable examples.

***Example 5.3.13:*** Let A = {$A_1 \cup A_2 \cup A_3 \cup A_4 \cup A_5$; $*_1, \ldots, *_5$} be a N-glsg of finite order where

$A_1$ = $\langle g \mid g^7 = 1\rangle$, cyclic group of order 7,
$A_2$ = $L_5(3) \times L_7(3)$ (loop of order 48),
$A_3$ = {$Z_4$, the semigroup under multiplication modulo 4},
$A_4$ = {$Z_5$ | for a, b $\in Z_5$, a * b = 2a + 3b (mod 5)} and
$A_5$ = $L_5(2) \times$ {$Z_3$, such that for a, b $\in Z_3$, a * b = a + 2b (mod 3)}.

Now o(A) = 7 + 48 + 4 + 5 + 18 = 82. Consider the proper subset M = {$M_1 \cup M_2$} of A where $M_1 = L_7(3)$ and $M_2 = L_5(2)$. Clearly M is a sub 2 loop of A of order 14. Now o(M) $\not|$ o(A).

Thus we see the N-glsg has no loop-loop P such that o(P) / o(A). So A is a free Lagrange t-loop. It is left for the reader to verify that A is also a free Lagrange K-group. Now we proceed on to see other examples.

***Example 5.3.14:*** Let A = {$A_1 \cup A_2 \cup A_3 \cup A_4$, $*_1, \ldots, *_4$} where

$A_1$ = $\{g \mid g^7 = 1\}$,
$A_2$ = $L_5(3) \times \{g \mid g^3 = 1\}$,
$A_3$ = $L_7(2) \times \{Z_2 = \{0, 1\}$ group under '+' modulo 2} and
$A_4$ = $L \times G$ where L = {e, a, b, c, d) is a loop and G = $\{x_1, x_2, x_3\}$ is a groupoid given by the following tables.

| * | e | a | b | c | d |
|---|---|---|---|---|---|
| e | e | a | b | c | d |
| a | a | e | c | d | b |
| b | b | d | a | e | c |
| c | c | b | d | a | e |
| d | d | c | e | b | a |

| + | $x_1$ | $x_2$ | $x_3$ |
|---|---|---|---|
| $x_1$ | $x_1$ | $x_3$ | $x_2$ |
| $x_2$ | $x_2$ | $x_1$ | $x_3$ |
| $x_3$ | $x_3$ | $x_2$ | $x_1$ |

L = {e, a, b, c, d} and G = $\{x_1, x_2, x_3\}$, o(A) = 7 + 18 + 16 + 15 = 56. Now consider the subset $L' = L_1 \cup L_2 \cup L_3$ of A where $L_1 = L_5$ (3), $L_2 = L_7$ (2) and $L_3$ = L. Clearly L' is sub3-loop of A. o(L') = 6 + 8 + 5 = 19. 19 $\not|$ 56.

***Example 5.3.15:*** Let $A = \{A_1 \cup A_2 \cup A_3 \cup A_4, *_1, \ldots, *_4\}$ where

$A_1$ = $\{g \mid g^{12} = 1\}$,
$A_2$ = L = {(e, a, b, c, d) × $Z_2$ {= {0, 1}, group under '+' modulo 2},
$A_3$ = {$Z_3$, semigroup under multiplication modulo 3} and
$A_4$ = {$L_5$ (3) × {$Z_5$, for a, b ∈ $Z_5$, a * b = 2a + 3b (mod 5)} be a 4-glsg.

o(A) = 12 + 10 + 3 + 30 = 55. Consider the proper subset T = $\{T_1 \cup T_2\}$ where

$T_1 = L \times \{0\}$ and
$T_2 = L_5(3) \times \{0\}$.

T is a sub2-loop of A and o (T) = 5 + 6 = 11 and o (T) / o (A).

Further H has only one sub-t- loop so A is a Lagrange loop-loop.

Consider the proper subset $P = P_1 \cup P_2 \cup P_3$ where

$P_1 = \{1, g^4, g^8\}$,
$P_2 = \{e\} \times \{0, 1\}$ and
$P_3 = \{e, g_1\} \times \{0\}$.

P is a sub 3-group of A, o(P) = 3 + 2 + 2 = 7. o(P) $\not|$ o(A) i.e. 7 $\not|$ 55.

The other sub 3-groups of A are $S = S_1 \cup S_2 \cup S_3$ where

$S_1 = \{1, g^2, g^4, g^6, g^8, g^{10}\}$,
$S_2 = \{e\} \times \{0, 1\}$ and
$S_3 = \{e, g_2\} \times \{0\}$.

o(S) = 6 + 2 + 2 = 10. o(S) $\not|$ o(A) as 8 $\not|$ 55.

Thus we see none of the order of the sub-3-group divides the order of A. So A is a Lagrange free sub-3-group so it is evident a N-glsg can be a Lagrange for some structure, Lagrange free for some other structure.

**DEFINITION 5.3.11:** *Let $A = \{A_1 \cup A_2 \cup \ldots \cup A_N; *_1, \ldots, *_N\}$ be a N-glsg of finite order. Let $P = \{P_1 \cup \ldots \cup P_s\}$ be a proper subset of A where P is a sub-s-semigroup of A. If the order of every sub–s-semigroup of A divides the order of A then we call A to be Lagrange sub-s-semigroup. If atleast the order of one of the sub-s-semigroup divides the order of A then we call A a weakly Lagrange sub-s-semigroup. If the order of no sub-s-semigroup divides the order of A then we call A the Lagrange free sub-s-semigroup.*

***Example 5.3.16:*** Let A = {$A_1 \cup A_2 \cup A_3 \cup A_4 \cup A_5$; $*_1, \ldots, *_5$} be a finite 5-glsg where

$A_1$ = {g | $g^5$ = 1},
$A_2$ = {$L_5(3)$, the loop of order 6},
$A_3$ = {$Z_{12}$, the semigroup under multiplication modulo 12},
$A_4$ = {$Z_{15}$, the semigroup under multiplication modulo 15}
and
$A_5$ = {$Z_5$ / for a, b $\in$ $Z_5$, a * b = a + 4 b (mod 5)} × {$Z_2$ = {0 1} semigroup under product modulo 2} = $B_1 \times B_2$.

o(A) = 5 + 6 + 12 + 15 + 10 = 48. Consider the subset G = {$G_1 \cup G_2 \cup G_3$} of A where

$G_1$ = {0, 2, 4, 8, 10},
$G_2$ = {0, 5, 10} and
$G_3$ = {0} × {0 1}.

o (G) = 5 + 3 + 2 = 10. o(G) $\not|$ o (A) i.e. 10 $\not|$ 48. Take S = $S_1 \cup S_2 \cup S_3$ a proper subset of A where

$S_1$ = {0, 2, 4, 8, 10},
$S_2$ = {0, 3, 6, 9, 12} and
$S_3$ = {0} × {0, 1}.

S is a 3-subsemigroup of A. o(S) = 5 + 5 + 2 = 12, o (S) / o (A) i.e. 12/48. Thus A is only a weakly Lagrange sub3-semigroup and not a Lagrange sub3-semigroup.

Now we proceed on to define the Lagrange groupoids-groupoids of the N-glsg.

**DEFINITION 5.3.12:** *Let A = {$A_1 \cup A_2 \cup \ldots \cup A_N$; $*_1, \ldots, *_N$} be a N-glsg of finite order. Suppose A has groupoids-groupoids say G = $G_{i_1} \cup \ldots \cup G_{i_r}$ such that the order of every sub r-groupoid divides the order of A, them we call A to be a*

*Lagrange sub-r-groupoid. If there exists atleast one groupoids-groupoids in A such that o(G) / o(A) then we call A a weakly Lagrange groupoid-groupoid. If A has groupoid-groupoid G but none of their order divide order of G then we call A a Lagrange free groupoid-groupoid.*

Now we illustrate these situations by the following examples.

***Example 5.3.17:*** Let A = $\{A_1 \cup A_2 \cup A_3 \cup A_4 \cup A_5; *_1, \ldots, *_N\}$ be a 5-glsg of finite order where

$A_1$ = $\{g \mid g^3 = 1\}$ cyclic group of order 3,
$A_2$ = $L_5$ (3) a loop of order 6,
$A_3$ = $\{Z_{10}$, a semigroup under multiplication modulo 10$\}$,
$A_4$ = $\{Z_{12}$, such that for a, b $\in Z_{12}$, a * b = 3a + 9b (mod 12)$\}$ and
$A_5$ = $\{Z_6$, such that for a, b $\in Z_6$, a * b = 2a + b (mod 6)$\}$.

Clearly o (A) = 3 + 6 + 10 + 12 + 6 = 37 a prime.

We see A has sub2 groupoids given by L = $(L_1 \cup L_2)$ where

$L_1 = \{0, 3, 6, 9\} \subset Z_{12}$ and
$L_2 = \{0, 2, 4\} \subset Z_6$.

o(L) = 4 + 3 = 7, 7 $\nmid$ 37. Thus A is a Lagrange free sub-r-groupoid as order of A is a prime.

We have the following interesting result namely. If A is a N-glsg of finite order n where n is a prime then A is Lagrange free for all sub-r-structures. We now give an example in which A is a weakly Lagrange groupoid-groupoid.

***Example 5.3.18:*** Let A = $\{A_1 \cup A_2 \cup A_3 \cup A_4 \cup A_5; *_1, *_2, *_3, *_4, *_5\}$ a 5-glsg of finite order where

$A_1$ = $\{g \mid g^5 = 1\}$,
$A_2$ = $L_7$ (3) a loop of order 8,
$A_3$ = $\{Z_{10}$ such that it is a semigroup under multiplication modulo 10$\}$,

$A_4$ = $\{Z_{12}$, for a, b in $Z_{12}$, a * b = a + 5b (mod 12)} and
$A_5$ = $\{Z_{15}$, for a, b in $Z_{15}$, a * b = 2a + 3b (mod 15)};

where $A_4$ and $A_5$ are groupoids. o(A) = 5 + 8 + 10 + 12 + 15 = 50.

Take K = $\{K \cup K_2\}$ a proper subset of A where

$K_1$ = $\{0, 3, 6, 9\} \subset Z_{12}$ and
$K_2$ = $\{0, 5, 10\} \subset Z_{15}$.

K is a sub2-groupoid, o(K) = 4 + 3 = 7, 7 $\nmid$ 50.

Let S = $\{S_1 \cup S_2\}$ where

$$S_1 = \{0, 6\} \text{ and } S_2 = \{0, 5, 10\}$$

subgroupoids of $A_4$ and $A_5$ respectively. Then S is a sub2-groupoid of A.

o (S) = 2 + 3 = 5; o(S) / o (A) i.e. 5/50.

Thus A is only a weakly Lagrange sub-2-groupoid. Now we proceed on to define Cauchy element of A, N-glsg when A is of finite order.

**DEFINITION 5.3.13:** *Let A = $\{A_1 \cup A_2 \cup \ldots \cup A_N, *_1, \ldots, *_N\}$be a N-glsg of finite order suppose for some a $\in A_i$ and $a^t$ = identity (t > 1 and t is the least positive number) and if t / o (A) then A is called the Cauchy element of A.*

It is important to note that A may have elements x of finite order such that, $x^n$ = identity till n $\nmid$ o (A).

So we make the following definition.

**DEFINITION 5.3.14:** *Let A = $\{A_1 \cup A_2 \cup \ldots \cup A_N; *_1, \ldots, *_N\}$ be a N-glsg of finite order. Suppose every element a in A which is of finite order $a^t$ is equal to identity is such that t / o(A) then we call A to be a Cauchy N-glsg.*

*Note:* It may happen that A has elements x which are not finite order i.e. $x^n = 0$, $x^t = x$ or so on.

Now if A has atleast one Cauchy element then we call A to be a weakly Cauchy N-glsg. If A has no Cauchy element then we call A to be a Cauchy free N-glsg. (In all these case we assume A to be of finite order).

***Example 5.3.19:*** Let A = $\{A_1 \cup A_2 \cup A_3 \cup A_4 \cup A_5; *_1, *_2, *_3, *_4, *_5\}$ where

$A_1$ = $\{g \mid g^{12} = 1\}$ cyclic group of order 12,
$A_2$ = $L_5(3)$ is the loop of order 6,
$A_3$ = $L_7(4)$ is the loop of order 8,
$A_4$ = $\{Z_{12}$, the semigroup under multiplication modulo 12$\}$
and
$A_5$ = $\{Z_{10}$, for a, b $\in Z_{10}$ a * b = 3a + 7b (mod 10)$\}$.

o(A) = 12 + 6 + 8 + 12 + 10 = 48.

It is easily verified that a is a Cauchy N-glsg. Now we give an example of a weakly Cauchy glsg.

***Example 5.3.20:*** Let A = $\{A_1 \cup A_2 \cup A_3 \cup A_4 \cup A_5; *_1, *_2, \ldots, *_5\}$ where

$A_1$ = $\{g \mid g^7 = 1\}$ a cyclic group of order 7,
$A_2$ = $\{L_5(3)\}$ loop of order 6,
$A_3$ = $L_7(4)$, loop of order 7,
$A_4$ = $\{Z_{10}$, the semigroup under multiplication modulo 10$\}$
and
$A_5$ = $\{Z_{14} \mid$ for a, b $\in Z_{14}$, a * b = a + 6b (mod 14)$\}$.

o (A) = 7 + 6 + 8 + 10 + 14 = 45. Every element in $A_1$ is of order 7, $7 \nmid$ . 45.

Every element in $A_2$ and $A_3$ are of order 2, $2 \nmid 45$. Clearly A is Cauchy free.

***Example 5.3.21:*** Let A = $\{A_1 \cup A_2 \cup A_3 \cup A_4; *_1, *_2, *_3, *_4\}$ be a N-glsg where

$A_1$ = $\{g \mid g^9 = 1\}$,
$A_2$ = $L_5$ (3) a loop of order 6,
$A_3$ = $\{Z_9 \mid Z_9$ is a semigroup under multiplication modulo 9$\}$ and
$A_4$ = $\{Z_8$ / for a, b $\in Z_8$, a * b = 2a + 6 b (mod 8)$\}$ ,

o(A) = 9 + 6 + 9 + 8 = 32.

Now all elements in $L_5(3)$ is of order 2 and 2/ o(A) i.e. 2/32. No element of $A_1$ divides 32. Thus A is only a weakly Lagrange N-glsg.

Next we proceed in to analyze about p-Sylow sub N-glsg of a finite N-glsg.

**DEFINITION 5.3.15:** *Let A = $\{A_1 \cup A_2 \cup \ldots \cup A_N; *_1, \ldots, *_N\}$ be a N-glsg of finite order (say n).*

*Suppose p is a prime such that $p^\alpha \nmid o(A)$ but $p^{\alpha+1} \nmid o(A)$ and if A has a sub N-glsg P of order $p^\alpha$ then we call P a p-Sylow sub N-glsg of A.*

*If for every prime p such that $p^\alpha$ / o(A) and $p^{\alpha+1} \nmid o(A)$ we have a p-Sylow sub N-glsg in A then we call A itself a Sylow N-glsg.*

*If A has atleast one p Sylow N-glsg then we call A a weakly Sylow N-glsg. If A has no p-Sylow sub N glsg then we call A to be a Sylow free N glsg.*

Now it is interesting to note that for a given prime p such that $p^\alpha$ / o(A) and $p^{\alpha+1} \nmid$ o(A).

One may have a sub N-glsg of order $p^{\alpha+t}$ ($t \geq 1$). Such situations do not in general occur only in groups but in all other N-structures such occurrences are very common. It is left for the reader to illustrate these situation by examples.

**DEFINITION 5.3.16:** *Let A = $\{A_1 \cup A_2 \cup \ldots \cup A_N, *_1, \ldots, *_N\}$ be a N-glsg. A is said to be a Smarandache N-glsg (S-N-glsg) if* *some $A_i$ are S-loops, some of $A_j$ are S-semigroups and some of $A_K$ are S-groupoids.*

**DEFINITION 5.3.17:** *Let $A = \{A_1 \cup A_2 \cup \ldots \cup A_N, *_1, *_2, \ldots, *_N\}$ be a N-glsg. A proper subset $P = \{P_1 \cup P_2 \cup \ldots \cup P_N, *_1, \ldots, *_N\}$ of A is said to be a Smarandache sub N-glsg if P itself is a S-N-glsg under the operations of A.*

We just illustrate by the following example:

***Example 5.3.22:*** Let $A = \{A_1 \cup A_2 \cup A_3 \cup A_4, *_1, *_2, *_3, *_4\}$ be a 4 glsg where

$A_1$ = $S_4$,
$A_2$ = $L_7(3)$; S-loop of order 8,
$A_3$ = S(3); S-semigroup of order 27 and
$A_4$ = $\{Z_4 \times Z_{10}$; where $Z_4$ is a groupoid such that $a * b = 2a + b (\text{mod } 4)$ for $a, b \in Z_4$ and $Z_{10}$ is a semigroup under multiplication modulo 10$\}$.

A is S-4-glsg.

A study of S-N- glsg can be carried out systematically as it has been done for N-glsg.

**Chapter Six**

# PROBLEMS

In this chapter we give problems for the reader to solve. Most of the problems are simple, only a few of them are difficult. Also, in this book problems have been given in each and every chapter as discussions and analysis are carried out. A researcher can make use of them and make it a point to solve these problems, which will certainly enable her/him to understand the topic better.

1. Define sub N-semigroup and illustrate it by an example.

2. Find all sub N-semigroups of $S = \{S_1 \cup S_2 \cup S_3 \cup S_4; *_1, *_2, *_3, *_4\}$ where

   $S_1$ = $S(3)$ ,
   $S_2$ = $\{Z_{12}, *_2$ is a semigroup under multiplication modulo 12$\}$,
   $S_3$ = $\left\{\begin{pmatrix} a & b \\ c & d \end{pmatrix} \mid a, b, c, d \in Z_4\right\}$, semigroup under matrix multiplication and
   $S_4$ = $\{(x, y) / x, y \in Z_6\}$ semigroup under component-wise multiplication.

3. Define order of a N-semigroup. Give an example of a N-semigroup of infinite order.

4. What is the order of the 4-semigroup given in Problem 2?

5. Define Lagrange sub N-semigroup. Illustrate it by an example.

6. Define Lagrange-N-semigroup. Give an example of a Lagrange 4-semigroup.

7. Define weakly Lagrange N-semigroup. Find a 7-semigroup which is a weakly Lagrange 7-semigroup.

8. Define N-ideals of a N-semigroup. Find N-ideals of the 3-semigroup, $\{S = S_1 \cup S_2 \cup S_3, *_1, *_2, *_3\}$ where

$S_1$ = S(3),
$S_2$ = $\{Z_{15}$, semigroup under multiplication modulo 15$\}$ and

$$S_3 = \left\{ \begin{pmatrix} a & b \\ c & d \end{pmatrix} \middle| a, b, c, d \in Z_6 \right\}.$$

9. Can Cayley's theorem be extended for N-semigroup if one defines a symmetric N-semigroup as $S = \{S(N_1) \cup \ldots \cup S(N_N)\}$ where $S(N_i)$ is the symmetric semigroup of all mappings of the set $(1, 2, \ldots, N_i)$ to itself?

10. Give an example of a S-N-semigroup.

11. Define sub S-N-semigroup and illustrate it with examples.

12. Find all the sub S-N-semigroups of $S = \{S_1 \cup S_2 \cup S_3 \cup S_4, *_1, *_2, *_3, *_4\}$ where $S_1 = A_4$, $S_2 = S(3)$, $S_3 = D_{2.7}$ and $S_4 = \{Z_{24}$, semigroup under multiplication modulo 24$\}$.

13. Define S-Lagrange sub N-semigroup of a S-N-semigroup.

14. Does the example given in Problem 12, have S-Lagrange sub N-semigroups?

15. Define S-Lagrange N-semigroup and illustrate it with examples.

16. Can a necessary and sufficient condition be given for a finite S-N semigroup to be Lagrange S-N-semigroup?

17. Obtain a modified form of Cayley's theorem for S-N-semigroups.

18. Define S-N-ideal in a N-semigroup.

19. Find S-N-ideals of the S-N-semigroup $S = \{S_1 \cup S_2 \cup S_3 \cup S_4, *_1, *_2, *_3, *_4\}$ where $S_1 = S(3)$, $S_2 = D_{2.6}$, $S_3 = A_4$ and $S_4 = Z_{12}$, semigroup under multiplication modulo 12.

20. Give an example to show that every sub S-N-semigroup in general is not an S-N-ideal.

21. Define S-homomorphism of S-N-groups. Illustrate it by an example.

22. Is the principal isotope of a N-loop a N-loop?

23. Does the 4-loop given by $L = \{L_1 \cup L_2 \cup L_3 \cup L_4, *_1, *_2, *_3, *_4\}$ where $L_1 = A_4$, $L_2 = S_3$, $L_3 = L_7(5)$ and $L_4 = L_5(2)$ have a principal isotope? Justify your claim.

24. Give example of S-commutative 7-loop?

25. Give an example of a S-weakly commutative 6-loop?

26. Does their exist a S-N-loop which is not S-weakly commutative? Prove your claim.

27. Give an example of a S-N-loop which has no S-commutator S-sub N-loop.

28. Give an example of a S-N-loop which has a non-trivial S-commutator sub N-loop.

29. Give an example of a S-5-loop and find its S-associator.

30. Give an example of S-3-loop which has no pseudo commutative pair.

31. Give an example of S-5-loop which has at least one S-pseudo commutative pair.

32. Give a non-trivial example of a S-pseudo commutative N-loop.

33. Give an example of a S-N-loop L which has more than one S-pseudo commutator subloop.

34. Give an example of a S-pseudo associative sub-N-loop (N = 5).

35. Can a S-N-loop L have more than one $S(A\,L_P^N)$?

36. Give an example of a S-6-loop which has a

    a. Non-trivial S-Moufang center.
    b. No S-Moufang center.

37. Can a necessary and sufficient condition be found on the N-loop L so that L is a Lagrange N-loop?

38. Define middle, right and left nucleus of a N-loop L. Illustrate them by examples. Can one always define middle, right and left nucleus for any N-loop L?

39. Derive class equations in case of N-groups.

40. Can a class equation be defined in case of N-group semigroup? Justify your claim.

41. Define normalizer of an element in case of a N-group. What will be the algebraic structure of a normalizer of a N-group?

42. Define normalizer and centralizer of the N-group $G = G_1 \cup G_2 \cup G_3 \cup G_4$ where

$G_1$ = $S_5$, the symmetric group of degree 5,
$G_2$ = $D_{2.7} = \{a, b / a^2 = b^7 = 1, bab = a\}$ the dihedral group of order 14,
$G_3$ = $\langle g / g^8 = 1\rangle$ and
$G_4$ = $A_6$ the alternating group of degree 6.

What is the structure of them?

43. Give an example of S-N-loop which has more than one S-Moufang center.

44. Can you characterize 5-N-group to be Lagrange N-group?

45. Obtain some interesting results on N-gslg structures.

46. Define S-N-glsg and illustrate it with example.

47. Define S-sub N-glsg.

48. Define

a. Lagrange S-sub N-glsg
b. p-Sylow S-sub N-glsg
c. Cauchy element in S-N-glsg.

49. Give an example of S-N-glsg which is a

a. Lagrange S-N-glsg
b. Cauchy S-N-glsg and
c. Sylow S-N- glsg.

50. Obtain a necessary and sufficient condition on the group and semigroup so that a S-group semigroup is

  a. Lagrange.
  b. Sylow.
  c. Cauchy.

51. Let L = $\{L_1 \cup L_2 \cup L_3, *_1, *_2, *_3\}$ where $L_1 = \{Z_{10}$, with $*_1$ defined on $Z_{10}$ as a $*_1$ b = 2a + 3b (mod 10) where a, b $\in Z_{10}\}$, $L_1$ is a groupoid,
$L_2$ is the loop given by the following table:

| $*_2$ | e | $a_1$ | $a_2$ | $a_3$ | $a_4$ | $a_5$ |
|---|---|---|---|---|---|---|
| e | e | $a_1$ | $a_2$ | $a_3$ | $a_4$ | $a_5$ |
| $a_1$ | $a_1$ | e | $a_3$ | $a_5$ | $a_2$ | $a_4$ |
| $a_2$ | $a_2$ | $a_5$ | e | $a_4$ | $a_1$ | $a_3$ |
| $a_3$ | $a_3$ | $a_1$ | $a_1$ | e | $a_5$ | $a_2$ |
| $a_4$ | $a_4$ | $a_3$ | $a_5$ | $a_2$ | e | $a_1$ |
| $a_5$ | $a_5$ | $a_2$ | $a_4$ | $a_1$ | $a_3$ | e |

$L_3$ is the loop given by the following table:

| $x_3$ | e | a | b | c | d |
|---|---|---|---|---|---|
| e | e | a | b | c | d |
| a | a | e | c | d | b |
| b | b | d | a | e | c |
| c | c | b | d | a | e |
| d | d | c | e | b | a |

Find all its substructures.

# REFERENCES

1. Albert, A.A., *Non-associative algebra I, II*, Ann. Math. (2), **43**, 685-707, (1942).

2. Birkhoff, G. and Bartee, T.C. *Modern Applied Algebra,* Mc-Graw Hill, New York, (1970).

3. Bruck, R. H., *A survey of binary systems*, Springer-Verlag, (1958).

4. Bruck, R.H, *Some theorems on Moufang loops*, Math. Z., **73**, 59-78 (1960).

5. Castillo J., *The Smarandache Semigroup*, International Conference on Combinatorial Methods in Mathematics, II Meeting of the project 'Algebra, Geometria e Combinatoria', Faculdade de Ciencias da Universidade do Porto, Portugal, 9-11 July 1998.

6. Chang Quan, Zhang, *Inner commutative rings*, Sictiuan Dascue Xuebao (Special issue), **26**, 95-97 (1989).

7. Chein, Orin and Goodaire, Edgar G., *Loops whose loop rings in characteristic 2 are alternative*, Comm. Algebra, **18**, 659-668 (1990).

8. Chein, Orin, and Goodaire, Edgar G., *Moufang loops with unique identity commutator (associator, square)*, J. Algebra, **130**, 369-384 (1990).

9. Chein, Orin, and Pflugfelder, H.O., *The smallest Moufang loop*, Arch. Math., **22**, 573-576 (1971).

10. Chein.O, Pflugfelder.H.O and Smith.J.D.H, (eds), *Quasigroups and loops: Theory and applications*, Sigma Series in Pure Maths, Vol. 8, Heldermann Verlag, (1990).

11. Chein, Orin, Kinyon, Michael. K., Rajah, Andrew and Vojlechovsky, Peter, *Loops and the Lagrange property*, (2002). http://lanl.arxiv.org/pdf/math.GR/0205141

12. Goodaire, E.G., and Parmenter, M.M., *Semisimplicity of alternative loop rings*, Acta. Math. Hung, **50**. 241-247 (1987).

13. Hall, Marshall, *Theory of Groups*. The Macmillan Company, New York, (1961).

14. Hashiguchi. K, Ichihara, S. and Jimbo, S., *Formal languages over binoids*, J. Autom Lang Comb, **5**, 219-234 (2000).

15. Herstein., I.N., *Topics in Algebra*, Wiley Eastern Limited, (1975).

16. Ivan, Nivan and Zukerman. H. S., *Introduction to number theory*, Wiley Eastern Limited, (1984).

17. Lang, S., *Algebra*, Addison Wesley, (1967).

18. Maggu, P.L., *On introduction of Bigroup concept with its applications in industry*, Pure App. Math Sci., **39**, 171-173 (1994).

19. Maggu, P.L., and Rajeev Kumar, *On sub-bigroup and its applications*, Pure Appl. Math Sci., **43**, 85-88 (1996).

20. Michael.K.Kinyon and Phillips.J.D, *Commutants of Bol loops of odd order*, (2002). http://lanl.arxiv.org/pdf/math.GR/0207119

21. Michael.K.Kinyon and Oliver Jones, *Loops and semidirect products*, (2000). http://lanl.arxiv.org/pdf/math.GR/9907085 (To appear in Communications in Algebra)

22. Raul, Padilla, *Smarandache Algebraic Structures*, Smarandache Notions Journal, **9**, 36-38 (1998).

23. Singh, S.V., *On a new class of loops and loop rings*, Ph.D. thesis IIT (Madras), guided by Vasantha. W.B., (1994).

24. a. Smarandache, Florentin, *Special Algebraic Structures*, in Collected Papers, Abaddaba, Oradea, **3**, 78-81 (2000).
b. Smarandache F., Multi structures and Multi spaces, (1969) www.gallup.unm.edu/~smarandache/transdis.txt

25. Solarin, A.R.T., and Sharma B.L., *On the identities of Bol-Moufang type*, Kyungpook Math. J., **28**, 51-62 (1988).

26. Tim Hsu, *Class 2 Moufang loops small Frattini Moufang loops and code loops*, (1996). http://lanl.arxiv.org/pdf/math.GR/9611214

27. Vasantha Kandasamy, W. B., *Fuzzy subloops of some special loops*, Proc. 26th Iranian Math. Conf., 33-37 (1995).

28. Vasantha Kandasamy, W. B., *On ordered groupoids and groupoid rings*, J. Math. Comp. Sci., **9**, 145-147 (1996).

29. Vasantha Kandasamy, W. B. and Meiyappan, D., *Bigroup and Fuzzy bigroup*, Bol. Soc. Paran Mat, **18**, 59-63 (1998).

30. Vasantha Kandasamy, W. B., *On Quasi loops*, Octogon, **6**, 63-65 (1998).

31. Vasantha Kandasamy, W. B., *On a new class of Jordan loops and their loop rings*, J. Bihar Math. Soc., **19**, 71-75 (1999).

32. Vasantha Kandasamy, W. B. and Singh S.V., *Loops and their applications to proper edge colouring of the graph $K_{2n}$*, Algebra and its applications, edited by Tariq et al., Narosa Pub., 273-284 (2001).

33. Vasantha Kandasamy, W. B., *Biloops*, U. Sci. Phy. Sci., **14**, 127-130 (2002).

34. Vasantha Kandasamy, W. B., *Groupoids and Smarandache groupoids*, American Research Press, Rehoboth, (2002). http://www.gallup.unm.edu/~smarandache/Vasantha-Book2.pdf

35. Vasantha Kandasamy, W. B., *On Smarandache Cosets*, (2002). http://www.gallup.unm.edu/~smaranandache/pseudoideals.pdf

36. Vasantha Kandasamy, W. B., *Smarandache groupoids*, (2002). http://www.gallup.unm.edu/~smarandache/Groupoids.pdf

37. Vasantha Kandasamy, W. B., *Smarandache loops*, Smarandache Notions Journal, **13**, 252-258 (2002). http://www.gallup.unm.edu/~smarandache/Loops.pdf

38. Vasantha Kandasamy, W. B., *Smarandache Loops*, American Research Press, Rehoboth, NM, (2002). http://www.gallup.unm.edu/~smarandache/Vasantha-Book4.pdf

39. Vasantha Kandasamy, W. B., *Bialgebraic Structures and Smarandache Bialgebraic Structures*, American Research Press, Rehoboth, NM, (2002). http://www.gallup.unm.edu/~smarandache/NearRings.pdf

40. Vasantha Kandasamy, W. B., *Smarandache Semigroups*, American Research Press, Rehoboth, NM, (2002). http://www.gallup.unm.edu/~smarandache/Vasantha-Book1.pdf

# INDEX

**A**

**B**

**C**

**D**

**F**

**G**

**I**

**K**

**L**

**M**

**N**

**P**

**Q**

**R**

**S**

## SMARANDACHE NOTIONS

**T**

**U**

**W**

**Z**

## ABOUT THE AUTHORS

**Dr.W.B.Vasantha Kandasamy** is an Associate Professor in the Department of Mathematics, Indian Institute of Technology Madras, Chennai, where she lives with her husband Dr.K.Kandasamy and daughters Meena and Kama. Her current interests include Smarandache algebraic structures, fuzzy theory, coding/ communication theory. In the past decade she has guided 11 Ph.D. scholars in the different fields of non-associative algebras, algebraic coding theory, transportation theory, fuzzy groups, and applications of fuzzy theory of the problems faced in chemical industries and cement industries. Currently, four Ph.D. scholars are working under her guidance.

She has to her credit 612 research papers of which 209 are individually authored. Apart from this, she and her students have presented around 329 papers in national and international conferences. She teaches both undergraduate and post-graduate students and has guided over 45 M.Sc. and M.Tech. projects. She has worked in collaboration projects with the Indian Space Research Organization and with the Tamil Nadu State AIDS Control Society. This is her 24th book.

She can be contacted at vasantha@iitm.ac.in
You can visit her work on the web at: http://mat.iitm.ac.in/~wbv

---

**Dr.Florentin Smarandache** is an Associate Professor of Mathematics at the University of New Mexico in USA. He published over 75 books and 100 articles and notes in mathematics, physics, philosophy, psychology, literature, rebus. In mathematics his research is in number theory, non-Euclidean geometry, synthetic geometry, algebraic structures, statistics, neutrosophic logic and set (generalizations of fuzzy logic and set respectively), neutrosophic probability (generalization of classical and imprecise probability). Also, small contributions to nuclear and particle physics, information fusion, neutrosophy (a generalization of dialectics), law of sensations and stimuli, etc.

He can be contacted at smarand@unm.edu